\documentclass[fleqn,11pt]{article}
\usepackage{amssymb,amsmath,amsfonts}
\usepackage{color}
\usepackage{graphicx}
\usepackage{latexsym}
\usepackage{amsmath}
\usepackage{amssymb}
\usepackage{graphics}
\usepackage[dvips]{epsfig}
\usepackage[margin=1in]{geometry}

\numberwithin{equation}{subsection}
\newtheorem{thm}{Theorem}[section]

\newtheorem{lem}[thm]{Lemma}
\newtheorem{prop}[thm]{Proposition}
\newtheorem{rmk}[thm]{Remark}
\newtheorem{defi}[thm]{Definition}

\newcommand{\ea}{\epsilon}
\newcommand{\ta}{\theta}
\newcommand{\ka}{\kappa}
\newcommand{\da}{\delta}
\newcommand{\la}{\lambda}
\renewcommand{\aa}{\alpha}

\newcommand{\pl}{\partial}
\newcommand{\sa}{\sigma}
\newcommand{\oa}{\omega}
\newcommand{\ga}{\gamma}

\newcommand{\iy}{\infty}

\newcommand{\lt}{\left}
\newcommand{\rt}{\right}
\newcommand{\be}{\begin{equation}}
\newcommand{\bs}{\begin{split}}
\newcommand{\es}{\end{split}}
\newcommand{\ee}{\end{equation}}
\newcommand{\bee}{\begin{equation*}}
\newcommand{\eee}{\end{equation*}}

\newcommand{\ef}{\eqref}

\begin{document}
\begin{center}
\large{ \bf On Nonlinear
  Asymptotic Stability of the Lane-Emden Solutions for the Viscous
  Gaseous Star Problem }
\end{center}
\begin{center} Tao Luo,  Zhouping Xin,  Huihui Zeng 
\end{center}

\begin{abstract} This paper proves the nonlinear asymptotic stability of the Lane-Emden solutions for
spherically symmetric motions of viscous gaseous stars if the adiabatic constant $\gamma$ lies in the stability
range $(4/3, 2)$. It is shown that for small perturbations of a Lane-Emden solution with same
mass, there exists a unique global (in time)  strong solution to the vacuum free boundary problem of the compressible
Navier-Stokes-Poisson system with spherical symmetry for  viscous stars, and the solution captures the precise
physical behavior that the sound speed is $C^{{1}/{2}}$-H$\ddot{\rm o}$lder continuous  across the vacuum boundary provided that $\gamma$ lies in $(4/3, 2)$. The key is to establish the global-in-time regularity uniformly up to the vacuum boundary, which ensures  the large time asymptotic uniform convergence of the evolving vacuum boundary, density and velocity to those of  the Lane-Emden solution with detailed convergence rates, and detailed large time behaviors of solutions near the vacuum boundary. In particular, it is shown that every spherical surface moving with the fluid  converges to the sphere  enclosing the same mass inside the domain of the Lane-Emden solution with a uniform convergence rate and  the large time asymptotic states for the vacuum free boundary problem \ef{103}  are  determined by the initial mass distribution and the total mass. To overcome the difficulty caused by the degeneracy and
singular behavior near  the vacuum free boundary and coordinates singularity at the symmetry center,  the main ingredients of the analysis consist of  combinations of some new weighted
nonlinear functionals (involving both lower-order and higher-order derivatives) and space-time weighted energy estimates.
The constructions of these weighted nonlinear functionals  and space-time weights depend crucially on the structures of the Lane-Emden
solution, the balance of pressure and gravitation,  and the dissipation. Finally, the uniform boundedness of the acceleration of the
vacuum boundary is also proved.
\end{abstract}
\tableofcontents

\section{Introduction}

\subsection{Problem}
In the fundamental hydrodynamical setting (cf. \cite{ch}),  the evolving boundary of a viscous gaseous star (the interface of fluids and vacuum states)  can be
modeled by the following free boundary problem of  the compressible Navier-Stokes-Poisson equations:
\begin{equation}\label{0.1}\begin{split}
&   \rho_t  + {\rm div}(\rho {\bf u}) = 0 &  {\rm in}& \ \ \Omega(t), \\
 &    (\rho {\bf u})_t  + {\rm div}(\rho {\bf u}\otimes {\bf u})+{\rm div}\mathfrak{S }= - \rho \nabla_{\bf x} \Psi  & {\rm in}& \ \ \Omega(t),\\
 &\rho>0 &{\rm in }  & \ \ \Omega(t),\\
 &\rho=0 \ \ {\rm and} \ \  \mathfrak{S}{\bf n}={\bf 0}   &    {\rm on}& \  \ \Gamma(t):=\pl \Omega(t),\\
 &    \mathcal{V}(\Gamma(t))={\bf u}\cdot {\bf n}, & &\\
&(\rho,{\bf u})=(\rho_0, {\bf u}_0) & {\rm on} & \ \  \Omega:= \Omega(0).
 \end{split} \end{equation}
 Here $({\bf x},t)\in \mathbb{R}^3\times [0,\iy)$,  $\rho $, ${\bf u} $, $\mathfrak{S}$ and $\Psi$ denote, respectively, the space and time variable, density, velocity, stress tensor and gravitational potential; $\Omega(t)\subset \mathbb{R}^3$, $\Gamma(t)$, $\mathcal{V}(\Gamma(t))$ and ${\bf n}$ represent, respectively, the changing volume occupied by a fluid at time $t$, moving interface of fluids and vacuum states, normal velocity of $\Gamma(t)$ and exterior unit normal vector to $\Gamma(t)$. The gravitational potential is described by
 $$\Psi({\bf x}, t)=-G\int_{\Omega(t)} \frac{\rho({\bf y}, t)}{|{\bf x}-{\bf y}|}d{\bf y} \ \ {\rm satisfying} \  \  \Delta \Psi=4\pi G \rho \  \ {\rm in} \ \  \Omega(t)
$$
with  the gravitational constant $G$ taken to be unity for convenience.
The  stress tensor  is given by
$$
\mathfrak{S}=pI_3-\lambda_1 \left(\nabla {\bf u}+\nabla {\bf u}^t-\frac{2}{3}({\rm div} {\bf u}) I_3\right)-\lambda_2({\rm div}  {\bf u})I_3, $$
 where $I_3$ is the $3\times 3$ identical matrix, $p$ is the pressure of the gas, $\lambda_1>0$ is the shear viscosity, $\lambda_2>0$ is the bulk viscosity, and $\nabla {\bf u}^t$ denotes the transpose of $\nabla {\bf u}$.
 We consider the polytropic gases for which the equation of state is given by
$$
 p=p(\rho)=K\rho^{\gamma}, $$
 where $K>0$ is a constant set to be unity for convenience, $\gamma>1$ is the adiabatic exponent.

For a non-rotating gaseous star,  it is important to consider  spherically symmetric motions since the stable  equilibrium configurations, which minimize the energy among all possible configurations (cf. \cite{liebyau}), are spherically symmetric, called Lane-Emden solutions.  In this work, we are concerned with the three-dimensional spherically symmetric solutions to the free boundary problem \eqref{0.1} and its nonlinear asymptotic stability toward the Lane-Emden solutions. The aim is to prove the global-in-time regularity uniformly up to the vacuum boundary of  solutions when $4/3<\gamma<2$  (the stable index) capturing an interesting behavior called  the physical vacuum  (cf. \cite{10,10',jangmas,jm,tpliudamping,13,ya}) which states that the sound speed $c=\sqrt{p'(\rho)}$ is $C^{ {1}/{2}}$-H$\ddot{\rm o}$lder continuous near the vacuum boundary, as long as the initial datum is a suitably small perturbation  of the Lane-Emden solution  with the same total mass.
Furthermore, we establish  the large time asymptotic convergence of the global strong solution, in particular, the convergence of  the vacuum boundary and  the  density,  to the  the Lane-Emden solutions with the detailed convergence rate as the time goes to infinity.

In the spherically symmetric setting, that is, $\Omega(t)$ is a ball with the  changing radius $R(t)$,
$$
 \rho({\bf x}, t) = \rho(r, t) \  \ {\rm and} \ \  {\bf u}({\bf x}, t) =  u(r, t) {\bf x} /r  \ \ {\rm with} \ \   r=|{\bf x}| \in \lt(0, R(t)\rt); $$
system  \eqref{0.1} can then be rewritten as
\begin{subequations}\label{103}
\begin{align}
&  (r^2\rho)_t+ (r^2\rho u)_r=0  & {\rm in } & \  \  \lt(0, \  R(t)\rt) , \label{103a}\\
&\rho( u_t +u u_r)+  p_r+ {4\pi\rho}r^{-2}\int_0^r\rho(s,t) s^2ds=\mu \lt(\frac{(r^2 u)_r}{r^2} \rt)_r & {\rm in } & \  \ \lt(0, \  R(t)\rt), \label{103b}\\
& \rho>0  & {\rm in } & \  \  \lt[0, \  R(t)\rt), \label{103c}\\
& \rho=0  \ \ {\rm and} \  \  \frac{4}{3}\la_1\lt( u_r-\frac{u}{r}\rt) + \la_2 \lt(  u_r+2\frac{u}{r}\rt)=0  & {\rm for} &  \   \ r=R(t), \label{103d}\\
& \dot R(t)=u(R(t), t) \ \ {\rm with} \ \ R(0)=R_0,  \ \    u(0,t)=0, &    & \label{103e}\\
&  (\rho, u) = (\rho_0, u_0)   & {\rm on } & \ \ (0, \ R_0), \label{103f}
\end{align}
\end{subequations}
where $\mu=4\lambda_1/3+\lambda_2>0$ is the viscosity constant. \ef{103c} and \ef{103d} state that $r=R(t)$ is the vacuum free boundary at which the normal stress $\mathfrak{S}{\bf n}=0$ reduces to
$$p-\frac{4}{3}\la_1\lt(  u_r-\frac{u}{r}\rt) - \la_2 \lt(  u_r+2\frac{u}{r}\rt)=0 \ \  {\rm for} \ \ r=R(t), \ \ t\ge 0;$$
\ef{103e} describes that  the free boundary issues from $r=R_0$ and moves with the fluid velocity, and
the center of the symmetry does not move. The initial domain is taken to be a  ball  $\{0\le r\le R_0\}$, and the initial density is assumed to satisfy
the following  condition:
\be\label{156}
\rho_0(r)>0 \ \ {\rm for} \ \  0\le r<R_0 , \ \ \rho_0(R_0)=0 \ \ {\rm and}  \ \
  -\iy<   \lt(\rho_0^{\ga-1}\rt)_r <0 \  \ {\rm at} \  \ r=R_0;
   \ee
so
\be\label{physicalvacuum}\rho_0^{\ga-1}(r) \sim R_0-r \ {\rm as~} r {\rm~ close~ to~} R_0, \ee
that is, the initial sound speed is $C^{\frac{1}{2}}$-H${\rm \ddot{o}}$lder continuous across the vacuum boundary. The unknowns here are $\rho,\ u$ and $R(t)$.

 The requirement  \eqref{156} for the initial density  near the vacuum boundary is motivated by that of the Lane-Emden solution, $\bar\rho$, (cf. \cite{ch,linss})   which solves
\be\label{le} \pl_r(\bar \rho^{\gamma})+ {4\pi}r^{-2} \bar{\rho} \int_0^r \bar \rho(s) s^2ds=0.
\ee
The solutions to \eqref{le} can be characterized by the values of $\gamma$
(cf. \cite{linss}) for given finite total mass $M > 0$, if $\gamma \in (6/5, 2)$, there exists at least one compactly supported solution.  For $\gamma \in (4/3, 2)$, every solution is compactly supported and unique.
If $\gamma= 6/5$, the unique solution admits an explicit expression, and it has infinite support.
On the other hand, for $\gamma\in  (1, 6/5)$, there are no solutions with finite total mass.  For $\gamma>6/5$, let $\bar R$ be the radius of the stationary star giving by the  Lane-Emden solution, then it holds (cf. \cite{linss,makino})
\be\label{pvforle} \bar\rho^{\ga-1}(r) \sim \bar R-r \ {\rm as~} r {\rm~ close~ to~} \bar R.\ee

\subsection{Motivations and goals}
The problem of nonlinear asymptotic stability of Lane-Emden solutions  is of fundamental importance in both astrophysics and the theory of nonlinear PDEs. It is believed by astrophysicists that Lane-Emden solutions are stable for $ {4}/{3}<\gamma<2$ (cf. \cite{ch,tokusky}) since they minimize the total energy among all the possible configurations. The main aim of this paper is to
justify rigorously  the precise  sense of this stability. In fact, we prove for the viscous gaseous star  with $ {4}/{3}<\gamma<2$ , the Lane-Emden solution is
strongly stable in the sense that it is asymptotically nonlinear stable. The first step for this purpose is to prove the global
existence of strong solutions. However,  due to the high degeneracy of system \eqref{103} caused by the behavior \eqref{156} near the vacuum boundary, it is a very challenging problem even for the local-in-time existence theory.   Indeed, the local-in-time well-posedness of smooth solutions to   vacuum free boundary problems with the behavior that the sound speed is $C^{ {1}/{2}}$-H$\ddot{\rm o}$lder continuous  across vacuum boundaries was only established recently for compressible inviscid flows (cf. \cite{10, 10', jangmas, jm}) (see also
 \cite{LXZ} for a local-in-time well-posedness theory in a new functional space for the three-dimensional compressible Euler-Poisson equations in  spherically symmetric motions). For the vacuum free boundary problem \eqref{103} of the compressible Navier-Stokes-Poisson equations featuring the behavior \eqref{156} near the vacuum boundary,  a local-in-time well-posedness theory of strong solutions was established    in \cite{jangnsp}. In order to obtain the nonlinear asymptotic stability
of Lane-Emden solutions, it turns out that suitable  estimates for higher order  derivatives uniformly up to the vacuum boundary are necessary. Indeed, this turns out to be  essential to prove the convergence of the evolving vacuum boundary and the uniform convergence of the density to those of   Lane-Emden solutions, in addition to the uniform convergence of the velocity.  We show the global-in-time regularity of solutions  when $4/3<\gamma<2$  capturing the behavior \eqref{physicalvacuum} (or \eqref{pvforle}) when the initial data are small perturbations of and have the same total mass as the stationary solution, $\bar\rho$, given by \eqref{le}.
 It should be remarked that   the regularity estimates near boundaries are notoriously
 difficult. This is particularly so for the vacuum boundary problem \ef{103} due to the high degeneracy caused by the singular behavior of \eqref{156} near  vacuum states.

Our nonlinear asymptotic stability results can be stated more precisely as follows.  Suppose that the initial datum $(\rho_0, u_0, R_0)$  is a small perturbation of the Lane-Emden solution $(\bar \rho, 0, \bar R)$ in a suitable sense (see Theorem \ref{mainthm1}) and has the same total mass,
$$\int_0^{R_0}r^2 \rho_0(r)dr= \int_0^{\bar R} r^2 \bar\rho(r)dr,$$
then there is a unique global-in-time strong solution
$(\rho, u, R(t))$ $(0\le t<+\infty)$ to \ef{103} which is regular uniformly  up to the vacuum boundary $r=R(t)$. Moreover, let $r(x, t)$ be the radius of the ball inside  $B_{R(t)}({\bf 0})$ satisfying:
 \be\label{ma1}   r_t(x, t)=u(r(x, t), t)  \ \ {\rm and} \ \    r(x, 0)=r_0(x) \ \  {\rm for} \ \  0\le x\le \bar R, \ee
 \be\label{ma2} \int_0^{r_0(x)}s^2\rho_0(s)ds=\int_0^xs^2\bar\rho(s)ds \  \  {\rm for} \ 0\le x\le \bar R.\ee
Then
\be\label{ma4} \lim_{t\to \infty} \|\lt(r(x, t)-x, \ \rho(r(x, t), t)-\bar\rho(x), \ u(r(x, t), t) \rt)\|_{L^{\infty}_x ([0, \bar R])}= 0 \ee
with some detailed convergence rates.
Notice that \ef{ma1} means that  the  sphere $r=r(x, t)$  with the initial position $r=r_0(x)$ is  moving with the fluid  and \ef{ma2} means that the initial mass inside the ball $B_{r_0(x)}({\bf 0})$ is the same as that of the Lane-Emden solution inside the ball $B_{x}({\bf 0})$ for $0\le x\le \bar R$.  It follows from the conservation of mass that the mass inside the ball $B_{r(x, t)}({\bf 0})$ at the instant $t$ is the same as
that of the Lane-Emden solution inside the ball $B_{x}({\bf 0})$ for $0\le x\le \bar R$. In particular,  the vacuum boundary is given by
$$ R(t)=r(\bar R, t).$$
The convergence of $r(x,t)$ to $x$ in \ef{ma4} means that every spherical surface
moving with the fluid converges to that inside the domain of Lane-Emden solution enclosing the same mass, in particular, the evolving vacuum boundary $R(t) $ converges to the vacuum boundary $\bar R$ as time goes to infinity. This also
gives the large time asymptotic convergence of every particle moving with the fluid since the motion is radial.
Moreover, the convergence \ef{ma4} means that the large time asymptotic states for the free boundary problem \ef{103}  are  determined completely by the initial mass distribution and total mass. Besides the convergence \ef{ma4}, we also establish  convergence rates of  higher-order norms involving  derivatives, and show that the vacuum boundary $R(t)$ has the regularity of $W^{2, \infty}([0 ,\infty)) $ under a compatibility condition of the initial data with the boundary condition which implies that the acceleration of the vacuum boundary is uniformly bounded for $t\in [0, \infty)$. (Indeed, one may check from the proof that every particle moving with the fluid has the bounded acceleration for $t\in [0, \infty)$.)  These results give a rather clear and complete characterization of the behavior of solutions both in large time and near the vacuum boundary.

The results obtained in the present work  are among  few results of   global {\it strong} solutions to vacuum free boundary problems of compressible fluids capturing the singular behavior of \eqref{156}, which is difficult and challenging due to the degeneracy caused by the physical vacuum and   coordinates singularity at the center of the symmetry.  We overcome this difficulty by establishing  higher-order estimates involving the second-order derivatives of the velocity field, together with   decay estimates of lower-order norms. This is achieved  by combining some new weighted
nonlinear functionals (involving both lower-order and higher-order derivatives) and space-time weighted energy estimates.
The constructions of these weighted nonlinear functionals  and space-time weights depend crucially on the structure of  Lane-Emden solutions (in particular, the behavior \ef{pvforle} near the vacuum boundary), the balance between the pressure and self-gravitation, and the dissipation of the viscosity.

We sketch here the main ideas and methods used in this article. The original free boundary problem \ef{103} is reduced to an initial boundary value problem on a fixed domain $x\in [0, \bar R]$ by  the Lagrangian particle trajectory formulation \ef{ma1} and \ef{ma2} with $\bar R$ being the radius of the Lane-Emden solution so that the domain of the Lane-Emden solution becomes the reference domain. In this formulation, essentially the basic unknown is the particle trajectory $r(x, t)$ defined by \ef{ma1} and \ef{ma2} (more precisely, the radius of each evolving surface inside the evolving domain, which is called  the particle trajectory for simplicity here and from now on), by which the density and velocity are determined. For  problem \ef{103}, this formulation is preferred because one can use it to trace each particle in the evolving domain, in particular, the evolving vacuum boundary.
(Indeed, the approach of using the Lagrangian particle trajectory formulation, i.e., using  the flow map of the Eulerian velocity field, to reduce free boundary problems of compressible fluids set on time-dependent domains to fixed time-independent
domains was first adopted in \cite{10}, and latter on in \cite{10', jm, 17'}.)
To make this strategy work, a crucial point  is to obtain the positive lower and upper bounds of the derivative of the particle trajectory, which are derived by a pointwise estimate away from the center of  symmetry, and an interior $L^2$-estimate of its second derivative.
Various multipliers  are applied to establish the decay estimates of lower-order norms and the regularity near the vacuum boundary.  In order to obtain the higher-order estimates, we first study the problem obtained by differentiating the original problem in the tangential direction, and then get weighted estimates for the second derivatives of the velocity field using the viscosity term. However, due to the degeneracy of \eqref{156},  the dissipation  of the viscosity alone is not enough for the global-in-time estimates,  and we have to make full use of the balance between the  pressure and  gravitation.  For this purpose, we decompose the gradient of the pressure as two parts, the first part is to balance the gravitation and the second part  is an anti-derivative of the viscosity  along the particle trajectory with a degenerate weight to match the viscosity.    The main ideas and strategy of establishing the decay estimates and high-order regularity
estimates will be given in Section \ref{sec3.2}.

\subsection{Review of related  works}
There have been extensive works on the studies of the Euler-Poisson and the Navier-Stoke-Poisson equations with vacuum, especially in recent years. We will concentrate on those closely related to the stability of vacuum dynamics.
The stability problem has been important in the theory of gaseous stars which has been studied extensively  by astrophysicists  (cf.  \cite{ch,weinberg,lebovitz1}). The linear stability of Lane-Emden solutions was studied in \cite{linss}.   A conditional nonlinear Lyapunov type
stability theory of stationary solutions  for $\gamma> 4/3$ was established  in \cite{rein} using a variational approach,  by assuming the existence of global solutions of the Cauchy problem for  the three-dimensional compressible Euler-Poisson equations (the same type of nonlinear stability results for rotating
stars were given  by \cite{luosmoller1,luosmoller2}.  For $\gamma\in (6/5, 4/3)$, the nonlinear dynamical instability of Lane-Emden solutions was proved by \cite{17'} and \cite{jangtice} in the framework of free boundary problems for Euler-Poisson systems   and Navier-Stokes-Poisson equations, respectively.  A nonlinear instability for $\gamma= {6}/{5}$ was proved by \cite{jang65}.  For $\gamma= {4}/{3}$, an instability was identified in \cite{DLYY} that a small
perturbation can cause part of the mass to go off to infinity for  inviscid flows.

It should be noted that the stability result in \cite{rein} is  in the framework of initial value problems  in the entire $\mathbb{R}^3$-space and  involves only a Lyapunov functional which is essentially equivalent to a $L^p$-norm of difference of solutions, and the vacuum boundary cannot be traced. Another interesting work is on the   vacuum free boundary problem of   modified compressible Navier-Stokes-Poisson equations with spherical symmetry  (cf. \cite{fangzhang1}), where the existence of a global weak solution was proved for a reduced  initial boundary value problem after using  the Lagrangian mass coordinates,   under some constraints  on the ratio of the  coefficients of the shear viscosity and bulk viscosity.   In contrast to the strong stability result in \ef{ma4}, for the global weak solutions obtained in \cite{fangzhang1}, only the uniform convergence of the velocity $u(r, t)$ is proved, due to the lack of regularity near the vacuum boundary.  The ideas and techniques developed in this paper can be applied to this modified  compressible Navier-Stokes-Poisson equations to obtain a strong stability result as in \ef{ma4}. Indeed, our global-in-time regularity gives not only the decay estimates for the weighted norms  $\|\bar{\rho}^{\gamma/2}(r_x(x, t)-1)\|_{L^{2}_x([0, \bar R])}$
and $ \|x v_x \|_{L^{2}_x([0, \bar R])}  $
 as in \cite{fangzhang1}, but also the decay
estimates of the unweighted norms of
$\| r_x(x, t)-1 \|_{L^{2}_x([0, \bar R])}$,  $ \|v_x\|_{L^{\infty}_x([0, \bar R])}$
 and some uniform estimates on the second derivatives valid up to the vacuum boundary, which are crucial to the nonlinear asymptotic stability for
this modified model.    Furthermore, our theory holds without the restrictions on the viscosity coefficients  as in \cite{fangzhang1}.  This  will be reported in a forthcoming paper (cf. \cite{LXZ2}).

The results obtained in the present work are for spherically symmetric motions. For the problem of general three dimensional
perturbations, in addition to motions in the normal direction, one will have to estimate motions in the tangential direction. Besides this,
the evolution of the geometry of the free surface will be estimated also. We believe the ideas and techniques developed in this paper
will be useful to the study of the problem of  general three dimensional perturbations.

We conclude the introduction  by noting that there are also other  prior results on   free boundary problems involving vacuum for compressible Navier-Stokes equations besides the ones aforementioned.  For the one-dimensional motions, there are many results concerning global weak solutions to free boundary problems of the Navier-Stokes equations, one may refer to   \cite{Okada,Okada3,LXY,fangzhang,JXZ,duan,YYZ,yangzhu,JXZ,zhu}  and references therein. As for the spherically symmetric motions, global existence and stability of weak solutions were obtained in \cite{Okada1,OSM} to  compressible Navier-Stokes equations for gases  surrounding a solid ball (a hard core) without self-gravitation (see also \cite{Chengq} for a compressible heat-conducting flow).  However, those results are restricted to cut-off domains excluding a neighborhood of the origin. It should be noted  that for a
modified system of Navier-Stokes equations, a global existence of weak solutions with spherical symmetry containing the origin was established in \cite{GLX} for which the density does not vanish on the boundary.  For a class of free boundary problems of compressible Navier-Stokes-Poisson equations away from vacuum states, the readers may refer to \cite{94, 95} for the local-in-time well-posedness results and \cite{96} for linearized stability results of stationary solutions.

\vskip 0.25cm
The rest of the paper is organized as follows. In Section \ref{sec2}, we reduce the free boundary problem \ef{103} to a fixed domain by using the Lagrangian formulation and state the global existence and nonlinear asymptotic stability theorems. The {\it a priori} estimates are derived in Section \ref{sec3}, whose main ideas and steps are outlined in Section \ref{sec3.2}. The global existence theorem is proved in Section \ref{sec3}.  Section \ref{sec4} is devoted to the proof of the nonlinear asymptotic stability theorem in the original Eulerian coordinates which is an easy consequence of the estimates obtained in Section \ref{sec3}.  The local existence of strong solutions in our function space framework given in Section \ref{sec2} is proved in Appendix,  Part I. In Appendix, Part II, a linearized analysis is given to illustrate some ideas for the original nonlinear problem.

\section{ Lagrangian formulation and main results}\label{sec2}
\subsection{Lagrangian formulation}
In this subsection, we adopt the  the Lagrangian particle trajectory formulation as  first used in \cite{10} for inviscid flows  and latter on in \cite{10', jm, 17'} to reduce the original free boundary problem \ef{103} to an initial boundary value problem on the fixed domain $x\in [0, \bar R]$.

For this purpose,  we first recall some properties of  Lane-Emden solutions. For $\ga\in(4/3, 2)$, it is known that for any given finite positive total mass, there exists a unique solution to equation \eqref{le} whose support is compact (cf. \cite{linss}). Without abusing notations,  $x$ will  denote the distance from the origin for the Lane-Emden solution. Therefore, for any $M\in(0,\infty)$, there exists a unique function $\bar\rho(x)$ such that
\begin{equation}\label{lex1}
\bar\rho_0:=\bar\rho(0)>0, \ \  \bar\rho(x)>0 \ \  {\rm for} \ \ x\in \lt(0,\ \bar R\rt),   \ \  \bar\rho\lt(\bar R\rt)=0, \ \  M=\int_0^{\bar R} 4\pi \bar \rho(s) s^2ds  ;
 \end{equation}
 \begin{equation*}\label{newle}
-\iy<\bar\rho_x<0 \ \ {\rm for} \ \  x\in (0,\ \bar R) \ \ {\rm and} \ \ \bar\rho(x) \le \bar\rho_0 \ \  {\rm for} \ \ x\in \lt(0,\ \bar R\rt);
\end{equation*}
\begin{equation}\label{rhox}
\left(\bar{\rho}^\ga \right)_x=-x \phi \bar{\rho}, \ \ {\rm where} \  \ \phi:= x^{-3}\int_0^x 4\pi \bar\rho(s)  s^2 ds \in \left[M/{\bar R}^3, \   4\pi \bar\rho_0/3\right] ;
\end{equation}
for a certain finite positive constant $\bar R$ (indeed, $\bar R$ is determined by $M$ and $\ga$). Note that
$$
\left(\bar{\rho}^{\ga-1} \right)_x=\ga^{-1}(\ga-1)    {\bar \rho}^{-1} \left(\bar{\rho}^{\ga} \right)_x
= - \ga^{-1}(\ga-1) x\phi.
$$
It then follows from \eqref{lex1} and \eqref{rhox} that $\bar\rho$ satisfies the physical vacuum condition, that is,
$$
 \bar{\rho}^{\ga-1}(x)  \sim  \bar R- x \ \ {\rm as} \ \ x \ \ {\rm~ close~ to~} \bar R.
$$
More precisely, there exists a constant $C$ depending on $M$ and $\ga$ such that
\begin{equation}\label{phy}
C^{-1} \lt( \bar R- x \rt) \le \bar{\rho}^{\ga-1}(x)  \le C \lt( \bar R- x \rt), \ \ x\in \lt(0,\ \bar R\rt).
\end{equation}

The  particle trajectory Lagrangian formulation for  \eqref{103} is given as follows.   Let  $x$ be the reference variable and define the Lagrangian variable $r(x, t)$ by
\be\label{Aug9-1}
  r_t(x, t)= u(r(x, t), t) \ \ {\rm for} \  \ t>0  \ \  {\rm and} \ \ r(x,0)=r_0(x), \ \  x\in I:=\lt(0, \bar R\rt)  .
\ee
Here  $r_0(x)$ is the initial position which maps  $ \bar I \to \lt [0, R_0\rt]$  satisfying
\begin{equation}\label{rox} \int_0^{r_0(x)} \rho_0(s) s^2ds = \int_0^x   \bar\rho(s)  s^2 ds, \ \ x\in \bar I, \end{equation}
so that
\begin{equation}\label{choice}
\rho_0(r_0(x))r_0^2(x)r_0'(x)=\bar\rho(x)x^2, \ x\in \bar I.\end{equation}
(Indeed, \ef{rox} means that the initial mass in the ball with the radius $r_0(x)$ is the same as that of the Lane-Emden solution in the ball
with the radius $x$. Then smoothness of $r_0(x)$ at $x=\bar R$ is equivalent to that the initial density $\rho_0$ has the same behavior
near $R_0$ as that of $\bar\rho$ near $\bar R$.)
The choice of  $r_0$  can be described by
\be\label{r000} r_0(x)=\psi^{-1}(\xi(x)), \ \ 0\le x\le \bar R; \ee
where $\xi$ and $\psi$ are one-to-one mappings, defined by
$$\xi:  (0, \bar R) \to (0, M): \ x \mapsto \int_0^x s^2\bar \rho (s)d s  \ \ {\rm and} \ \
  \psi:  (0, R_0) \to (0, M): \ z \mapsto \int_0^z s^2 \rho_0 (s)ds. $$
Moreover $r_0(x)$ is an increasing  function and
the initial total mass has to be the same as that for $\bar \rho$, that is,
\begin{equation}\label{samemass}
\int_0^{R_0}4\pi  \rho_0(s) s^2 ds=\int_0^{r_0(\bar R)}4\pi  \rho_0(s) s^2 ds =
\int_0^{\bar R} 4\pi  \bar \rho (s) s^2 ds=M,
\end{equation}
to ensure that $r_0$ is a diffeomorphism from $\bar I$ to $\lt [0, R_0\rt]$.
In view of  \eqref{103a}, we see
\begin{equation}\label{masswithin}
\int_0^{r(x, t)}\rho(s, t)s^2ds=\int_0^{r_0(x)} \rho_0(s) s^2ds, \ \  x\in I.
\end{equation}

Define the Lagrangian density and velocity respectively by
$$f(x, t)=\rho(r(x, t), t) \ \ {\rm and} \ \  v(x, t)=u(r(x,t), t).$$
Then  the Lagrangian version of  \eqref{103a} and \ef{103b} can be written on the reference domain $I$ as
\begin{subequations}\label{new419} \begin{align}
& (r^2f)_t +r^2f\frac{v_x}{r_x}=0    & {\rm in}& \  \    I\times (0, T],\label{new419a}\\
& f  v_t+ \frac{ (f^{\gamma})_x}{  r_x}+ {4\pi f} {r^{-2}} \int_0^{r_0(x)} \rho_0(s) s^2ds = \frac{\mu}{r_x } \left(\frac{ (r^2v)_x}{r^2  r_x}\right)_x \ \  &{\rm in}& \  \   I\times (0, T]. \label{new419b}
\end{align}
\end{subequations}
Solving \eqref{new419a} gives   that
$$f(x, t)r^2(x, t)  r_x(x, t)= \rho_0(r_0(x)) r_0^2(x)  r_{0x}(x), \ \  x\in I.$$
Therefore,
$$
f(x, t)=  \frac{x^2\bar \rho(x)}{r^2(x, t)  r_x(x, t)}  \ \ {\rm for} \   \ x\in  I,
$$
 due to   \ef{choice}.
So, \eqref{103}  can be written on the reference domain $I=\lt(0, \bar R\rt)$ as
\begin{subequations}\label{419} \begin{align}
& \bar\rho\left( \frac{x}{r}\right)^2  v_t   +  \left[    \left(\frac{x^2}{r^2}\frac{\bar\rho}{ r_x}\right)^\ga    \right]_x +  \frac{x^2}{r^4}    \bar\rho \int_0^x 4\pi y^2\bar\rho(y)   dy =\mu  \left(\frac{ (r^2v)_x}{r^2  r_x}\right)_x  & {\rm in} & \    I\times (0, T],  \label{419a}\\
&  v(0, t)=0, \ \   \mathfrak{B}(\bar R, t)=0    & {\rm on}  & \ [0,T],\label{419b}\\
& (r,\ v)(x, 0) = \lt(r_0(x), \  u_0(r_0(x)) \rt)  &  {\rm on}  & \     I \times \{t=0\}, \label{419c}
\end{align}
\end{subequations}
where $\mathfrak{B}$ is the normal stress at the boundary given by
\begin{equation}\label{bdry1}
\mathfrak{B}:=\frac{4}{3}\la_1\left(\frac{v_x}{r_x}-\frac{v}{r}\right)  + \la_2 \left(\frac{v_x}{r_x}+2\frac{v}{r}\right)
=\frac{4}{3}\la_1\frac{r}{r_x}\left(\frac{v}{r}\right)_x  + \la_2 \frac{\left(r^2 v\right)_x}{r_x r^2}.
\end{equation}
Moreover, it can be derived from \ef{Aug9-1}, \ef{rox} and \ef{419b} that
\be\label{Aug9-2}
r(0,t)=r_0(0)+\int_0^t v(0,s) ds =0 \ \  {\rm on} \ \ [0, T].
\ee

{\bf Notation}. Throughout the rest of this paper,  $c$ and $C$ will be used to denote  generic positive constants which are independent of time $t$ but may depend on $\gamma$, $\lambda_1$, $\lambda_2$, $M$ and the bounds of $\bar \rho$ such as  $\bar\rho(0)$ and $\bar\rho(\bar R/2)$;
and we will use the following notations:
$$\int: =\int_{I},  \ \ \|\cdot\| :=\|\cdot\|_{L^2(I)}, \ \ \|\cdot\|_{L^p}:=\|\cdot\|_{L^p(I)} \   (p=1, \iy),  \ \    {\rm and} \ \ \|\cdot\|_{H^s}:=\|\cdot\|_{H^s(I)} \   (s=1, 2) .$$

\subsection{Strong solutions and  functionals}\label{sec2.2}
A strong solution to  problem \ef{419} is defined as follows.
\begin{defi}\label{definitionss}  $v\in C\lt([0, T]; H^2_{loc}([0, \bar R))\rt)\cap C\lt([0, T];  W^{1, \infty}(I)\rt)$ with
\be\label{r}
r(x, t)=r_0(x)+\int_0^t v(x, s)ds  \ \ for \ \    (x, t)\in I\times [0, T] ,
\ee
satisfying the initial condition \ef{419c} is called a strong solution of problem \ef{419} in $[0, T] $, if\\
1) $c_1\le r_x(x, t) \le c_2$, $(x, t)\in I\times [0, T]$, for some positive constants  $c_1$ and $c_2$;  \\
2) $\bar\rho ^{ - {1}/{2}} \lt[   ({r^2r_x}  )^{-1}{(r^2 v)_x} \rt]_x\in C([0, T]; L^2(I))$  and  $\bar\rho ^{\ga- {1}/{2}} \lt(r_{xx },   \lt(r/x\rt)_x\rt)\in C^1([0, T]; L^2(I))$;\\
3) $\bar\rho^{ {1}/{2}} v \in C^1([0, T]; L^2(I))$;\\
4) $v(0, t)=0$ and $\mathfrak{B}(\bar R, t)=0$  hold   in   the  sense   of $W^{1, \infty}$-trace and $H^1$-trace, respectively,   for   $t\in [0, T]$;\\
5) \ef{419a}   holds    for  $ (x, t)\in I\times [0, T]$,  a.e..
\end{defi}

\begin{rmk}\label{9.20} Let $(r,v)$ be a strong solution of problem \ef{419} defined in Definition \ref{definitionss}, then it holds that for any $a\in (0 ,\bar R)$,
\begin{align}
& r\in C^1\lt([0, T];  W^{1, \infty}(I)\rt)  \cap C^1 \lt([0, T];  H^2\lt([0,a]\rt)\rt),  \ \  {r}/{x}  \in C^1 \lt([0, T];  H^1(I)\rt), \label{hhreg1}\\
& v \in C \lt([0, T];  W^{1, \infty}(I)\rt)  \cap C  \lt([0, T];  H^2\lt([0,a]\rt)\rt), \ \   {v}/{x}  \in C  \lt([0, T];  H^1(I)\rt), \label{hhreg2}\\
&  \lt(  v/r, \  {v_x}/{r_x}, \   \mathfrak{B}  \rt)  \in C\lt([0, T];  H^1(I)\rt). \label{hhreg3}
\end{align}
The arguments for \ef{hhreg1}-\ef{hhreg3} go as follows. It gives from $c_1\le r_x(x, t) \le c_2$ and $r(0,t)=0$ that $c_1\le  x^{-1} {r(x, t)}  \le c_2$. As a consequence of the facts $v\in C\lt([0, T];  W^{1, \infty}(I)\rt)$ and $v(0,t)=0$, one has that $v/x\in C\lt([0, T];  L^\iy(I)\rt)$; which, together with \ef{r} and  $c_1\le r/x \le c_2$, gives  $r/x\in C^1\lt([0, T];  L^\iy(I)\rt)$. Moreover, it follows from  \ef{r},  $\bar\rho^{\ga- {1}/{2}} \lt(r_{xx },   \lt(r/x\rt)_x\rt)\in C^1([0, T]; L^2(I))$, and $\bar\rho(x) \ge \bar\rho (a)$ for $x\in [0, a]$ that
\be
\lt( r_{xx}, \ \lt(r/x\rt)_x \rt) \in
C^1\lt([0, T]; L^2 \lt([0,a]\rt)\rt) \ \ {\rm and} \ \ \lt( v_{xx}, \ \lt(v/x\rt)_x \rt) \in
C\lt([0, T];  L^2 \lt([0,a]\rt)\rt).
\ee
So, \ef{hhreg1} and \ef{hhreg2} hold, since $|(r/x)_x|\le a^{-1}(|r_x|+|r/x|)$ for $x\in [a, \bar R]$. It follows from $v/r=(v/x)(x/r)$, $c_1\le  r/x  \le c_2$, \ef{hhreg1} and \ef{hhreg2} that $v/r \in C\lt([0, T];  H^1(I)\rt)$. Due to 2) of Definition \ref{definitionss}, it holds that $\lt(v_x/r_x + 2 v/r \rt)_x\in C([0, T]; L^2(I)) $. These imply $(v_x/r_x)_x\in C([0, T]; L^2(I))$. So, \ef{hhreg3} holds.
\end{rmk}

In what follows, we give some remarks on the above definition of strong solutions to problem \ef{419}.  First of all, in order to make the transformation $x \mapsto r$ invertible to define  particle trajectories, it is essential to require the positive lower and upper bounds for $r_x$, i.e., $1)$ of Definition \ref{definitionss}. The requirement $v\in C([0, +\infty); W^{1, \infty}(I))$ is  also essential, because it is also related to the well-definiteness of particle trajectories defined by  \ef{Aug9-1}, whose uniqueness for the given initial values $r_0(x)$ is not ensured without the Lipschitz continuity of $u$ in $r$, and the bound of $u_r$ and $v_x$ are related by the identity $u_r(r(x, t), t)=v_x(x, t)/r_x(x, t)$.  The above definition of strong solutions guarantees that each term in equation $\ef{419a}$ is in $ C([0, T]; L^2(I))$. More than this, each term multiplied by $\bar\rho^{-{1}/{2}}$ is in $ C([0, T]; L^2(I))$.  Indeed, it is quite natural to require $\bar\rho^{{1}/{2}}v_t \in C([0, T]; L^2(I))$ other than $\bar\rho v_t \in C([0, T]; L^2(I))$ from the kinetic energy point of view. The kinetic energy of the system is  $\|x\bar\rho^{{1}/{2}}v(\cdot, t)\|$,  which is expected to be bounded and continuous in time. (Indeed, this can be justified by applying the multiplier $r^2v$ to equation \ef{419a}.) By studying the problem obtained by differentiating the original one with respect to $t$, one may expect $\|x\bar\rho^{ {1}/{2}}v_t(\cdot, t)\|$ is bounded and continuous in time. We can improve this to the boundedness and continuity in time of  $\|\bar\rho^{{1}/{2}}v_t(\cdot, t)\|$  by using the viscosity, which will be shown later. This leads us naturally to require that each term in $\ef{419a}$ multiplied by $\bar\rho^{-{1}/{2}}$ is in $ C([0, T]; L^2(I))$.

To show the well-posedness of strong solutions to problem \ef{419} defined in Definition \ref{definitionss}, we introduce the following higher-order functional:
\be\label{mathmarch}
\mathfrak{E}(t)=\lt\|(r_x-1,  \ v_x)(\cdot,t)\rt\|_{L^\iy}^2  +  \lt\|   \bar\rho^{\ga-{1}/{2}}(r_{xx},\ (r/x)_x)(\cdot, t) \rt\|^2  + \lt\|   \bar\rho^{{1}/{2}} v_t(\cdot, t) \rt\|^2,  \ \  t\ge 0.
\ee
We will prove the global existence and uniqueness of strong solutions satisfying  $\mathfrak{E}(t)\le C\mathfrak{E}(0)$ $(t\ge 0)$ for some constant $C>0$ independent of $t$, and some decay estimates as in $i)$ of Theorem \ref{mainthm},  provided that $\mathfrak{E}(0)$ is suitably small and the  following compatibility condition of the initial data with the boundary conditions holds.
\be\label{compatibility}
v(0, 0)=0 \ \  {\rm and} \ \  \mathfrak{B}(\bar R, 0)=0.
\ee

In order to gain further regularity of strong solutions obtained via $\mathfrak{E}(t)$, we introduce the following  functional:
\begin{align}
&\mathfrak{F}_\alpha(t)=\mathfrak{E}(t)+  \lt\|\bar\rho^{(2\ga-1-\aa)/2} r_{xx}(\cdot, t) \rt\|^2, \ \ \aa\in [0, 2\ga-1],  \ \ t\ge 0.\notag
\end{align}
Besides $\mathfrak{E}(0)$ is small, if the initial data are assumed to satisfy further regularity that $\mathfrak{F}_\alpha(0)<\infty$ for $0<\alpha\le2\ga-1$, we will prove further regularity and decay (than those stated in $i)$ of Theorem \ref{mainthm}) of strong solutions with decay rates which may depend on $\alpha$ in various norms, as shown in $ii)$ and $iii)$ of Theorem \ref{mainthm}. Indeed, the index $\alpha$ indicates the behavior of solutions near the vacuum boundary, which has influence on the decay rates of solutions. In particular, we will prove the  regularity that $v_{xx}(\cdot, t)\in L^2(I)$ for all $t\ge 0$ and the decay of $\|v_{xx}(\cdot, t)\|$, if
the initial data satisfy additional regularity that  $r_{xx}(\cdot, 0)\in L^2(I)$ (i.e.,  $\mathfrak{F}_{2\ga-1}(0)<\iy$).

Some remarks are given on the smallness requirement of  functional $\mathfrak{E}(t)$.  For problem \ef{419}, $r(x, t)=x$ is the equilibrium solution. The basic stability requirement is on the smallness of $r(x,t)-x$ for all time in certain topology. The smallness of the quantity $r_x(x, t)-1$, the $x$-derivative of $r(x,t)-x$,  in $L^{\infty}$-norm, ensures that $r_x(x, t)$ is bounded from below and above by positive constants. This also gives the smallness of $\|r/x-1\|_{L^{\infty}}$ and thus the  lower and positive bounds for $r/x$ due to the condition $r(0, t)=0$ (see \ef{Aug9-2}). So, it is essential to obtain the smallness of $\|r_x-1\|_{L^{\infty}}$.  Since $v=0$ is the equilibrium for problem \ef{419},
it is required that  $\|v(\cdot, t)\|_{W^{1, \infty}(I)}$ to be small for all $t>0$ if this holds true initially, from the stability point of view. It is crucial to derive the smallness of $L^{\infty}$-bound for $v_x$, which also gives the bound for $v/x$ due to the  condition that $v(0, t)=0$. Our basic strategy is to use the weighted $L^2$-energy method together with some pointwise estimates to bound $\|r_x-1\|_{L^{\infty}}$ and $\|v_x\|_{L^{\infty}}$. For example, we may bound $\|r_x-1\|_{L^{\infty}([0, \bar R/2])}$ by use of $\|r_{xx}\|_{L^2([0, \bar R/2])}$ which can be bounded by $\| \bar\rho^{\ga-1/2} r_{xx}  \|$. Therefore, that the functional $\mathfrak{E}(t)$ is small for all time is essential for our stability analysis of problem \ef{419}.

We give here some consequences of  the smallness of $\mathfrak{E}(t)$ ($t\ge 0$), which ensure that  each term in $\ef{419a}$ multiplied by $\bar\rho^{-{1}/{2}}$ is in $L^2(I)$ ($t\ge 0$), and that the boundary condition $v(0,t)=0$ and $\mathfrak{B}(\bar R, t)=0$ are achieved in the sense of $W^{1,\iy}$-trace and $H^1$-trace, respectively.  For $t\ge 0$ and $a\in (0, \bar R)$, it holds that
\begin{align}
& \lt\|(r/x-1)(\cdot,t)\rt\|_{L^\iy}^2 \le \lt\|(r_x-1)(\cdot,t)\rt\|_{L^\iy}^2 \le \mathfrak{E}(t), \label{9.21.1}\\
& 1/2 \le r(x,t)/x \le 3/2 \ \ {\rm and} \ \  1/2 \le r_x(x,t) \le 3/2 , \ \ x\in I, \label{9.21.2}\\
& \lt\|(v/x)(\cdot,t)\rt\|_{L^\iy}^2 \le \lt\|v_x(\cdot,t)\rt\|_{L^\iy}^2 \le \mathfrak{E}(t) ,\label{9.21.3}\\
&\lt\|\bar\rho ^{ - {1}/{2}} \lt(v_x/r_x + 2 v/r \rt)_x (\cdot, t) \rt\|^2= \lt\|\bar\rho ^{ - {1}/{2}} \lt[   ({r^2r_x}  )^{-1}{(r^2 v)_x} \rt]_x (\cdot, t) \rt\|^2 \le C\mathfrak{E}(t),\label{9.21.4}\\
& \lt\|   \bar\rho^{\ga-{1}/{2}}(v_{xx}, \ (v/x)_x)(\cdot, t) \rt\|^2 \le C\mathfrak{E}(t), \label{9.21.5}\\
& \lt\|(v, r-x)(\cdot, t)\rt\|_{H^2([0, a])}^2 \le C(a)\mathfrak{E}(t), \ \ \lt\|\lt(\frac{r}{x}-1, \frac{v}{x}, \frac{v}{r}, \frac{v_x}{r_x}, \mathfrak{B}\rt)(\cdot, t)\rt\|_{H^1}^2 \le C\mathfrak{E}(t),\label{9.21.6}
\end{align}
where $C$ and $C(a)$ are positive constants independent of $t$.

The arguments for \ef{9.21.1}-\ef{9.21.6} go as follows.    \ef{9.21.1} and \ef{9.21.3} follow from $r(0,t)=0$ and $v(0,t)=0$, respectively; \ef{9.21.2} follows from the smallness of $\mathfrak{E}(t)$;  \ef{9.21.4} follows from \ef{9.21.1},  \ef{9.21.2}, the definition of $\mathfrak{E}(t)$ in \ef{mathmarch}, and equation \ef{419a} which can be rewritten as
$$\mu  \left(\frac{ (r^2v)_x}{r^2  r_x}\right)_x = \bar\rho\left( \frac{x}{r}\right)^2  v_t   + \bar\rho^\ga \left[    \left(\frac{x^2}{r^2}\frac{1}{ r_x}\right)^\ga    \right]_x + \lt[\left(\frac{x^2}{r^2}\frac{1}{ r_x}\right)^\ga  - \frac{x^4}{r^4} \rt]  \left(\bar{\rho}^\ga\right)_x   ; $$
\ef{9.21.5} follows from \ef{9.21.3},  \ef{9.21.4}, the definition of $\mathfrak{E}(t)$ in \ef{mathmarch},  and the estimate below:
\be\label{9.21.7}
\lt\|\bar\rho^{\ga-\frac{1}{2}}\lt(v_{xx}, \  \lt(\frac{v}{x}\rt)_x \rt) \rt\|^2 \le C \lt\| \lt(\frac{v_x}{r_x} + 2 \frac{v}{r}\rt)_x  \rt\|^2+ C \lt\|\lt(v_x,   \frac{v}{x}\rt)\rt\|^2_{L^\iy}\lt\|\bar\rho^{\ga-\frac{1}{2}}\lt(r_{xx}, \  \lt(\frac{r}{x}\rt)_x \rt) \rt\|^2  .
\ee
(Indeed, \ef{9.21.7} follows from \ef{tlg2}, \ef{tlg5}, \ef{9.21.2}, and \ef{weightvxx} which will be prove later.) In a similar argument as that for Remark \ref{9.20}, one can obtain \ef{9.21.6}.

\subsection{Main theorems and remarks}

The first  theorem of this paper is on  the global existence   and the regularity  of  strong solutions, which also gives the strong Lyapunov stability of the
Lane-Emden solution  in the functional $\mathfrak{E}(t)$:

\begin{thm}\label{mainthm1} Let $\ga\in(4/3,\ 2)$ and $\bar\rho$ be the Lane-Emden solution satisfying \ef{lex1}-\ef{rhox}. Assume that \ef{compatibility} holds and the initial density $\rho_0$ satisfies \eqref{156} and \ef{samemass}.
There exists a   constant $\bar\da >0$ such that if
$\mathfrak{E}(0)\le \bar\da,$
then   the  problem \eqref{419}  admits a unique strong solution  in $I\times[0, \iy)$ with
\be\label{keyconclusion} \mathfrak{E}(t)\le C\mathfrak{E}(0), \ \  t\ge 0, \ee
for  some constant $C$ independent of $t$.

Moreover, if $\mathfrak{F}_{2\ga-1}(0)<\iy$ (i.e., $\|r_{xx} (\cdot, 0)\|<\iy$), then the strong solution obtained above satisfies the following  further regularity estimates
\be\label{Aug23-1}
\|v_{xx}(\cdot, t)\|^2 \le C \mathfrak{F}_{2\ga-1} (0) \lt(1+ \mathfrak{F}_{2\ga-1} (0)   \rt)  \lt(1+ \mathfrak{E} (0)  \rt), \ \ t\ge 0,
\ee
for some constant $C$ independent of $t$; and
\be\label{Aug23-2}
\|r_{xx}(\cdot, t)\|^2 \le C \mathfrak{F}_{2\ga-1} (0)  +  C(T) \mathfrak{E}(0), \ \ t\in [0, T],
\ee
for some constants $C$ independent of $t$ and $C(T)$ depending on $T$.
\end{thm}

For any $t\ge 0$, since $r_x(x, t)>0$ for $x\in \bar I$, $r(x, t)$ defines a diffeomorphism from the reference domain  $\bar I$ to the changing domain $\{0\le r\le R(t)\}$ with the boundary
\be\label{vacuumboundary}
R(t)=r\lt(\bar R ,  t\rt).
\ee
 It also induces a diffeomorphism from the initial domain, $\bar B_{R_0}(0)$,
 to the evolving domain, $\bar B_{R(t)}(0)$,
  for all $t\ge 0$:
$${\bf x}\ne {\bf 0} \in \bar B_{R_0}(0)\to r\lt(r_0^{-1}(|{\bf x}|),  t\rt)|{\bf x}|^{-1}{{\bf x}}  \in  \bar B_{R(t)}(0), $$
where $r_0^{-1}$ is the inverse map of $r_0$ defined in \ef{r000}. Here
$$\bar B_{R_0}(0):= \{{\bf x}\in \mathbb{R}^3: |{\bf x}|\le R_0\} \ \  {\rm and} \  \ \bar B_{R(t)}(0):=\{{\bf x}\in \mathbb{R}^3: |{\bf x}|\le R(t)\}. $$
 Denote  the inverse of the map $r(x, t)$ by $\mathcal{R}_t$ for  $t\ge 0$ so that
$$
{\rm if~} \ \  r=r(x, t)\ \ {\rm  for ~} \ \  0\le r\le R(t), \ \   {\rm then~} \ \ x=\mathcal{R}_t(r).$$
For the strong solution $(r, v)$ obtained in Theorem \ref{mainthm1}, we set for $0\le r\le R(t)$ and $t\ge 0$,
\be\label{solution}
  \rho(r, t)=\frac{x^2\bar \rho(x)}{r^2(x, t)  r_x(x, t)}  \ \  {\rm and} \ \
   u(r, t)=v(x, t) \ \  {\rm with} \ \     x=\mathcal{R}_t(r).  \ \
   \ee
Then the triple $(\rho(r,t), u(r,t), R(t))$  ($t\ge 0$) defines a global strong solution to the free boundary problem \ef{103}. Furthermore,  we have the strong nonlinear asymptotic stability of the Lane-Emden solution as follows.

\begin{thm}\label{mainthm2} Under the assumptions in Theorem \ref{mainthm1}. Then the triple $(\rho, u, R(t))$   defined by
\ef{vacuumboundary} and \ef{solution} is the unique global strong solution to the free boundary problem \eqref{0.1} satisfying
$R\in W^{1, \infty}( [0, \ +\infty) ).$  Moreover, the solution satisfies the following  estimates.

i)  For any $0<\theta< {2(\ga-1)}/({3\ga})  $, there exists a positive constant $C(\theta)$  independent of $t$ such that for all $t\ge 0$,
\begin{align}
& \sup_{0\le x\le  \bar R} |r(x, t)-x|\le C(\theta) (1+t)^{-\frac{\ga-1}{\ga}+\frac{\ta}{2}}\sqrt {\mathfrak{E}(0)}, \label{rr1} \\
& \sup_{0\le r\le  R(t)} \lt|u(r, t)\rt| \le C(\theta) (1+t)^{-\frac{3\ga-2}{4\ga}+ \frac{\ta}{2} }\sqrt {\mathfrak{E}(0)}, \label{estthm1b} \\
& \sup_{0\le r\le  R(t)} \lt|\lt(u_r, \  r^{-1} u\rt)(r, t)\rt|  \le C(\theta) (1+t)^{-\frac{\ga-1}{2\ga}+ \frac{\ta}{2}}\sqrt {\mathfrak{E}(0)}, \label{estthm1b''} \\
& \sup_{0\le x\le  \bar R }\lt|\lt(\bar\rho(x) \rt)^{({3\ga-6})/{4}} \lt[\rho (r(x, t), t)-\bar\rho(x) \rt]\rt|\le
  C(\theta)    (1+t)^{-\frac{ \ga -1}{2\ga}+ \frac{\ta}{2} } \sqrt {\mathfrak{E}(0)}. \label{estthm1d}
\end{align}

ii) Suppose that $\mathfrak{F}_\alpha(0)<\iy$ for some $\alpha\in (0, \ga)$. Let $\theta$  be any constant satisfying
$0<\theta<  \min\lt\{{2(\ga-1)}/({3\ga}), \ \  2(\ga-\alpha)/\ga \rt\}$. Set $\kappa = (1/\ga) \min\{  {\alpha-(\ga-1)} , \   {\ga-1} \} -\theta$ when $\alpha\in (\ga-1, \ga)$, and $\kappa=0$ when $\alpha = \ga-1$. Then  there exists a positive constant  $ C(\alpha,\ta)$  independent of $t$  such that  for all $t\ge 0$,
\be\label{estthm1d8.23}   \sup_{0\le x\le  \bar R }\lt|\lt(\bar\rho(x) \rt)^{({3\ga-6}-\alpha)/{4}}\lt[\rho (r(x, t), t)-\bar\rho(x) \rt]\rt|\le
 C(\alpha,\ta)    (1+t)^{-\frac{ \ga -1}{2\ga}+ \frac{\ta}{2} } \sqrt {\mathfrak{F}_\alpha (0)};
\ee
and if $\alpha\in [\ga-1, \ga)$,
\begin{align}
&\sup_{0\le r\le  R(t)} \lt|u(r, t)\rt|\le C(\alpha,\ta) (1+t)^{-\frac{8\ga-5}{4\ga}-\frac{\ka}{4}+\frac{5}{4}\ta } \sqrt{ {\mathfrak{F}}_\alpha(0) +
  {\mathfrak{E}}(0) {\mathfrak{F}}_\alpha(0)}, \label{estthm2b} \\
& \sup_{0\le r\le  R(t)} \lt|(u_r, \ r^{-1} u)(r, t)\rt|\le C(\alpha,\ta) (1+t)^{-({1}/{2})\min\{b_1,  b_2\} } \sqrt{ {\mathfrak{F}}_\alpha(0) +
  {\mathfrak{E}}(0) {\mathfrak{F}}_\alpha(0)},\label{estthm2b'''} \\
&\sup_{0\le x\le  \bar R}\lt|\lt(\bar\rho(x) \rt)^{(\ga-2)/2}\lt[\rho (r(x, t), t)-\bar\rho(x) \rt]\rt| \notag \\
 &\qquad \le  C(\alpha,\ta)(1+t)^{-\frac{\ka}{2}- \frac{2\ga-1}{2\ga}+\frac{\ta}{2}} \sqrt{{\mathfrak{F}}_\alpha(0) +
  {\mathfrak{E}}(0) {\mathfrak{F}}_\alpha(0)}. \label{estthm2d}\end{align}
Here
\begin{align}
  b_1=&\min\lt\{ \max\lt\{  \lt(\frac{\kappa }{2}+\frac{4\ga-3}{2\ga}-\frac{3}{2}\ta\rt)\frac{\aa+1}{2\ga-1+\aa}, \  \frac{3 }{2}\kappa +\frac{2\ga-1}{2\ga}-\frac{ \ta}{2}\rt\}, \rt.\notag\\
& \lt. \qquad \ \ \frac{\kappa }{2}+\frac{4\ga-3}{2\ga}-\frac{3}{2}\ta  \rt\} + \frac{2\ga-1}{\ga} -\ta ,\label{newb1} \\
  b_2 =  &   \min\lt\{\frac{\kappa }{2}+\frac{4\ga-3}{2\ga}-\frac{3}{2}\ta, \ \frac{\kappa }{4}+\frac{10\ga-9}{4\ga}-\frac{9}{4}\ta\rt\} + \frac{\kappa }{2}+\frac{4\ga-3}{2\ga}-\frac{3}{2}\ta. \label{newb2}
\end{align}

iii) Furthermore, if $\|x\bar\rho^{ {1}/{2}} v_{tt}(\cdot, 0)\|^2   +|v_t(\bar R , 0)|<\infty$, then $R\in W^{2, \infty}([0 ,  +\infty))$ and
\be\label{accelaration}
|\ddot{R}(t)|\le |v_t(\bar R , 0)|+ C(\mathfrak{E}(0))^{1/4}\lt((\mathfrak{E}(0))^{1/4}
+\|x\bar\rho^{ {1}/{2}} v_{tt}(\cdot, 0)\|^{\frac{1}{2}}\rt), \  \
t\ge 0, \ee
for some constant $C$ independent of $t$.
 \end{thm}

\begin{rmk}  In $\mathfrak{E}(0)$, the term $\|\bar\rho^{{1}/{2}}v_t(\cdot, 0)\|$  can be given in terms of the initial data as follows: for $x\in I$,
$$\bar\rho^{\frac{1}{2}}  v_t (x, 0)  =  \bar\rho^{-\frac{1}{2}} \left( \frac{r}{x}\right)^2  \lt\{ \mu  \left(\frac{ (r^2v)_x}{r^2  r_x}\right)_x   - \bar\rho^\ga \left[    \left(\frac{x^2}{r^2}\frac{1}{ r_x}\right)^\ga    \right]_x - \lt[\left(\frac{x^2}{r^2}\frac{1}{ r_x}\right)^\ga  - \frac{x^4}{r^4} \rt]  \left(\bar{\rho}^\ga\right)_x  \rt\}(x,0)  .$$
Indeed, this is equivalent to equation \ef{419a} at $t=0$.
\end{rmk}

\begin{rmk}\label{rmk2'} The estimates in \ef{estthm1d}, \ef{estthm1d8.23} and \ef{estthm2d} yield  the uniform convergence with rates of the density to \ef{0.1} to that of the Lane-Emden solution for both large time and near the vacuum boundary since $\gamma<2$. \end{rmk}

\begin{rmk}\label{rmk3}The initial perturbation here includes three parts: the deviation of the initial domain from that of the Lane-Emden solution, the difference of initial density from that of the Lane-Emden solution, and the velocity.  Since the Lane-Emden solution is completely determined by the total mass $M$, our nonlinear asymptotic stability result shows that the time asymptotic state of the free boundary problem is determined by the total mass which is conserved in the time evolution. \end{rmk}

\begin{rmk}The condition $\|x\bar\rho^{ {1}/{2}} v_{tt}(\cdot, 0)\|^2   +|v_t(\bar R , 0)|<\infty$ in $iii)$ of Theorem \ref{mainthm2} to ensure $R\in W^{2, \infty}([0 ,  +\infty))$ (uniform boundedness of the acceleration of the vacuum boundary) is a higher-order compatibility condition of the initial data with the vacuum boundary.  Indeed, one may check from the proof that every particle moving with the fluid has the bounded acceleration for $t\in [0, \infty)$ if it does so initially.  \end{rmk}

\section{Proof of main results}\label{sec3}
\subsection{A  theorem with detailed estimates}
For the convenience of presentation, we set $\bar R=1$ and
$I=(0, \bar R)=(0, 1).$
Indeed,  we will prove the following results for the global strong solutions obtained in Theorem \ref{mainthm1}, which gives
not only the nonlinear asymptotic stability results stated in Theorem \ref{mainthm2}, but also detailed behavior of the solutions both
in large time and near the vacuum boundary and the origin.

\begin{thm}\label{mainthm} Let $v$ be the global strong solution to the problem \eqref{419} with $r$ given by \ef{r}  as obtained in Theorem \ref{mainthm1}.

i) Let $\theta$ and $\delta$ be any constants satisfying
$0<\theta<  {2(\ga-1)}/({3\ga})$  and $ \delta\in (0, 1)$.
Then there exist positive constants $C(\theta)$  and $C(\theta, \delta)$ independent of $t$ such that for all $t\ge 0$,
\begin{align}\label{estthm1}
& (1+t)^{\frac{2(\ga-1)}{\ga}-{\ta}}\lt\|(r-x)(\cdot,t)\rt\|_{L^\iy}^2
 +(1+t)^{ \frac{3\ga-2}{2\ga}-{\ta} }\lt\|(v,xv_x)(\cdot,t)\rt\|_{L^\iy}^2
 \notag\\
& + (1+t)^{\frac{\ga-1}{ \ga}- \ta }
 \lt\|\bar\rho^{({3\ga-2})/{4}}  \lt(r_x-1,  {r}/{x}-1 \rt)(\cdot,t)\rt\|_{L^\iy}^2+ (1+t)^{\frac{\ga-1}{ \ga}- \ta }
 \lt\|   \lt(v_x,  {v}/{x} \rt)(\cdot,t)\rt\|_{L^\iy}^2 \notag\\
  &
 +(1+t)^{\frac{2\ga-1}{\ga}-\ta}\lt(  \left\|\lt(x\bar\rho^{{1}/{2}}v_t,v,xv_x\rt) (\cdot, t) \right\|^2 +   \left\|\bar\rho^{{\ga}/{2}}\left(r-x, xr_x-x \right)(\cdot, t) \right\|^2 \rt)\notag\\
&+ (1+t)^{\frac{3(\ga-1)}{\ga}-{\ta}}\lt\|(r-x)(\cdot,t)\rt\|^2
+(1+t)^{\frac{\ga-1}{\ga}-\ta}\left\|\left(r_x-1, {r}/{x}-1, v_x, {v}/{x}, \bar\rho^{{1}/{2}}v_t\right)(\cdot,   t) \right\|^2 \notag\\
& + \left\|\bar\rho^{\frac{\ga\ta}{4}-\frac{\ga-1}{2}}\left(r-x, xr_x-x\right)(\cdot,   t) \right\|^2
  \le C(\theta)\mathfrak{E}(0)
\end{align}
and
\begin{align}\label{estthm1'}
&(1+t)^{\frac{\ga-1}{\ga}-\ta}\ \lt\|\lt(r_x-1,  {r}/{x}-1,  v_x,  {v}/{x}\rt)(\cdot,t)\rt\|_{H^1  \lt(\lt[0, \delta\rt]\rt)}^2
  \le C(\theta, \delta) \mathfrak{E}(0).
\end{align}

ii) Suppose that $\mathfrak{F}_\alpha(0)<\iy$ for some $\alpha\in (0, \ga)$. Let $\theta$ and $\delta$ be any constants satisfying
$0<\theta<  \min\lt\{{2(\ga-1)}/({3\ga}), \ \  2(\ga-\alpha)/\ga \rt\}$   and  $\delta\in (0, 1)$.
Then there exist positive constants $C(\alpha)$, $C( \alpha, \theta)$ and $C(\alpha, \theta, \delta)$ such that for all $t\ge 0$,
\begin{align}
& \mathfrak{F}_\alpha(t) \le C(\alpha) \mathfrak{F}_\alpha(0),   \label{} \\
&(1+t)^{({\ga-1})/{ \ga} - \ta} \lt\|\bar\rho^{ ({3\ga-2-\alpha})/4}  (r_x-1, r/x-1)(\cdot,t)\rt\|_{L^\iy}^2 \le C(\alpha,\ta)\mathfrak{F}_\alpha (0);\label{}
\end{align}
and if $\alpha\in [\ga-1, \ga)$,
\be\label{estthm2}\begin{split}
&(1+t)^{ \frac{\kappa }{2}+\frac{4\ga-3}{2\ga}-\frac{3}{2}\ta }   \lt\|\lt(   v_x, v/{x} , \bar\rho^{{1}/{2}} v_t \rt)(\cdot,t)\rt\|^2   + (1+t)^{\frac{8\ga-5}{4\ga}+\frac{\ka}{4}-\frac{5}{4}\ta } \lt\|v(\cdot,t)\rt\|_{L^\iy}^2   \\
 & +  (1+t)^{\frac{1}{2} b_1 }
\|xv_x(\cdot,t)\|^2_{L^\iy} +  (1+t)^{\frac{1}{2}\min\{b_1,  b_2\}  } \lt\| \lt( v_x, v/x\rt) (\cdot,t)\rt\|_{L^\iy}^2   \\
& + (1+t)^{\frac{\ka}{2}+ \frac{2\ga-1}{2\ga}-\frac{\ta}{2}} \lt\|\bar\rho^{ \ga/2 }  (r_x-1, r/x-1)(\cdot,t)\rt\|_{L^\iy}^2 \le C(  \alpha, \ta)\widetilde{\mathfrak{F}}_\alpha(0),
\end{split}\ee
\be\label{estthm2'}\begin{split}
& (1+t)^{\min\lt\{\frac{\kappa }{2}+\frac{4\ga-3}{2\ga}-\frac{3}{2}\ta, \ \frac{\kappa }{4}+\frac{10\ga-9}{4\ga}-\frac{9}{4}\ta\rt\}} \lt\|\lt(r_x-1, {r}/{x}-1,  v_x, v/{x}\rt)(\cdot,t)\rt\|^2_{H^1([0,\da])} \\
& + (1+t)^{ \frac{\kappa }{2}+\frac{4\ga-3}{2\ga}-\frac{3}{2}\ta }   \lt\|\lt(  r_x-1, r/x-1 \rt)(\cdot,t)\rt\|_{L^2([0,\da])}^2   \le C(  \alpha, \ta, \da)\widetilde{\mathfrak{F}}_\alpha(0).
\end{split}\ee
Here $\widetilde{\mathfrak{F}}_\alpha(0)= {\mathfrak{F}}_\alpha(0) +
  {\mathfrak{E}}(0) {\mathfrak{F}}_\alpha(0)$,  $\kappa=0$ when $\alpha = \ga-1$,  $\kappa = (1/\ga) \min\{  {\alpha-(\ga-1)} , \   {\ga-1} \} -\theta$ when $\alpha\in (\ga-1, \ga)$,
and $b_1$ and $b_2$ are given by \ef{newb1} and \ef{newb2}.

iii) Suppose that $\mathfrak{F}_\alpha(0)<\iy$ for some $\alpha\in [\ga, 2\ga-1]$. Let $\theta$  be any constant  satisfying $0<\theta<  {2(\ga-1)}/({3\ga})$.
Then there exist positive constants $C$ and $C(   \theta)$  such that for all $t\ge 0$,
\begin{align}
&  \mathfrak{F}_\alpha (t)
\le   C \mathfrak{F}_\alpha (0)  +  C(\ta) (1+t)^{(\alpha-\ga + \ta\ga)/(\alpha-1)}\mathfrak{E}(0) ;
\end{align}
and if $\alpha=2\ga-1$,
\begin{align}
& \lt\|  r_{xx}(\cdot, t)\rt\|^2 \le C \mathfrak{F}_{2\ga-1} (0)  +  C (\ta) (1+t)^{\frac{1}{2}+\frac{\ga}{2\ga-2}\ta} \mathfrak{E}(0), \label{finite}\\
&\lt\|\lt( v_{xx}, \ (v/x)_x \rt)(\cdot, t)\rt\|^2 \le C(\theta) (1+t)^{-\frac{7\ga-6}{4\ga} + 4 \ta }  \mathfrak{F}_{2\ga-1} (0) \lt(1+ \mathfrak{F}_{2\ga-1} (0)   \rt)  \lt(1+ \mathfrak{E} (0)  \rt)  . \label{decayof2ndderivative}
\end{align}

iv)  Suppose that $\|x\bar\rho^{ {1}/{2}} v_{tt}(\cdot, 0)\|^2<\infty  $.  Then there exists a positive constant  $C$  independent of $t$ such that for all $t\ge 0$,
\be\label{furthregularity}
\|x\bar\rho^{ {1}/{2}} v_{tt}(\cdot, t)\|^2+\int_0^{\infty}  \lt\|(v_{ss}, x v_{ssx})(\cdot,s)\rt\|^2 ds\le C \mathfrak{E}(0)+ C\|x\bar\rho^{ {1}/{2}} v_{tt}(\cdot, 0)\|^2.  \ee
\end{thm}

\subsection{Main ideas and the structure of the proof}\label{sec3.2}
The local existence and uniqueness of strong solutions to \ef{419} on some time interval $[0, T_*]$   are given in Appendix, Part I, by using a finite difference method as used in \cite{Okada,LiXY,LXY,Chengq}. In order to prove the global existence of  strong solutions, we need to derive the uniform-in-time boundedness of the nonlinear functional $\mathfrak{E}(t)$ defined in \ef{mathmarch}.  Our basic strategy is to use the weighted $L^2$-energy method together with some pointwise estimates.
For the weighted $L^2$-energy estimates, motivated by  the linearized analysis for problem  \ef{419} shown in Appendix, Part II, a natural functional should be $\mathcal{E}(t)$  given by
\begin{align}\label{mathcalE}
\mathcal{E}(t):=&\left\|(r-x, xr_x-x)(\cdot,t)  \right\|^2   +  \left\| \left(v, xv_x \right)(\cdot, t) \right\|^2 +\lt\|(r_x-1)(\cdot,t)\rt\|_{L^\iy\lt([1/2,\ 1]\rt)}^2  \notag\\
&
+  \lt\|   \bar\rho^{\ga- {1}/{2}} (r_{xx}, (r/x)_x)(\cdot, t) \rt\|^2  + \lt\|   \bar\rho^{1/2} v_t(\cdot, t) \rt\|^2.
\end{align}
(Indeed, $\mathfrak{E}(t)$ and $\mathcal{E}(t)$ are equivalent under some assumptions which will be shown   in Lemma \ref{boundsforrv}.)
However, much efforts are needed to pass from the linear to  nonlinear analysis. An important step for this is to identify the a priori bounds with which the basic bootstrap argument can work. It is found that the appropriate a priori assumption is that $|r_x(x, t)-1|$ and $|v_x|$ are suitably small to pass from the linear  to nonlinear analysis. To this end, we use a bootstrap argument by making the following  {\it a priori} assumptions.  Let $v$ be a strong solution to \ef{419} on $[0 , T]$  with
 $$r(x, t)=r_0(x)+\int_0^t v(x, \tau)d\tau, \ \ (x,t)\in [0, \ 1]\times[0,T].$$
The basic {\it a priori}  assumption  is that there exist suitably small fixed constants $\ea_0\in (0, {1}/{2}]$ and $\ea_1\in (0, 1]$ such that
\begin{equation}\label{aprirx}
\left|r_x(x,t)-1\right|  \le \ea_0 \ \  {\rm and} \  \
\lt|v_x(x,t)\rt|   \le \ea_1 \ \ {\rm for}  \ \ (x,t)\in I\times [0, T].
\end{equation}
It follows from  \ef{aprirx}  and the boundary condition $v(0,t)=0$ (so $r(0,t)=0$)  that \begin{equation}\label{rx}
 \left|x^{-1}r(x,t) -1 \right| \le  \left|r_x(x,t)-1\right|\le \ea_0 \ \ {\rm for}
 \ \  (x,t)\in [0, \ 1]\times[0,T],
\end{equation}
\begin{equation}\label{vx}
\lt|x^{-1} v(x,t) \rt| \le \lt|v_x(x,t)\rt|   \le  \ea_1 \ \ {\rm for}
 \ \ (x,t)\in [0, \ 1]\times[0,T].
\end{equation}
In particular, it holds that
\begin{equation}\label{Liy}
   {1}/{2}\le r_x(x,t)\le  {3}/{2}   \  \ {\rm and} \ \   {1}/{2}\le x^{-1}{r}(x,t)  \le  {3}/{2} \ \ {\rm for}
 \ \ (x,t)\in I\times [0, T].
\end{equation}

In order to close the argument, it is needed to bound the $L^{\infty}$-norms of $r_x-1$ and $v_x$ by the initial data.
The usual approach for this is to use energy estimates and  the Sobolev embedding, for example, for $r_x-1$,
$$\|(r_x -1)(\cdot,t)\|^2_{L^{\infty} }\le C \lt(\| (r_x -1)(\cdot, t)\|^2 +\|r_{xx}(\cdot, t)\|^2 \rt).$$
However, due to the strong degeneracy of the equation near the vacuum boundary,  the $L^2$-norm of $r_{xx}$ may grow in time (see \ef{finite}), the uniform $L^2$-bound for $r_{xx}$ is valid only for an interval of $x$ away from the vacuum boundary $x=1$ (see \ef{estthm1'}), say,  $\|r_{xx}\|_{L^2([0, 1/2])}$. The term $\|\bar\rho^{\ga- {1}/{2}}r_{xx}(\cdot, t)\|$ in the functional $\mathcal{E}(t)$ can only give  the bound of $|r_x(x, t)-1|$ for $x$ away from the vacuum boundary. This is the reason why we include $\|(r_x  -1)(\cdot, t)\|^2_{L^{\infty}([1/2,1])}$ in the the functional $\mathcal{E}(t)$.
Similar ideas apply to the estimate of  $\|v_x(\cdot, t)\|_{L^{\infty}}$, which is routinely bounded by $\|v_x(\cdot, t)\|_{H^1}$ via energy methods in the region away from the vacuum boundary $x=1$,  say $[0, {1}/{2}]$,  since the system is not degenerate there. However, in the region near the vacuum boundary, say, $[{1}/{2}, 1]$, it is not a easy task to obtain the $L^2$-bound  of the second derivative of $v$,  $\|v_{xx}(\cdot, t)\|_{L^2([{1}/{2}, 1])}$, again, due to the strong degeneracy of the system near the vacuum.  Indeed, our idea to bound $v_x$ near the vacuum boundary is to use the $L^2$-norm of the viscosity term to obtain the $L^2$-bound for $(v_x/r_x)_x$, instead of that for $v_{xx}$. Because $r_x$ has positive lower and upper bounds.

The $L^2$-norm of the viscosity term  is quite different from that of $v_{xx}$. This can be seen as follows:
$$\mathcal{V}:=\left(\frac{ (r^2v)_x}{r^2  r_x}\right)_{x}=\lt(\frac{v_x}{r_x}+2 \frac{v}{r}\rt)_x=\frac{1}{r_x} v_{xx}-\frac{1}{r_x^2} v_x r_{xx} +2\lt(\frac{v}{r}\rt)_x.$$
Under the {\it a priori} assumption that $r_x$ is close to $1$, the difference between the viscosity term $\mathcal{V}$ and $v_{xx}$ involves the second
derivative of $r$, $r_{xx}$, whose  $L^2$-norm may grow in time. This growth may be balanced by the decay of the $L^{\infty}$-norm of
$v_x$ so that the $L^2$-norm of $r_x^{-2}{v_xr_{xx}}  $ is bounded. It is possible to make an ansatz of the decay of the $L^{\infty}$-norm of
$v_x$ and then close the argument.  However it is quite intricate to do so. On the level of the second derivatives, our strategy is to close the estimates by using the $L^2$-estimate of the viscosity term $\mathcal{V}$, instead of that of $v_{xx}$. It should be emphasized that it is enough to get the $L^2$-estimate for the viscosity term $\mathcal{V}$  to close the estimates and obtain the necessary bounds for the global existence of strong solutions and assurance of the well-definiteness of the boundary condition, as explained in Section \ref{sec2.2}.  After proving the global existence of strong solutions and basic decay estimates, we can improve further regularity that $ v_{xx} \in L^2(I)$ for all $t\ge 0$ and the decay of $\| v_{xx}   (\cdot, t)\|$ if the initial data satisfy additional regularity that  $r_{xx}(\cdot, 0) \in L^2(I)$ (i.e.,  $\mathfrak{F}_{2\ga-1}(0)<\iy$).

Due to equation \ef{viscosityequation}, the equivalent form of equation \ef{419a},   the $L^2$-norm of the viscosity term $\mathcal{V}$ can be bounded by $C\mathcal{E}(t)$. (More precisely, $\|\bar\rho^{-1/2}\mathcal{V}(\cdot, t)\| \le C \mathcal{E}(t)$.)
So, it suffices to show that the higher-order functional $\mathcal{E}(t)$ defined by \ef{mathcalE} is bounded uniformly by the initial data, i.e.,
$$\mathcal{E}(t)\le C \mathcal{E}(0) \ \ {\rm for} \ \ {\rm all} \ \ t\in [0, T].$$
We outline the main steps for this as follows.

{\em Lower-order estimates}. The  key elements in our analysis are the weighted estimates by applying various multipliers to the following equation:
 $$
  \bar\rho\left( \frac{x}{r}\right)^2  v_t   +  \left[    \left(\frac{x^2}{r^2}\frac{\bar\rho}{ r_x}\right)^\ga    \right]_x - \frac{x^4}{r^4}  \left(\bar{\rho}^\ga\right)_x  =      \mathfrak{B}_x + 4\la_1 \left(\frac{v}{r}\right)_x ,
$$
which is equivalent to $\eqref{419a}$. Here $\mathfrak{B}$ is defined in \eqref{bdry1}.   In Lemma \ref{lem1}, we use the multiplier $r^2 v$ to  yield the bound for the basic energy
$$\left\|x\bar\rho^{{1}/{2}} v (\cdot, t) \right\|^2 + \left\|x\bar\rho^{{\ga}/{2}}\left( {r}/{x}-1, r_x-1 \right)(\cdot, t) \right\|^2 +\int_0^t \lt\|(v,xv_x)(\cdot,s)\rt\|^2ds
 . $$
The multiplier $r^3-x^3$ plays an important role in the proof of Lemma \ref{lem2},  in which a bound is obtained for
$$ \left\|x\left( {r}/{x}-1, r_x-1 \right)(\cdot, t) \right\|^2 +\int_0^t \left\|x\bar\rho^{{\ga}/{2}}\left( {r}/{x}-1, r_x-1 \right)(\cdot, s) \right\|^2 ds
 , $$
which refines the weighted estimate of $\|x\bar\rho^{{\gamma}/{2}}\left( {r}/{x}-1, r_x-1\right)\|$
obtained in the basic energy estimates. We also show the  decay estimates for the
basic energy in Lemma \ref{lem2} by establishing a bound for
 $$(1+t)\lt(\left\|x\bar\rho^{{1}/{2}} v (\cdot, t) \right\|^2 + \left\|\bar\rho^{{\ga}/{2}}\left({r}-x, xr_x-x \right)(\cdot, t) \right\|^2\rt) +\int_0^t (1+s) \lt\|(v,xv_x)(\cdot,s)\rt\|^2ds. $$
With those estimates, we are able to bound  $|r_x-1|$ away from the origin in Lemma \ref{lem3}, by noting  the fact that $\mathfrak{B}$ can be written as the time derivative of a function so that one can  integrate equation \ef{nsp1} with respect to both $x$ and $t$ to get the desired estimates, where the monotonicity of the Lane-Emden density plays an important role.
Adopting the multiplier $r^2v_t$, a bound for
$$(1+t)\lt(\left\|x\bar\rho^{{1}/{2}} v_t (\cdot, t) \right\|^2 + \left\| \left(v, xv_x \right)(\cdot, t) \right\|^2\rt) +\int_0^t (1+s) \lt\|(v_t,xv_{tx})(\cdot,s)\rt\|^2ds
  $$
is given in Lemma \ref{lem4} by studying the time differentiated problem of \ef{419}.

Further decay estimates and regularity are given in Lemma \ref{lem51}, which is important to the derivation of the decay for $\|r-x\|_{L^{\infty} }$  in \ef{estthm1}.  This in particular implies the convergence of the evolving boundary $r=R(t)$ to
that of the Lane-Emden stationary solution. Lemma \ref{lem51} also shows that  rates of time decay in various norms depend on $\ga$ indicating the balance between the pressure and self-gravitation, and the dissipation of the viscosity.
These further decay estimates are derived from  the following two multipliers:
\be\label{9.22.2}
\int_0^x \bar\rho^{-\beta}(y)(r^3-y^3)_ydy \ \  {\rm  and} \ \  \int_0^x \bar\rho^{-\beta}(y)(r^2 v)_ydy \ \ {\rm for} \ \  0<\beta<\gamma-1.
\ee
It should be noted that the first multiplier in \ef{9.22.2} is motivated by virial equations in the study of stellar dynamics and equilibria (cf. \cite{lebovitz2,tokusky}). To the best of our knowledge, those multipliers have not been used in previous literatures.

{\em Higher-order estimates}.   To get the $L^2$-estimate of the viscosity term $\mathcal{V}$, we write it as follows
$$ \mathcal{V} = \mathcal{G}_{xt},  \ \
{\rm where} \ \  \mathcal{G} : =  \ln r_x + 2 \ln \lt(\frac{r}{x}\rt),  $$
and rewrite $\eqref{419a}$ as, by virtue of \ef{rhox},
 \be\label{viscosityequation}\mu \mathcal{G}_{xt} +\ga \left(\frac{x^2 \bar{\rho}}{r^2 r_x } \right)^{\ga } \mathcal{G}_x
    =  \frac{x^2}{r^2} \bar{\rho}  v_t - \left[ \left(\frac{x^2 }{r^2 r_x } \right)^{\ga }  -\left(\frac{x }{r } \right)^{4}  \right]  x \phi \bar\rho   ,  \   {\rm where} \    \phi(x)= x^{-3}\int_0^x 4\pi \bar\rho(s)  s^2 ds.\ee
This form has advantage in estimating the second derivative (in $x$)  estimates here.  Indeed, the interplay among  the viscosity,  pressure and gravitational force can be seen easily. The gradient of the pressure is decomposed as follows:
 $$\lt[ \bar\rho^{\ga}\left(\frac{x^2 }{r^2 r_x } \right)^{\ga }\rt]_x=-\left(\frac{x^2 }{r^2 r_x } \right)^{\ga }x \phi \bar\rho-\ga \left(\frac{x^2 \bar{\rho}}{r^2 r_x } \right)^{\ga } \mathcal{G}_x.$$
 The first part of this decomposition  is used to balance the gravitational force  and  the second part  is the $t$-antiderivative of the viscosity, $\mathcal{G}_{xt}$, multiplied by a weight which is equivalent to $\bar\rho^{\ga} $. This weight is degenerate on the boundary, but strictly positive in the interior. The degeneracy near the vacuum boundary  is one of main obstacles in   higher-order estimates,  which is overcome by choosing suitable weights and multipliers as we will outline below, and a delicate use of the Hardy and weighted Sobolev inequalities.  (The method of using Hardy and weighted Sobolev inequalities to build regularity was  first adopted for the physical vacuum problem for inviscid flows in \cite{10}, and latter on in \cite{10',jm,17'}.) Indeed, in terms of $\mathcal{G}$,   the principal part of \ef{viscosityequation} is
$$\mu \mathcal{G}_{xt} +\ga \left(\frac{x^2}{r^2 r_x } \right)^{\ga }  \bar{\rho}^\ga \mathcal{G}_x,$$
which is linear in  $\mathcal{G}_x$ and with  a degenerate damping.  This structure leads to desirable estimates on $\mathcal{G}$ and their derivatives.

With the   lower-order estimates obtained already, we can derive in Lemma \ref{lem5} the uniform bound for
$\left\| \bar\rho^{\ga- {1}/{2}} (r_{xx}, (r/x)_x ) (\cdot, t) \rt\|$ (due to the bound for $\left\| \bar\rho^{\ga- {1}/{2}}  \mathcal{G}_x (\cdot, t) \rt\|$ and \ef{weightrxx})
and the decay estimates for
$ \left\| \bar\rho^{ {1}/{2}} v_t(\cdot, t) \rt\|$.
 This completes the proof of the uniform-in-time bounds for the higher-order energy functional $\mathcal{E}(t)$,  which also verifies the
{\it a priori} assumptions \ef{rx} and \ef{vx}  due to the equivalence of $\mathcal{E}(t)$ and $\mathfrak{E}(t)$ shown in Lemma \ref{boundsforrv}, and consequently, the global existence of the strong solution  is obtained.

With the  decay estimates for the lower-order norms in Lemma \ref{lem51} and the higher-order estimates in Lemma \ref{lem5},  we prove the decay estimates of
$\lt\|\lt(r_x-1, \  {r}/{x}-1\rt)(\cdot,t)\rt\| $,  $\lt\|(v, \ v_x)(\cdot,t)\rt\|_{L^\iy}$  and $ \lt\|(r-x)(\cdot,t)\rt\|_{L^\iy}$
in Lemma \ref{lem312}, with which part $i)$ of Theorem \ref{mainthm} is proved.

The second part of the higher-order estimates will be given in Section \ref{sec3.4.3}, in which the faster decay estimates are given under the assumption of the finiteness of $\mathfrak{F}_\alpha(0)$,  $\alpha\in (0, \ga)$.
A key  ingredient in the proof is  to use the new multiplier $x^2\bar\rho^{2\ga-2} \mathfrak{P}_t$, where
$$\mathfrak{P}(x,t)=\ga \left(\frac{x^2 \bar{\rho}}{r^2 r_x } \right)^{\ga } \mathcal{G}_x + \left[ \left(\frac{x^2 }{r^2 r_x } \right)^{\ga }  -\left(\frac{x }{r } \right)^{4}  \right]  x \phi \bar\rho,$$
the sum of  the gradient of the pressure and gravitation force.
The proof of $ii)$ of Theorem \ref{mainthm} consists of Lemmas \ref{lem8.19}-\ref{lem44}.

Parts $iii)$ and  $iv)$ of Theorem \ref{mainthm} on the further regularity of solutions  are proved in Section \ref{sec3.4.4}. In particular, the bound for $\|r_{xx}(\cdot,t)\|$ and the decay for $\|v_{xx}(\cdot,t)\|$ are given in Lemma \ref{lem10}.

\begin{rmk} The insight for the above estimates and the basic decay rates can be gained from the linearized analysis for problem \ef{419}, which is shown in Appendix, Part II.  The linearized analysis gives the ideas for the basic decay rates. The further decay rates of various norms in Theorem \ref{mainthm} (based on Lemma \ref{lem51}) are obtained by applying various multipliers,
reflecting the fact  that the asymptotic stability mechanism is due to the balance between the pressure and self-gravitation, and the dissipation of the viscosity, and also the decay rates depend on the behavior of the initial data near the vacuum boundary.  In particular, for the multipliers in \ef{9.22.2}, we choose $\beta<\ga-1$ so that
$\int \bar\rho^{-\beta}dx<\infty$ due to the physical vacuum behavior that $\bar\rho(x)\sim (1 -x)^{{1}/({\ga-1})}$  near the vacuum boundary. This has an effect on the decay rates of various norms given in  Lemma \ref{lem51} which depend on $\gamma$, and indicates that the physical vacuum behavior may be the reason for the dependence of decay rates on $\ga$.
\end{rmk}

In the rest of this article, we will frequently use the following weighted Sobolev embedding and general version of the Hardy inequality, whose proof can be  found  in \cite{KM}, to build the regularity. This idea was first used in \cite{10} for inviscid flows, and latter on in \cite{10',jm,17'}.

\begin{lem}\label{lemebedding} {\rm (weighted Sobolev embedding)} Let $d$ denote the distance function to the boundary $\pl I$. Then the weighted  Sobolev space $H^1_d(I)$, given by
$$ H^1_d(I) := \lt\{d F\in L^2(I): \ \  \int_I  d^2 \lt( |F|^2 + |F_x|^2 \rt)dx<\infty\rt\},$$
satisfies the  embedding:
\be\label{sobolev}
H^1_d(I)\hookrightarrow L^2(I).
\ee
\end{lem}
\begin{lem}\label{hardy} {\rm (Hardy  inequality)} Let $k>1$ be a given real number and $g$ be a function satisfying
$$
\int_0^{1/2} x^k\lt(g^2 + g_x^2\rt) dx < \iy,
$$
then it holds that
\be\label{hardyorigin}
\int_0^{1/2} x^{k-2} g^2 dx \le c \int_0^{1/2} x^k \lt( g^2 + g_x^2 \rt)  dx,
\ee
where $c$ is a generic constant independent of $g$.
\end{lem}
As a consequence of Lemma \ref{hardy}, one has
\be\label{hardybdry}
\int_{1/2}^{1} (1-x)^{k-2} g^2 dx \le c \int_{1/2}^{1} (1-x)^k \lt( g^2 + g_x^2 \rt)  dx,
\ee
provided that the right-hand side is finite.

\subsection{Lower-order estimates}\label{sec3.3}
In this and next subsections, we derive the {\it a priori } estimates for the strong solution $(r,v)$ on the time interval $[0, T]$ defined in Definition \ref{definitionss},
under the assumption \ef{rx} and \ef{vx}.  We start with   the lower-order estimates in this subsection, for which we rewrite  equation \eqref{419a} as
 \begin{equation}\label{nsp1}\begin{split}
& \bar\rho\left( \frac{x}{r}\right)^2  v_t   +  \left[    \left(\frac{x^2}{r^2}\frac{\bar\rho}{ r_x}\right)^\ga    \right]_x - \frac{x^4}{r^4}  \left(\bar{\rho}^\ga\right)_x  =      \mathfrak{B}_x + 4\la_1 \left(\frac{v}{r}\right)_x .
\end{split}
\end{equation}
Here $\mathfrak{B}$ is defined in \eqref{bdry1}. Recall the boundary conditions \eqref{419b} and \eqref{lex1}:
\be\label{Aug7bdry}
\bar\rho(1)=0, \ \  v(0,t)=0 \ \ {\rm and} \ \ \mathfrak{B}(1,t)=0.
\ee

First, we estimate the basic energy, for which the condition $\gamma>4/3$ is crucial.

\begin{lem}\label{lem1} Suppose that \ef{rx} holds for a suitably small positive number $\epsilon_0$. Then,
 \begin{equation}\label{lem1est}\begin{split}
  &  \left\|x\bar\rho^{\frac{1}{2}} v (\cdot, t) \right\|^2 + (3\ga -4) \left\|x\bar\rho^{\frac{\ga}{2}}\left(\frac{r}{x}-1, r_x-1 \right)(\cdot, t) \right\|^2
    +  \sa \int_0^t  \left\|\left(v, xv_x \right) (\cdot,s) \right\|^2  ds  \\
    \le &c  \lt(
   \left\|x\bar\rho^{\frac{1}{2}} v(\cdot, 0) \right\|^2  + \left\|x \bar\rho^{\frac{\ga}{2}}\left(\frac{r_0}{x}-1, r_{0x}-1\right)  \right\|^2 \rt), \ \ \ \  0\le t\le T,
\end{split}
\end{equation}
where $\sa=\min\left\{2\la_1/3,\ \la_2 \right\}$.
\end{lem}
{\em Proof}. Multiplying equation \eqref{nsp1} by $r^2 v$ and integrating the product with respect to spatial variable, we have, using the integration by parts and   boundary condition \ef{Aug7bdry}, that
 \begin{equation}\label{hz1}\begin{split}
 \frac{d}{dt}\int \tilde{\eta}(x,t) dx  = - \int \mathfrak{B}  \left(r^2 v\right)_x dx  + 4 \la_1  \int r^2 v \left(\frac{v}{r}\right)_x dx ,
\end{split}
\end{equation}
where
\begin{equation}\label{hheta}
\tilde{\eta}(x,t):=\frac{1}{2} x^2 \bar{\rho} v^2  + x^2\bar{\rho}^\ga\left[\frac{1}{\ga-1}\left(\frac{x}{r}\right)^{2\ga-2}\left(\frac{1}{r_x}\right)^{\ga-1}
+\left(\frac{x}{r}\right)^{2}r_x - 4 \frac{x}{r}  \right].
\end{equation}
By the Taylor expansion, the quantity $[\cdot]$ in $\tilde\eta$ can be rewritten as
 $$
 \frac{4-3\ga}{\ga-1} + (2-\ga)\left( \frac{r}{x} - r_x  \right)^2  +  \frac{3\ga-4}{2}\left[ 2\left(\frac{r}{x}-1\right)^2  +  \left(r_x-1\right)^2 \right] +\widetilde{Q} ,
$$
where $\widetilde{Q}$ represents the cubic terms which can be bounded by
$$ |\widetilde{Q}|\le   c \left( \left|r_x-1\right|^3 +  \left| \frac{r}{x}-1\right|^3 \right)
\le c \ea_0 \left( \left|r_x-1\right|^2 +  \left| \frac{r}{x}-1\right|^2   \right),$$
due to \ef{Liy} and \ef{rx}.
This implies that for $\ga\in ( {4}/{3} ,2]$,
$$
   \frac{1}{\ga-1}\left(\frac{x}{r}\right)^{2\ga-2}\left(\frac{1}{r_x}\right)^{\ga-1}
+\left(\frac{x}{r}\right)^{2}r_x - 4 \frac{x}{r}
\ge   \frac{4-3\ga}{\ga-1} +  \frac{3\ga-4}{4} \left[2\left(\frac{r}{x}-1\right)^2  + \left(r_x-1\right)^2 \right],
$$
provided that $\ea_0$ is less than a constant depending on $3\ga-4$. Set
\begin{equation}\label{etadefn}
{\eta}(x,t):=\tilde{\eta}(x,t)- \frac{4-3\ga}{\ga-1}x^2\bar{\rho}^\ga.
\end{equation}
Then the above calculations imply that
\begin{equation}\label{etalower}
{\eta}(x,t) \ge \frac{1}{2} x^2 \bar{\rho} v^2 + \frac{3\ga-4}{4} x^2\bar{\rho}^\ga\left[2\left(\frac{r}{x}-1\right)^2  + \left(r_x-1\right)^2 \right],
\end{equation}
\begin{equation}\label{etaup}
{\eta}(x,t) \le \frac{1}{2} x^2 \bar{\rho} v^2 + c x^2\bar{\rho}^\ga \left[\left(\frac{r}{x}-1\right)^2  + \left(r_x-1\right)^2 \right].
\end{equation}
Clearly,  \eqref{hz1} and \ef{bdry1} show that
\begin{equation}\label{hz0}\begin{split}
 \frac{d}{dt}\int  {\eta}(x,t) dx  = &-\frac{4}{3}\la_1   \int \frac{r^4}{r_x}\left|\left(\frac{v}{r}\right)_x\right|^2  dx
  -\la_2    \int \frac{1}{r_xr^2} \left| \left(r^2 v\right)_x \right|^2dx.
\end{split}
\end{equation}
Note that
$$
 \frac{r^4}{r_x}\left|\left(\frac{v}{r}\right)_x\right|^2  = \frac{r^2}{r_x}v_x^2+r_x v^2 - 2rvv_x \ \ {\rm and} \ \
 \ \ \frac{\left| \left(r^2 v\right)_x \right|^2}{r_xr^2} =\frac{r^2}{r_x}v_x^2+4r_x v^2 +4 rvv_x.
$$
We  obtain
\begin{equation}\label{heg1}\begin{split}
\frac{d}{dt} \int  {\eta}(x,t) dx \le - 3\sa   \int \left[ \frac{r^2}{r_x}v_x^2+ 2r_x v^2  \right]dx,
\end{split}
\end{equation}
where $\sa=\min\left\{2\la_1/3,\ \la_2 \right\}$; and
 \begin{equation}\label{eg1}\begin{split}
 \int  {\eta}(x,t) dx + 3\sa  \int_0^t\int \left[ \frac{r^2}{r_x}v_x^2+ 2r_x v^2  \right]dxds \le  \int  {\eta}(x,0) dx, \ \ t\in [0,T].
\end{split}
\end{equation}
This, together with \ef{etalower}, \ef{etaup} and \ef{Liy}, implies \ef{lem1est}.

\hfill $\Box$

In the following lemma,   we use the multiplier $r^3-x^3$ which is motivated by  virial equations in the study of  stellar dynamics and equilibria (cf. \cite{lebovitz2,tokusky})  to  refine the weighted estimate of $\|x\bar\rho^{{\gamma}/{2}}\left({r}/{x}-1, r_x-1\right)\|$ obtained in Lemma \ref{lem1} by improving the estimates near the vacuum, and give  the decay estimates for the basic energy.
\begin{lem}\label{lem2} Suppose that \ef{rx} holds for a suitably small positive number $\epsilon_0$. Then,
\begin{equation}\label{lem2est}\begin{split}
 &\sa \left\|x\left(\frac{r}{x}-1, r_x-1\right)(\cdot,   t) \right\|^2 + (3\ga-4)\int_0^t \left\|x\bar\rho^{\frac{\ga}{2}}\left(\frac{r}{x}-1, r_x-1\right)(\cdot, s) \right\|^2 ds \\ \le  &   C \lt(\left\|x\left(\frac{r_0}{x}-1, r_{0x}-1\right)  \right\|^2 +  \left\|x\bar\rho^{\frac{1}{2}} v(\cdot, 0)\right\|^2 \rt), \ \ \ \ 0\le t\le T,
\end{split}
\end{equation}
and
\begin{equation}\label{lem2est'}\begin{split}
 &(1+t)   \left\|x\bar\rho^{\frac{1}{2}} v (\cdot, t) \right\|^2 + (3\ga-4) (1+t)  \left\|x\bar\rho^{\frac{\ga}{2}}\left(\frac{r}{x}-1, r_x-1 \right)(\cdot, t) \right\|^2 \\
 & + (1+t)^{\frac{2\ga-2}{\ga}} \left\|(r-x)(\cdot, t) \right\|^2
    + \sa \int_0^t   (1+s) \left\|\left(v, xv_x \right) (\cdot,s) \right\|^2   ds  \\
   & \le    C \lt(\left\|x\left(\frac{r_0}{x}-1, r_{0x}-1\right)  \right\|^2 +  \left\|x\bar\rho^{\frac{1}{2}} v(\cdot, 0)\right\|^2 \rt),  \ \ \ \ 0\le t\le T,
\end{split}
\end{equation}
where $\sa=\min\left\{2\la_1/3,\ \la_2 \right\}$.
\end{lem}
{\em Proof}. The proof consists of two steps. With the basic energy estimate obtained in the previous lemma, we can achieve the estimate for $\|x(r_x-1, {r}/{x}-1\|$ by a moment argument in Step 1. It should be pointed out the double integral obtained in Step 1 will play a crucial role in the derivation of  higher-order estimates later.  In Step 2, we show
the time decay estimates for the
basic energy.

{\em Step 1}.  Multiplying \eqref{nsp1} by $r^3-x^3$ and integrating the resulting equation with respect to the spatial variable, we have, with the help of the integration by parts and  boundary condition  \eqref{Aug7bdry}, that
 \begin{equation*}\label{decAug4}\begin{split}
&\int     \bar{\rho}^\ga \left\{       \left[\frac{x^4}{r^4} \left(r^3-x^3\right)\right]_x  - \left(\frac{x^2}{r^2 r_x} \right)^\ga    \left(r^3-x^3\right)_x   \right\} dx \\
= &    -   \int  \left[ \mathfrak{B} \left(r^3-x^3\right)_x  - 4\la_1 \left(\frac{v}{r}\right)_x \left(r^3-x^3\right) \right] dx  -\int x^3 \bar\rho   v_t \left( \frac{r}{x}-\frac{x^2}{r^2}\right)dx .
\end{split}
\end{equation*}
Notice that
\begin{equation*}\label{}\begin{split}
 &\mathfrak{B} \left(r^3-x^3\right)_x  - 4\la_1 \left(\frac{v}{r}\right)_x \left(r^3-x^3\right)
 =  3\mathfrak{B} \left(r^2 r_x-x^2\right) - 4\la_1\left(\frac{v_x}{r}-\frac{vr_x}{r^2}\right)\left(r^3-x^3\right) \\
=& 4\la_1 x^2\left(\frac{v}{r}-\frac{v_x}{r_x}+\frac{xv_x}{r}
-\frac{xvr_x}{r^2}\right) + 3\la_2 x^2 \left( \frac{r^2}{x^2} v_x + 2 \frac{r v}{x^2 }r_x - \frac{v_x}{r_x} -2\frac{v}{r}\right)\\
=&4\la_1 x^2\left[\ln\left(\frac{r}{x r_x}\right) +\frac{xr_x}{r}
 -1\right]_t + 3\la_2 x^2 \left[ \frac{r^2}{x^2}r_x - \ln \left(\frac{r^2}{x^2}r_x \right) -1\right]_t;
\end{split}
\end{equation*}
which implies
$$
  \int  \left[ \mathfrak{B} \left(r^3-x^3\right)_x  - 4\la_1 \left(\frac{v}{r}\right)_x \left(r^3-x^3\right) \right] dx
   =   \frac{d}{dt} \int x^2\left[ 4\la_1  \Phi_1\left(\frac{r}{x r_x}\right) + 3\la_2 \Phi_2 \left(\frac{r^2}{x^2}r_x \right) \right]dx,
$$
where
\be\label{phi12}\Phi_1(z):=\ln z + z^{-1} -1 \ \ {\rm and} \ \  \Phi_2(z):=z-\ln z -1.\ee
Notice also that
$$
\int x^3 \bar\rho   v_t \left( \frac{r}{x}-\frac{x^2}{r^2}\right)dx =\frac{d}{dt}\int x^3 \bar\rho   v  \left( \frac{r}{x}-\frac{x^2}{r^2}\right)dx - \int x^2 \bar\rho   v^2  \left( 1+2 \frac{x^3}{r^3}\right)dx.
$$
Then, we set
\be\label{toto6}
\eta_0:=  x^2 \tilde{\eta}_0
+x^3 \bar\rho   v  \left( \frac{r}{x}-\frac{x^2}{r^2}\right)
\ \  {\rm and} \ \  \tilde{\eta}_0 := 4\la_1  \Phi_1\left(\frac{r}{x r_x}\right) + 3\la_2 \Phi_2 \left(\frac{r^2}{x^2}r_x \right),
\ee
and obtain
\begin{align}\label{toto7}
&\frac{d}{dt}\int \eta_0 (x, t)dx+\int     \bar{\rho}^\ga \left\{       \left[\frac{x^4}{r^4} \left(r^3-x^3\right)\right]_x  - \left(\frac{x^2}{r^2 r_x} \right)^\ga    \left(r^3-x^3\right)_x   \right\} dx \notag\\
&= \int x^2 \bar\rho   v^2  \left( 1+2 \frac{x^3}{r^3}\right)dx.
\end{align}

Noting that the quantity $\{\cdot\}$ on the left-hand side of \eqref{toto7} can be rewritten as
$$
  x^{2}\lt[3\left(\frac{x^2}{r^2r_x}\right)^{\ga}   -3\left(\frac{x^2}{r^2r_x}\right)^{ \ga-1 } -   \left(\frac{x}{r}\right)^2r_x
  +4\left(\frac{x}{r}\right)^5 r_x - 7\left(\frac{x}{r}\right)^4 + 4 \frac{x}{r}\rt],
$$
we can then show, using a similar way as to the derivation of \eqref{etalower}, that
\begin{align}\label{dec1}
&\int     \bar{\rho}^\ga \left\{       \left[\frac{x^4}{r^4} \left(r^3-x^3\right)\right]_x  - \left(\frac{x^2}{r^2 r_x} \right)^\ga    \left(r^3-x^3\right)_x   \right\} dx \notag \\
&\ge \frac{3(3\ga-4)}{2} \int x^2\bar{\rho}^\ga\left[2\left(\frac{r}{x}-1\right)^2  + \left(r_x-1\right)^2 \right] dx ,
\end{align}
when \ef{rx} holds for a small $\epsilon_0$.
It   follows from \eqref{toto7} and \ef{dec1}  that
\be\label{bye}
  \frac{d}{dt} \int \eta_0(x,t) dx
  + \frac{3(3\ga-4)}{4} \int x^2\bar{\rho}^\ga\left[2\left(\frac{r}{x}-1\right)^2  + \left(r_x-1\right)^2 \right] dx
  \le   C  \int    v^2  dx.
\ee
This implies, with the aid of \eqref{eg1} and \ef{Liy},  that
\begin{equation}\label{newegrx}\begin{split}
 &\int \eta_0(x,t) dx
  + \frac{3(3\ga-4)}{4}\int_0^t \int x^2\bar{\rho}^\ga\left[2\left(\frac{r}{x}-1\right)^2  + \left(r_x-1\right)^2 \right] dxds \\
  \le & C\int (\eta+\eta_0)(x,0) dx.
\end{split}
\end{equation}

It remains to analyze $\tilde{\eta}_0$. By the Taylor expansion and \ef{rx}, one may get that
\begin{equation*}\label{eta0lower}\begin{split}
 \tilde{\eta}_0 \ge  \frac{1}{4} \left[ 4\la_1   \left(\frac{r}{x r_x}-1\right)^2 + 3\la_2  \left(\frac{r^2}{x^2}r_x -1 \right)^2 \right]
 \ge \frac{3}{4}  \sa \left[2  \left(\frac{r}{x r_x}-1\right)^2 +  \left(\frac{r^2}{x^2}r_x -1 \right)^2 \right],
\end{split}\end{equation*}
where $\sa=\min\left\{2\la_1/3,\ \la_2 \right\}$; and
\begin{equation}\label{eta0up}\begin{split}
 \tilde{\eta}_0\le  2 \left[ 4\la_1   \left(\frac{r}{x r_x}-1\right)^2 + 3\la_2  \left(\frac{r^2}{x^2}r_x -1 \right)^2 \right]\le C \lt[\left(\frac{r}{x}-1\right)^2 +(r_x-1)^2 \rt];
\end{split}\end{equation}
provided that $\ea_0$ in \ef{rx} is suitably small.  Notice that
$$
\left(\frac{r}{x}-1\right)^2 \le  \left(\frac{r^3}{x^3}-1\right)^2 =   \left(\frac{r}{x r_x} \frac{r^2}{x^2}r_x -1\right)^2
\le C  \left[\left(\frac{r}{x r_x}-1\right)^2 +  \left(\frac{r^2}{x^2}r_x -1 \right)^2 \right] \le C \sa^{-1} \tilde \eta_0
$$
and also
$$
\left(r_x-1\right)^2 \le C  \tilde \eta_0 + C\left(\frac{r}{x}-1\right)^2  \le C \sa^{-1} \tilde \eta_0.
$$
We then achieve, with the help of \ef{newegrx} and \ef{eg1},  that
\begin{equation}\label{we}\begin{split}
 \sa \int x^2\lt[\left(\frac{r}{x}-1\right)^2 +  \left(r_x-1\right)^2\rt] dx
  + (3\ga-4)\int_0^t \int x^2\bar{\rho}^\ga\left[ \left(\frac{r}{x}-1\right)^2  \rt. \\
  \lt. + \left(r_x-1\right)^2 \right] dxds
 \le  C\int (\eta+\eta_0)(x,0) dx
   .
\end{split}
\end{equation}
This, together with \ef{etaup} and \ef{eta0up}, implies \ef{lem2est}.

{\em Step 2}. We are ready to show the time decay of the basic energy. Let $\eta$ be given by \ef{etadefn}. It follows from \eqref{heg1} that
$$
  (1+t)  \int \eta(x,t)dx + 3 \sa  \int_0^t (1+s) \int    \left(\frac{r^2}{r_x} v_{x}^2 + 2 r_x v^2  \right) dx ds
 \le   \int \eta(x,0)dx +   \int_0^t   \int \eta(x,s)dx ds.
$$
In view of \eqref{etaup}, \eqref{lem1est} and \eqref{lem2est}, one has that
\begin{equation*}\label{}\begin{split}
  \int_0^t   \int \eta(x,s)dx ds
  \le  & C \int_0^t \int v^2 dxds + C\int_0^t   \int  x^2\bar{\rho}^\ga \left[\left(\frac{r}{x}-1\right)^2  + \left(r_x-1\right)^2 \right]dxds  \\
  \le  &   C \int (\eta_0 + \eta)(x,0)dx.
\end{split}
\end{equation*}
So, it holds that
\begin{equation}\label{5-0}\begin{split}
(1+t)  \int \eta(x,t)dx + 3 \sa  \int_0^t (1+s) \int    \left(\frac{r^2}{r_x} v_{x}^2 + 2 r_x v^2  \right) dx ds
  \le  C \int (\eta_0 + \eta)(x,0)dx.
\end{split}
\end{equation}
This, together with \ef{etaup} and \ef{eta0up}, implies
\begin{equation*}\label{}\begin{split}
 &(1+t)   \left\|x\bar\rho^{\frac{1}{2}} v (\cdot, t) \right\|^2 + (3\ga-4) (1+t)  \left\|x\bar\rho^{\frac{\ga}{2}}\left(\frac{r}{x}-1, r_x-1 \right)(\cdot, t) \right\|^2
   \\
  & + \sa \int_0^t   (1+s) \left\|\left(v, xv_x \right) (\cdot,s) \right\|^2   ds  \le    C \lt(\left\|x\left(\frac{r_0}{x}-1, r_{0x}-1\right)  \right\|^2 +  \left\|x\bar\rho^{\frac{1}{2}} v(\cdot, 0)\right\|^2 \rt).
\end{split}
\end{equation*}
Since $x\bar\rho^{\gamma-1}$ is equivalent to the distance function, $dist(x, \partial I)$, it then follows from the Sobolev embedding \ef{sobolev}, \ef{phy} and the H${\rm \ddot{o}}$lder inequality that
\begin{align}\label{rminusx}
  &\int(r-x)^2(x,t)dx \le \int x^2 \bar\rho^{2(\ga-1)}\left((r-x)^2+( r_x-1)^2 \right)(x,t)dx\notag \\
   \le & \lt(\int x^2 \left((r-x)^2+( r_x-1)^2  \right)(x, t) dx\rt)^{\frac{2-\ga}{\ga}}
 \lt( \int x^2 \bar\rho^{\ga}\left((r-x)^2+( r_x-1)^2 \right)dx \rt)^{\frac{2\ga-2}{\ga}}\notag\\
 \le & C(1+t)^{-\frac{2\ga-2}{\ga}}\lt(\left\|x\left(\frac{r_0}{x}-1, r_{0x}-1\right)  \right\|^2 +  \left\|x\bar\rho^{\frac{1}{2}} v(\cdot, 0)\right\|^2 \rt).
\end{align}
This finishes the proof of
\ef{lem2est'}.

\hfill $\Box$

With the estimates obtained so far, we are able to derive a pointwise bound for $\lt|r/{x}-1\rt|$ and $\lt|r_x-1\rt|$ away from the origin, by realizing that      equation \ef{nsp1} can be integrated with respect to both $x$ and $t$. It should be noted that  the monotonicity of  the Lane-Emden density which decreases  in the  radial direction outward plays an important role for this estimate.

\begin{lem}\label{lem3} Let $I_2=[1/2,1]$. For a suitably small constant $\ea_0$ in \ef{rx}, it holds that for $(x,t)\in I_2\times [0,T]$,
 \begin{equation}\label{lem3est}\begin{split}
   & \lt|x^{-1} {r(x,t)} -1\rt|+\lt|r_x(x,t)-1\rt| \\
   \le &  C \lt(\left\|x\bar\rho^{\frac{1}{2}} v(\cdot, 0) \right\|
   +  \left\|x\left(\frac{r_0}{x}-1, r_{0x}-1\right)  \right\|+ \lt\|r_{0x}-1\rt\|_{L^{\infty}\lt( I_2\rt)}\rt).
\end{split}
\end{equation}
\end{lem}
{\em Proof}. The proof consists two steps.

{\em Step 1  (bound  for $ {r}/{x}-1$)}.
 Notice that
$$
x(r-x)^2=\int_0^x \lt[y(r(y,t)-y)^2\rt]_y dy
\le \lt\|r-x\rt\|^2 + 2 \lt\|r-x\rt\|\lt\|x(r_x-1)\rt\|.
$$
This, together with  \ef{lem2est}, yields that for $x\in I_2$,
 \begin{equation}\label{we1}\begin{split}
\lt|\frac{r}{x}-1\rt|^2  \le 8 x(r-x)^2
\le & C\lt(\left\|x\bar\rho^{\frac{1}{2}} v(\cdot, 0) \right\|^2
   +  \left\|x\left(\frac{r_0}{x}-1, r_{0x}-1\right)  \right\|^2\rt).
\end{split}
\end{equation}

{\em Step  2  (bound for $r_x-1$)}. Integrating equation \ef{nsp1} over $[x, 1]$ and using the boundary condition \eqref{Aug7bdry}, one gets
\begin{equation}\label{ec1}\begin{split}
& \int_x^1 \bar\rho\left( \frac{y}{r}\right)^2  v_t dy  -   \left(\frac{x^2}{r^2}\frac{\bar\rho}{ r_x}\right)^\ga   - \int_x^1 \frac{y^4}{r^4}  \left(\bar{\rho}^\ga\right)_y dy =     - \mathfrak{B}+ 4\la_1 \int_x^1 \left(\frac{v}{r}\right)_ydy ;
\end{split}
\end{equation}
where $\mathfrak{B}$,  defined by \eqref{bdry1}, can be rewritten as
$$\mathfrak{B}= \mu\lt( \ln r_x\rt)_t -\lt(\frac{4}{3}\la_1-2\la_2\rt) \lt( \ln r\rt)_t.$$
So,  \ef{ec1} is equivalent to
\begin{equation*}\label{}\begin{split}
 \mu \lt( \ln r_x\rt)_t = & \left(\frac{x^2}{r^2}\frac{\bar\rho}{ r_x}\right)^\ga  - \lt( \int_x^1 \bar\rho \frac{y^2}{r^2}  v  dy \rt)_t -2  \int_x^1 \bar\rho\frac{ y^2 }{r^3}  v^2  dy + \int_x^1 \frac{y^4}{r^4}  \left(\bar{\rho}^\ga\right)_y dy \\
 &+ \lt(\frac{4}{3}\la_1-2\la_2\rt) \lt( \ln r\rt)_t+ 4\la_1 \int_x^1 \lt( \ln r\rt)_{yt} dy .
\end{split}
\end{equation*}
Integrate it with respect to the temporal variable to obtain
$$
\mu\ln \lt(\frac{r_x}{r_{0x}}\rt)  =  \int_0^t \left(\frac{x^2}{r^2}\frac{\bar\rho}{ r_x}\right)^\ga ds + \mathfrak{L}    ,
$$
where
\begin{equation*}\label{}\begin{split}
\mathfrak{L} =& - \lt. \int_x^1 \bar\rho \frac{y^2}{r^2}  v  dy \rt|_0^t -2 \int_0^t \int_x^1 \bar\rho\frac{ y^2 }{r^3}  v^2  dyds  + \int_0^t \int_x^1 \frac{y^4}{r^4}  \left(\bar{\rho}^\ga\right)_y dy ds \\
& + \lt(\frac{4}{3}\la_1-2\la_2\rt)   \ln \lt(\frac{r}{r_0}\rt)+ 4\la_1 \ln \lt(\frac{r(1,t) }{r_0(1)}\frac{r_0(x) }{r(x,t)}\rt)  ;
\end{split}
\end{equation*}
which implies that
\begin{equation}\label{ec2}\begin{split}
 r_x  = r_{0x} \exp\lt\{ \frac{1}{\mu} \mathfrak{A} \rt\} \exp\lt\{ \frac{1}{\mu}  \mathfrak{L} \rt\}  , \ \ {\rm where} \ \ \mathfrak{A}= \int_0^t \left(\frac{x^2}{r^2}\frac{\bar\rho}{ r_x}\right)^\ga ds .
\end{split}
\end{equation}
On the other hand, direct calculations show, by virtue of \ef{ec2},  that
$$\mathfrak{A}_t=\left(\frac{x^2}{r^2}\frac{\bar\rho}{ r_x}\right)^\ga
=\left(\frac{x^2}{r^2} \frac{\bar\rho}{r_{0x}} \right)^\ga  \exp\lt\{ - \frac{\ga}{\mu} \mathfrak{A}\rt\} \exp\lt\{ -\frac{\ga}{\mu}  \mathfrak{L} \rt\},$$
so that
$$ \exp\lt\{  \frac{\ga}{\mu} \mathfrak{A}\rt\}
=1+ \int_0^t \frac{\ga}{\mu}\left(\frac{x^2}{r^2} \frac{\bar\rho}{r_{0x}} \right)^\ga   \exp\lt\{ -\frac{\ga}{\mu}  \mathfrak{L} \rt\} d\tau.$$
It then follows from \ef{ec2} that
\begin{equation}\label{ec3}\begin{split}
 r_x =&r_{0x}  \lt[1+ \int_0^t \frac{\ga}{\mu}\left(\frac{x^2}{r^2} \frac{\bar\rho}{r_{0x}} \right)^\ga   \exp\lt\{ -\frac{\ga}{\mu}  \mathfrak{L}_1  \rt\} \exp\lt\{ -\frac{\ga}{\mu}  \int_0^\tau \int_x^1 \frac{y^4}{r^4}  \left(\bar{\rho}^\ga\right)_y dy ds \rt\}d\tau\rt]^{1/\ga} \\
  & \times \exp\lt\{ \frac{1}{\mu}  \mathfrak{L}_1 \rt\}\exp\lt\{ \frac{1}{\mu}  \int_0^t \int_x^1 \frac{y^4}{r^4}  \left(\bar{\rho}^\ga\right)_y dy ds \rt\},
\end{split}
\end{equation}
where
$$\mathfrak{L}_1=\mathfrak{L}- \int_0^t \int_x^1 \frac{y^4}{r^4}  \left(\bar{\rho}^\ga\right)_y dy ds. $$
In view of \ef{lem1est} and \ef{we1}, one can get  that for $x\ge 1/2$,
\begin{equation*}\label{}\begin{split}
\lt|\mathfrak{L}_1\rt|
\le &C   \lt(\int x^2  \bar\rho   v^2  dx \int \bar\rho dx\rt)^{1/2} + C   \lt(\int x^2  \bar\rho u^2_0(r_0(x))  dx \int \bar\rho dx\rt)^{1/2}\\
&+C   \int_0^t \int    v^2  dyds
 +C \lt\|\frac{r}{x}-1\rt\|_{L^\iy\lt(I_2\times[0,T]\rt)}
 \le C\tilde{\mathfrak{e}}
\end{split}
\end{equation*}
where
$$\tilde{\mathfrak{e}}= \left\|x\bar\rho^{\frac{1}{2}} v(\cdot, 0) \right\|
   +  \left\|x\left(\frac{r_0}{x}-1, r_{0x}-1\right)  \right\| .$$
It therefore follows from \ef{ec3} and \ef{we1}  that
\begin{equation*}\label{}\begin{split}
 r_x \le &r_{0x}  \lt[1+ (1+ C\mathfrak{e} ) \int_0^t \frac{\ga}{\mu}{\bar\rho}^\ga    \exp\lt\{ -\frac{\ga}{\mu}  \int_0^\tau \int_x^1 \frac{y^4}{r^4}  \left(\bar{\rho}^\ga\right)_y dy ds \rt\}d\tau\rt]^{1/\ga} \\
  &  \times(1+ C\mathfrak{e} ) \exp\lt\{ \frac{1}{\mu}  \int_0^t \int_x^1 \frac{y^4}{r^4}  \left(\bar{\rho}^\ga\right)_y dy ds \rt\}\\
\le &r_{0x}\lt(1+C\mathfrak{e} \rt) \lt[\exp\lt\{ \frac{\ga}{\mu}  \int_0^t \int_x^1 \frac{y^4}{r^4}  \left(\bar{\rho}^\ga\right)_y dy ds \rt\} \rt.\\
      &\lt.+ \lt(1+C\mathfrak{e} \rt)\int_0^t \frac{\ga}{\mu}\bar\rho^\ga  \exp\lt\{ \frac{\ga}{\mu}  \int_\tau^t \int_x^1 \frac{y^4}{r^4}  \left(\bar{\rho}^\ga\right)_y dy ds \rt\} d\tau\rt]^{1/\ga},
\end{split}\end{equation*}
where
$\mathfrak{e}= \tilde{\mathfrak{e}}+\lt\|r_{0x}-1\rt\|_{L^{\infty}(I_2)}.$
Observe that $\left(\bar{\rho}^\ga\right)_x<0$. So, one can derive from \ef{we1} that
\begin{equation*}\label{}\begin{split}
 r_x \le &r_{0x}\lt(1+C\mathfrak{e} \rt) \lt[\exp\lt\{ \frac{\ga}{\mu}  \int_0^t \int_x^1 \lt(1-C\mathfrak{e} \rt) \left(\bar{\rho}^\ga\right)_y dy ds \rt\} \rt.\\
      &\lt.+ \lt(1+C\mathfrak{e} \rt)\int_0^t \frac{\ga}{\mu}\bar\rho^\ga  \exp\lt\{ \frac{\ga}{\mu}  \int_\tau^t \int_x^1 \lt(1-C\mathfrak{e} \rt)  \left(\bar{\rho}^\ga\right)_y dy ds \rt\} d\tau\rt]^{1/\ga}\\
  \le&  r_{0x}\lt(1+C\mathfrak{e} \rt) \lt[\exp\lt\{- \frac{\ga}{\mu} (1-C\mathfrak{e})  \bar{\rho}^\ga t \rt\} \rt.\\
      &\lt.+ \lt(1+C\mathfrak{e} \rt)\int_0^t \frac{\ga}{\mu}\bar\rho^\ga  \exp\lt\{- \frac{\ga}{\mu} (1-C\mathfrak{e})  \bar{\rho}^\ga (t-\tau) \rt\} d\tau\rt]^{1/\ga} \\
\le&  r_{0x}\lt(1+C\mathfrak{e} \rt) \lt[\exp\lt\{- \frac{\ga}{\mu} (1-C\mathfrak{e})  \bar{\rho}^\ga t \rt\} + \frac{1+C\mathfrak{e}}{1-C\mathfrak{e} }\lt. \exp\lt\{- \frac{\ga}{\mu} (1-C\mathfrak{e})  \bar{\rho}^\ga (t-\tau) \rt\}\rt|_{\tau=0   }^t\rt]^{1/\ga} \\
 \le&  r_{0x}\lt(1+C\mathfrak{e} \rt) \lt\{\exp\lt\{- \frac{\ga}{\mu} (1-C\mathfrak{e})  \bar{\rho}^\ga t \rt\} \lt[1-\frac{1+C\mathfrak{e} }{1-C\mathfrak{e}}\rt]+ \frac{1+C\mathfrak{e} }{1-C\mathfrak{e}}\rt\}^{1/\ga}
 \le  r_{0x} \lt( 1+C\mathfrak{e} \rt).
\end{split}
\end{equation*}
  Similarly,
$r_x\ge r_{0x} \lt( 1-C\mathfrak{e} \rt). $
These two estimates, together with \ef{we1}, imply \ef{lem3est}.

\hfill $\Box$

The following lemma gives the decay estimates for the weighted norms of both the time and spatial derivatives of $v$.

\begin{lem}\label{lem4} Let \ef{rx} and \ef{vx} be true. Then it holds that,   for  $0\le t\le T$,
\begin{equation}\label{lem4est}\begin{split}
 &(1+t)    \left\|\lt(x\bar\rho^{\frac{1}{2}} v_t,  v, x v_x \right)(\cdot, t) \right\|^2
    + \int_0^t   (1+s) \left\|\left(v_t, xv_{tx},  \right) (\cdot,s) \right\|^2   ds \\
   & \le    C \lt( \left\|\lt(x\bar\rho^{\frac{1}{2}} v_t,  v, x v_x \right)(\cdot, 0) \right\|^2  + \left\|x\left(\frac{r_0}{x}-1, r_{0x}-1\right)  \right\|^2   \rt).
\end{split}
\end{equation}
\end{lem}
{\em Proof}. Multiplying equation \eqref{nsp1} by $r^2$ and differentiating the resulting equation with respect to $t$, we obtain
\begin{equation}\label{nsp1t}\begin{split}
 &\bar\rho x^2  v_{tt}  -\ga r^2 \left[  \left(\frac{x^2}{r^2}\frac{\bar\rho}{ r_x}\right)^\ga \left(2\frac{v}{r}+\frac{v_x}{r_x}\right)   \right]_{x}
+ 2 rv \left[  \left(\frac{x^2}{r^2}\frac{\bar\rho}{ r_x}\right)^\ga \right]_{x}
+2 \frac{x^4}{r^3} v \left(\bar{\rho}^\ga\right)_x \\
   = & r^2 \left[\mathfrak{B}_{xt} + 4\la_1 \left(\frac{v}{r}\right)_{xt} \right]+ 2rv \left[\mathfrak{B}_{x} + 4\la_1 \left(\frac{v}{r}\right)_{x}\right].
\end{split}
\end{equation}
Set
\begin{equation}\label{eta1}\begin{split}
\eta_1(x,t):=&\frac{1}{2}x^2 \bar{\rho} v_t^2 +\left(\frac{x^2}{r^2}\frac{\bar\rho}{ r_x}\right)^\ga\left[(2\ga-1)r_x v^2
+2(\ga-1)r vv_x +\frac{\ga}{2}\frac{r^2}{r_x} v_x^2
\right] \\
&-\bar\rho^\ga  \lt[\lt(4\frac{x^3}{r^3}-3\frac{x^4}{r^4}r_x\rt) v^2 +  2\frac{x^4}{r^3}vv_x \rt].
\end{split}\end{equation}
Following the estimates for $\eta$ defined in \ef{etadefn}, we can show that,  for $\ga\in (4/3, 2]$,
\begin{equation}\label{eta1lower}
{\eta_1}(x,t) \ge \frac{1}{2} x^2 \bar{\rho} v_t^2 + \frac{3\ga-4}{4} x^2\bar{\rho}^\ga\left[2\left(\frac{v}{x} \right)^2  + v_x^2 \right],
\end{equation}
\begin{equation}\label{eta1up}
{\eta_1}(x,t) \le \frac{1}{2} x^2 \bar{\rho} v_t^2 + c x^2\bar{\rho}^\ga \left[\left(\frac{v}{x} \right)^2  + v_x^2 \right],
\end{equation}
provided that \ef{rx}  holds with $\ea_0$ being suitably small.
Multiplying \ef{nsp1t}  by $  v_t$ and integrating the product with respect to the spatial variable, we have, using the integration by parts and boundary condition \eqref{Aug7bdry}, that
\begin{equation}\label{hhevt}\begin{split}
 &\frac{d}{dt}\int \eta_1(x,t) dx  +\int \left[\mathfrak{B}_t  \left(r^2 v_t\right)_x  - 4 \la_1  r^2 v_t \left(\frac{v}{r}\right)_{xt} \right]dx
= \mathfrak{I}_1+\mathfrak{I}_2  ,
\end{split}
\end{equation}
where
\begin{equation*}\label{df}\begin{split}
 \mathfrak{I}_1:= &
    - \int \left[ \mathfrak{B}  (2rvv_t )_x - 8 \la_1    \left(\frac{v}{r}\right)_{x} rvv_t  \right] dx, \\
\mathfrak{ I}_2:=& (2\ga-1)\int  \left[\left(\frac{x^2}{r^2}\frac{\bar\rho}{ r_x}\right)^\ga r_x \right]_tv^2 dx
+2(\ga-1)\int \left[\left(\frac{x^2}{r^2}\frac{\bar\rho}{ r_x}\right)^\ga r \right]_tvv_x dx \\
&+\frac{\ga}{2}\int \left[\left(\frac{x^2}{r^2}\frac{\bar\rho}{ r_x}\right)^\ga \frac{r^2}{r_x} \right]_t v_x^2
 dx - \int \bar\rho^\ga  \lt[\lt(4\frac{x^3}{r^3}-3\frac{x^4}{r^4}r_x\rt)_t v^2 +  2\lt(\frac{x^4}{r^3}\rt)_tvv_x \rt] dx.
\end{split}
\end{equation*}

 The second term on the left-hand side of \eqref{hhevt} can be estimated as follows. Notice that
 $$
 \mathfrak{B}_t=\frac{4}{3}\la_1\frac{r}{r_x}\left(\frac{v_t}{r}\right)_x  + \la_2 \frac{\left(r^2 v_t\right)_x}{r_x r^2} + \bar{\mathfrak{B}},
 $$
where
 \begin{equation}\label{barB}\begin{split}
 \bar{\mathfrak{B}}:=&\frac{4}{3}\la_1 \left[\left(\frac{v}{r}\right)^2-\left(\frac{v_x}{r_x}\right)^2 \right]
 -\la_2  \left[2\left(\frac{v}{r}\right)^2+\left(\frac{v_x}{r_x}\right)^2 \right].
 \end{split}
\end{equation}
Thus,
 \begin{equation}\label{bfevt4}\begin{split}
 &\int \left[\mathfrak{B}_t  \left(r^2 v_t\right)_x  - 4 \la_1  r^2 v_t \left(\frac{v}{r}\right)_{xt} \right]dx\\
=&   \int  \left\{ \frac{4}{3}\la_1\left(\frac{v_t}{r}\right)_x \left[\frac{r}{r_x}\left(r^2 v_t\right)_x- 3 r^2 v_t\right] +\la_2\frac{ \left| \left(r^2 v_t\right)_x \right|^2}{r_xr^2}\right\}dx  -\mathfrak{ I}_3\\
\ge &  3\sa   \int \left[ \frac{r^2}{r_x}v_{tx}^2+ 2r_x v_t^2  \right]dx -\mathfrak{I}_3,
   \end{split}
\end{equation}
where $\sa=\min\left\{2\la_1/3, \  \la_2 \right\}$  and
 $$
 \mathfrak{I}_3:= -  \int  \bar{\mathfrak{B}}\left(r^2 v_t\right)_x dx - 4\la_1  \int r^2 v_t \left(\frac{v^2}{r^2}\right)_x dx .  $$
So, \eqref{hhevt} implies that
$$
 \frac{d}{dt}\int \eta_1(x,t) dx   + 3\sa \int \left( \frac{r^2}{r_x}v_{tx}^2+ 2r_x v_t^2  \right)dx
\le  \mathfrak{ I}_1 +\mathfrak{I}_2 + \mathfrak{I}_3.
$$

For $\mathfrak{ I}_1$ and $\mathfrak{ I}_3$, it follows from \eqref{Liy}, \eqref{vx} and the Cauchy inequality that
 $$
\mathfrak{ I}_1+\mathfrak{ I}_3 \le
\sa   \int \left( \frac{r^2}{r_x}v_{tx}^2+ 2r_x v_t^2  \right)dx  +  C\sa^{-1} \ea_1^2 \int    \left(x^2 v_{x}^2 +  v^2  \right) dx.
$$
Similarly, $\mathfrak{ I}_2$ can be bounded by
$$
\mathfrak{ I}_2 \le   C  \ea_1    \int    \left(x^2 v_{x}^2 + v^2  \right) dx.
$$
So, we arrive at the following estimate
\begin{equation}\label{later1}\begin{split}
 \frac{d}{dt}\int \eta_1(x,t) dx   + 2\sa \int \left( \frac{r^2}{r_x}v_{tx}^2+ 2r_x v_t^2  \right)dx
\le  C\int    \left(x^2 v_{x}^2 + v^2  \right) dx ,
\end{split}
\end{equation}
provided that \ef{vx} holds for $\ea_1\le 1$.  This, together with \eqref{eg1}, implies that
\begin{equation}\label{eg2}\begin{split}
  \int \eta_1(x,t) dx +\sa\int_0^t \int \left( x^2 v_{sx}^2+  v_s^2  \right)dxds
\le  \int \eta_1(x,0) dx  + C\int \eta(x,0)dx
\end{split}
\end{equation}
and
\begin{equation*}\label{5-1}\begin{split}
&  (1+t)  \int \eta_1(x,t)dx +  \sa  \int_0^t (1+s) \int    \left(x^2 v_{sx}^2 + v_s^2  \right) dx ds\\
 \le & \int \eta_1(x,0)dx +   \int_0^t   \int \eta_1(x,s)dx ds
 + C\int_0^t (1+s)\int    \left(x^2 v_{x}^2 + v^2  \right) dxds\\
\le & \int \eta_1(x,0)dx +  C \int_0^t   \int v_s^2 dx ds
 + C\int_0^t (1+s) \int    \left(x^2 v_{x}^2 + v^2  \right) dxds.
\end{split}
\end{equation*}
Here \ef{eta1up} has been used.
This, together with \ef{lem2est'}, \ef{eg2} and \eqref{eta1up}, implies
\begin{equation}\label{tttt}\begin{split}
 &(1+t)  \lt(   \left\|x\bar\rho^{\frac{1}{2}} v_t (\cdot, t) \right\|^2 +  \left\| \bar\rho^{\frac{\ga}{2}}\left(v, xv_x \right)(\cdot, t) \right\|^2\rt)
    + \int_0^t   (1+s) \left\|\left(v_s, xv_{sx},  \right) (\cdot,s) \right\|^2   ds \\
   & \le    C \lt(\left\|x\bar\rho^{\frac{1}{2}} (v_t,v) (\cdot, 0) \right\|^2 +  \left\| \bar\rho^{\frac{\ga}{2}}\left(v, xv_x \right)(\cdot, 0) \right\|^2 + \left\|x\left(\frac{r_0}{x}-1, r_{0x}-1\right)  \right\|^2   \rt).
\end{split}
\end{equation}

 Observe that
$$
(1+t)v^2(x,t)
\le v^2(x,0)+\int_0^t (1+s)\lt[2v^2(x,s) +  v_s^2(x,s) \rt] ds.
$$
Integrate the above inequality with respect to the spatial variable to give
$$
(1+t)\int v^2(x,t)dx
\le \int v^2(x,0)dx +\int_0^t (1+s) \int \lt[2v^2(x,s) +  v_s^2(x,s) \rt]dx ds.
$$
Similarly, it holds that
$$
(1+t)\int x^2 v_x^2(x,t)dx
\le \int x^2 v^2(x,0)dx +\int_0^t (1+s) \int \lt[2x^2 v^2(x,s) +  x^2 v_s^2(x,s) \rt]dx ds.
$$
This, together with \ef{lem2est'} and \ef{tttt}, implies that
\begin{equation}\label{jump}\begin{split}
 &(1+t)  \lt(   \left\|x\bar\rho^{\frac{1}{2}} v_t (\cdot, t) \right\|^2 +  \left\| \left(v, xv_x \right)(\cdot, t) \right\|^2\rt)
    + \int_0^t   (1+s) \left\|\left(v_s, xv_{sx},  \right) (\cdot,s) \right\|^2   ds \\
   & \le    C \lt(\left\|x\bar\rho^{\frac{1}{2}}v_t (\cdot, 0) \right\|^2 +  \left\|  \left(v, xv_x \right)(\cdot, 0) \right\|^2 + \left\|x\left(\frac{r_0}{x}-1, r_{0x}-1\right)  \right\|^2   \rt).
\end{split}
\end{equation}
This finishes the proof of \ef{lem4est}.

\hfill $\Box$

Next, we derive  further time decay estimates based on Lemmas \ref{lem1}, \ref{lem2} and \ref{lem4}  by using
two multipliers
$$\int_0^x \bar\rho^{-\beta}(y)(r^3-y^3)_ydy \ \  {\rm  and} \ \  \int_0^x \bar\rho^{-\beta}(y)(r^2 v)_ydy \ \ {\rm for} \ \  0<\beta<\gamma-1.$$
The key is to deal with the behavior of solutions near both the boundary and  geometrical singularity at the origin simultaneously.
The improved decay estimates
obtained in this lemma give the convergence of the evolving boundary $r=R(t)$ to
that of the Lane-Emden stationary solution.

\begin{lem}\label{lem51} Suppose that \ef{rx} and \ef{vx} hold.  Then for any $\theta\in \lt(0,\  {2(\ga-1)}/({3\ga})\rt)$, there exists
a constant $C(\theta)$ independent of $t$ such that
\begin{align}\label{estlem51}
&\left\|\bar\rho^{\frac{\ga\ta}{4}-\frac{\ga-1}{2}}\left(r-x, xr_x-x\right)(\cdot,   t) \right\|^2
+(1+t)^{\frac{\ga-1}{\ga}-\ta}\left\|\left(xr_x-x\right)(\cdot,   t) \right\|^2 +(1+t)^{\frac{3(\ga-1)}{\ga}-\theta} \notag\\
& \times\left\|(r-x)(\cdot, t) \right\|^2
+(1+t)^{\frac{2\ga-1}{\ga}-\ta}\lt(  \left\|\lt(x\bar\rho^{\frac{1}{2}}v_t,v,xv_x\rt) (\cdot, t) \right\|^2 +   \left\|\bar\rho^{\frac{\ga}{2}}\left(r-x, xr_x-x \right)(\cdot, t) \right\|^2 \rt)
 \notag  \\
&+ \int_0^t \lt[\left\|\bar\rho^{\frac{\ta\ga+2}{4}}\left(r-x, xr_x-x \right)(\cdot, s) \right\|^2 +  (1+s)^{\frac{\ga-1}{\ga}-\ta}\left\|\bar\rho^{\frac{\ga}{2}}\left(r-x, xr_x-x \right)(\cdot, s) \right\|^2\rt] ds\notag\\
  & + \int_0^t \lt[(1+s)^{\frac{2\ga-1}{\ga}-\ta} \left\|\left(v, xv_x, v_s, x v_{sx} \right) (\cdot,s) \right\|^2  +  (1+s)^{\frac{2\ga-1}{2\ga}-\frac{\ta}{2}} \left\|\bar\rho^{\frac{\ga\ta}{4}-\frac{\ga-1}{2}}\left(v, xv_x \right) (\cdot,s) \right\|^2 \rt]  ds \notag\\
 & \le C(\theta)\lt( \left\|\lt(  v, x  v_x ,  x\bar\rho^{\frac{1}{2}} v_t \rt)(\cdot, 0) \right\|^2   +  \left\| r_{0x}-1 \right\|_{L^\iy}^2\rt) , \  \ \  \ t\in [0, T].
\end{align}
Moreover, we have for any $a\in (0,1)$ and $\theta\in \lt(0,\  {2(\ga-1)}/({3\ga})\rt)$,
\begin{align}
(1+t)^{2(\ga-1)/\ga-\ta}\|(r-x)(\cdot,t)\|_{L^\iy([a,1])}^2
+(1+t)^{(2\ga-1)/\ga-\ta}\|v(\cdot,t)\|_{L^\iy([a,1])}^2 \notag\\
\le C(a,\theta)\lt( \left\|\lt(  v, x  v_x ,  x\bar\rho^{{1}/{2}} v_t \rt)(\cdot, 0) \right\|^2   +  \left\| r_{0x}-1 \right\|_{L^\iy}^2\rt) , \  \ \  \ t\in [0, T]. \label{8/12-1}
\end{align}
\end{lem}
{\em Proof}.  For any given  $\theta\in \lt(0,  \ {2(\ga-1)}/({3\ga})\rt)$,  we set
\begin{equation}\label{alphaiota}
\beta:=\gamma-1-\frac{\gamma\theta}{2}, \ \ \iota:=\frac{\theta}{2}, \  \ \kappa:=\frac{\beta}{\ga}-\iota,  \ \  \nu:=\frac{1}{2}\lt(1+\frac{\beta}{\ga}-\iota\rt),
\end{equation}
so that
 $0<\beta<\gamma-1$  and $  0<\iota< {\beta}/({2\gamma}).$
The proof of  this lemma consists of the following four steps.

{\em Step 1}. In this step, we prove that
\begin{align}\label{key}
&\int \lt[(1+t)^{ \nu}\bar\rho^\ga +1\rt] x^2\bar{\rho}^{-\beta}\left[\left(\frac{r}{x}-1\right)^2  + \left(r_x-1\right)^2 \right] (x,t) dx\notag\\
&+ \int_0^t  (1+s)^{ \nu}   \int \bar\rho^{-\beta} \left(x^2 v_x^2+ v^2  \right)dxds
 +  \int_0^t \int x^2\bar{\rho}^{\ga-\beta} \left[\left(\frac{r}{x}-1\right)^2  + \left(r_x-1\right)^2 \right] dxds\notag\\
& \le C \left\|  r_{0x}-1  \right\|_{L^\iy}^2
  +  C \sum_{i=1}^3 \int_0^t (1+s)^{ \nu}| K_i |ds+ C \sum_{i=1}^3 \int_0^t |L_i| ds,
 \end{align}
where
\begin{equation}\label{Li}\begin{split}
L_1= &-\int \bar\rho \frac{x^2}{r^2}  v_t   \left( \int_0^x \bar\rho^{-\beta}(r^3-y^3)_ydy \right) dx , \\
L_2=& \int     \bar{\rho}^\ga \lt(\frac{x^4}{r^4}\rt)_x \left[\bar\rho^{-\beta} \left(r^3-x^3\right) -  \int_0^x \bar\rho^{-\beta}(r^3-y^3)_ydy  \right] dx, \\
L_3=& 4\la_1    \int  \left(\frac{v}{r}\right)_x\left[\int_0^x \bar\rho^{-\beta}(r^3-y^3)_ydy - \bar\rho^{-\beta}\left(r^3-x^3\right) \right] dx;
\end{split}
\end{equation}
and
\begin{equation}\label{Ki}\begin{split}
K_1  = &-\int \bar\rho\frac{x^2}{r^2}v_t \lt(\int_0^x \bar\rho^{-\beta}(r^2 v)_ydy \rt)dx,\\
K_2= & \int \bar{\rho}^{\ga}\lt(\frac{x^4}{r^4}\rt)_x\lt[\bar{\rho}^{-\beta}{r^2}v
 - \int_0^x \bar\rho^{-\beta}(r^2 v)_ydy\rt]dx,\\
K_3= & 4 \la_1 \int \left(\frac{v}{r}\right)_x \lt[ \lt(\int_0^x \bar\rho^{-\beta}(r^2 v)_ydy \rt)-\bar\rho^{-\beta} r^2 v \rt]dx.
\end{split}
\end{equation}

To this end,   we multiply \eqref{nsp1}  by  $\int_0^x \bar\rho^{-\beta}(y)(r^3-y^3)_ydy$ and integrate the resulting equation with respect to the spatial variable to obtain, with the aid of the integration by parts and  the boundary condition \eqref{Aug7bdry}, that
 \begin{equation*}\label{}\begin{split}
& \int     \bar{\rho}^{\ga-\beta} \left\{       \left[\frac{x^4}{r^4} \left(r^3-x^3\right)\right]_x  - \left(\frac{x^2}{r^2 r_x} \right)^\ga    \left(r^3-x^3\right)_x   \right\} dx
\\
 &+\int \bar\rho^{-\beta} \left[ \mathfrak{B} \left(r^3-x^3\right)_x - 4\la_1 \left(\frac{v}{r}\right)_x \left(r^3-x^3\right) \right] dx
=\sum_{i=1}^3 L_i.
\end{split}
\end{equation*}
Noticing that
\bee\label{}\begin{split}
   &\int x^2\bar{\rho}^{-\beta} \left[\left(\frac{r}{x}-1\right)^2  + \left(r_x-1\right)^2 \right](x,0) dx \\
   &\le C \left\|\lt(r_0-x, \ x r_{0x}-x \rt)  \right\|_{L^\iy}^2 \int_0^1 (1-x)^{-\frac{\beta}{\ga-1}}dx
   \le C\frac{\left\|  r_{0x}-1 \right\|_{L^\iy}^2}{(\ga-1)-\beta}
\end{split}\eee
due to \ef{phy} and $r_0(0)=0$,  one can  obtain, following the derivation of \ef{we}, that
\begin{align}\label{beauty1}
   &\int x^2\bar{\rho}^{-\beta} \left[\left(\frac{r}{x}-1\right)^2  + \left(r_x-1\right)^2 \right](x,t) dx  +  \int_0^t \int x^2\bar{\rho}^{\ga-\beta} \left[\left(\frac{r}{x}-1\right)^2  + \left(r_x-1\right)^2 \right] dxds\notag\\
  & \le C \left\| r_{0x}-1  \right\|_{L^\iy}^2   + C \sum_{i=1}^3 \int_0^t |L_i| ds.
\end{align}
Next, multiplying equation \eqref{nsp1} by $\int_0^x \bar\rho^{-\beta}(y)(r^2 v)_ydy$, integrating the product with respect to spatial variable, and using the integration by parts and  the boundary condition \eqref{Aug7bdry},  one obtains,
$$
 \int \bar{\rho}^{\ga-\beta}\lt[\lt(\frac{x^4}{r^2}v\rt)_x-\left(\frac{x^2}{r^2r_x}\right)^{\ga} \left(r^2 v\right)_x \rt]dx
   +  \int \bar\rho^{-\beta}\lt[\mathfrak{B}  \left(r^2 v\right)_x dx - 4 \la_1  r^2 v \left(\frac{v}{r}\right)_x \rt]dx = \sum_{i=1}^3 K_i ,
$$
Following the derivation of \ef{heg1}, one can then obtain
$$
\frac{d}{dt} \int  {\eta}_2 (x,t) dx + 3\sa   \int \bar\rho^{-\beta} \left[ \frac{r^2}{r_x}v_x^2+ 2r_x v^2  \right]dx \le \sum_{i=1}^3 K_i,
$$
where
$$
{\eta}_2(x,t): =\bar\rho^{-\beta}\lt(\eta(x,t)- \frac{1}{2} x^2 \bar{\rho} v^2 \rt)\approx  x^2\bar{\rho}^{\ga-\beta}\left[\left(\frac{r}{x}-1\right)^2  + \left(r_x-1\right)^2 \right].
$$
Here and thereafter,  $f\approx g$ means that $C^{-1}g\le f\le C g$ with a generic positive constant $C$.
Multiplying the equation above by $(1+t)^{\nu}$ and integrating the product with respect to the temporal variable lead to
\begin{equation}\label{beauty}\begin{split}
&(1+t)^{\nu}\int x^2\bar{\rho}^{\ga-\beta}\left[\left(\frac{r}{x}-1\right)^2  + \left(r_x-1\right)^2 \right] (x,t) dx
+ \int_0^t  (1+s)^{\nu}  \int \bar\rho^{-\beta} \left(x^2 v_x^2+ v^2  \right)dxds\\
& \le C \left\|  r_{0x}-1  \right\|_{L^\iy}^2 +  C \sum_{i=1}^3 \int_0^t (1+s)^{\nu}| K_i |ds
 +  C\int_0^t  \int x^2\bar{\rho}^{\ga-\beta}\left[\left(\frac{r}{x}-1\right)^2  + \left(r_x-1\right)^2 \right] dxds,
\end{split}
\end{equation}
due to the fact ${\beta}/{\ga}<1$.
So,   estimate \ef{key}  follows by a suitable combination of \ef{beauty} and \ef{beauty1}.

{\em Step 2.} In this step, we show that
\be\label{14eye}\begin{split}
&\int \lt[(1+t)^{2\nu}\bar\rho^\ga
+(1+t)^{\kappa} +\bar{\rho}^{-\beta}\rt] x^2\left[\left(\frac{r}{x}-1\right)^2  + \left(r_x-1\right)^2 \right] (x,t) dx\\
&+(1+t)^{2\nu}\int x^2 \bar\rho v^2 (x,t) dx+\int_0^t   \int \lt[(1+s)^{\kappa} + \bar\rho^{-\beta} \rt]x^2\bar{\rho}^{\ga}  \left[\left(\frac{r}{x}-1\right)^2  + \left(r_x-1\right)^2 \right] dxds \\
&
+ \int_0^t \int \lt[(1+s)^{2\nu}+ (1+s)^{\nu}  \bar\rho^{-\beta} \rt] \left(x^2 v_x^2+ v^2  \right)dxds\\
&\le  C \lt( \lt\|x \bar\rho^{\frac{1}{2}}v(\cdot,0)\rt\|^2
+ \left\|  r_{0x}-1  \right\|_{L^\iy}^2\rt)
  +  C \sum_{i=1}^3 \int_0^t (1+s)^{\nu}| K_i |ds+ C \sum_{i=1}^3 \int_0^t |L_i| ds,
\end{split}\ee
where $L_i$ and $K_i$ ( $i=1, 2, 3$) are given by \ef{Li} and \ef{Ki}, respectively.

To prove \ef{14eye}, one can integrate the product of  $(1+t)^{2\nu}$ and \ef{heg1} with respect to the temporal variable to get
\begin{align*}\label{}
&(1+t)^{2\nu}\int \lt\{x^2 \bar\rho v^2 + x^2\bar{\rho}^{\ga}\left[\left(\frac{r}{x}-1\right)^2  + \left(r_x-1\right)^2 \right]\rt\} (x,t) dx\\
& +\int_0^t  (1+s)^{2\nu}   \int \left(x^2 v_x^2+ v^2  \right)dxds
 \le C \lt\|x \bar\rho^{\frac{1}{2}}v(\cdot,0)\rt\|^2 + C\left\|  r_{0x}-1  \right\|_{L^\iy}^2  \\
&+C\int_0^t ({1+s})\int v^2 dxds
+C\int_0^t (1+s)^{\frac{\beta}{\ga}-\iota} \int \ x^2\bar{\rho}^{\ga}\left[\left(\frac{r}{x}-1\right)^2  + \left(r_x-1\right)^2 \right]  dxds,
\end{align*}
since $\beta/\ga<1$.
Integrate the product of  $(1+t)^{\kappa}$ and \ef{bye} with respect to the temporal variable to give
\begin{align*}\label{}
& (1+t)^{\kappa}  \int x^2 \left[\left(\frac{r}{x}-1\right)^2  + \left(r_x-1\right)^2 \right](x,t) dx\\
& +  \int_0^t (1+s)^{\kappa} \int x^2\bar{\rho}^{\ga} \left[\left(\frac{r}{x}-1\right)^2  + \left(r_x-1\right)^2 \right] dxds  \le  C\left\|  r_{0x}-1  \right\|_{L^\iy}^2 \\
& \quad +
  C\int (1+s) \int v^2 dxds  + C \int (1+s)^{\frac{\beta}{\ga}-\iota-1}\int x^2 \left[\left(\frac{r}{x}-1\right)^2  + \left(r_x-1\right)^2 \right](x,t) dxds.
\end{align*}
The last term on the right-hand side of the inequality above is estimated as follows. It follows from the H$\ddot{\rm o}$lder inequality and the Young inequality that
\begin{align*}\label{}
 \int x^2 \left[\left(\frac{r}{x}-1\right)^2  + \left(r_x-1\right)^2 \right]  dx\le
 & \lt(\int x^2 \bar\rho^{\ga-\beta } \left[\left(\frac{r}{x}-1\right)^2  + \left(r_x-1\right)^2 \right] dx\rt)^{\frac{\beta }{\ga}} \\
& \times
\lt(\int x^2 \bar\rho^{-\beta } \left[\left(\frac{r}{x}-1\right)^2  + \left(r_x-1\right)^2 \right] dx\rt)^{\frac{\ga-\beta }{\ga}}
\end{align*}
and
\bee\label{}\begin{split}
&\int_0^t(1+s)^{\beta/\ga-\iota-1}\int x^2 \left[\left(\frac{r}{x}-1\right)^2  + \left(r_x-1\right)^2 \right]  dxds\\
&\le C\int_0^t \int x^2 \bar\rho^{\ga-\beta} \left[\left(\frac{r}{x}-1\right)^2  + \left(r_x-1\right)^2 \right] dxds
\\
&+
C\int_0^t (1+s)^{\lt(\beta/\ga-\iota-1\rt)\frac{\ga}{\ga-\beta}}ds
\sup_{s\in [0,t]}\int x^2 \bar\rho^{-\beta} \left[\left(\frac{r}{x}-1\right)^2  + \left(r_x-1\right)^2 \right](x,s) dx \\
&\le C\frac{\ga-\beta}{\iota \ga}\sup_{s\in [0,t]}\int x^2 \bar\rho^{-\beta} \left[\left(\frac{r}{x}-1\right)^2  + \left(r_x-1\right)^2 \right](x,s) dx\\
&+C\int_0^t \int x^2 \bar\rho^{\ga-\beta} \left[\left(\frac{r}{x}-1\right)^2  + \left(r_x-1\right)^2 \right] dxds.
\end{split}\eee
In a similar way as to deriving \ef{key}, we then have, noting
 \ef{lem2est'}, that
\be\label{bridge}\begin{split}
&(1+t)^{2\nu}\int \lt\{x^2 \bar\rho v^2 + x^2\bar{\rho}^{\ga}\left[\left(\frac{r}{x}-1\right)^2  + \left(r_x-1\right)^2 \right]\rt\} (x,t) dx\\
&+(1+t)^{\kappa}  \int x^2 \left[\left(\frac{r}{x}-1\right)^2  + \left(r_x-1\right)^2 \right](x,t) dx
+\int_0^t  (1+s)^{2\nu}  \int \left(x^2 v_x^2+ v^2  \right)dxds\\ &
+\int_0^t (1+s)^{\kappa}  \int x^2\bar{\rho}^{\ga} \left[\left(\frac{r}{x}-1\right)^2  + \left(r_x-1\right)^2 \right] dxds \\
& \le C \lt\|x \bar\rho^{\frac{1}{2}}v(\cdot,0)\rt\|^2 + C \left\|  r_{0x}-1  \right\|_{L^\iy}^2 + \sup_{s\in [0,t]}\int x^2 \bar\rho^{-\beta} \left[\left(\frac{r}{x}-1\right)^2  + \left(r_x-1\right)^2 \right](x,s) dx\\
&+C\int_0^t \int x^2 \bar\rho^{\ga-\beta} \left[\left(\frac{r}{x}-1\right)^2  + \left(r_x-1\right)^2 \right] dxds.
\end{split}\ee
Make a summation of  $k\times\ef{key}$ and \ef{bridge} with suitable large $k$ to give
\ef{14eye}.

{\em Step 3}.  We claim that
\begin{align}\label{step4}
&\int \lt[(1+t)^{\frac{2\ga-1}{\ga}-\ta}\bar\rho^\ga
+(1+t)^{\frac{\ga-1}{\ga}-\ta} +\bar{\rho}^{-\lt(\ga-1-\frac{1}{2}\ga\ta\rt)}\rt] x^2\left[\left(\frac{r}{x}-1\right)^2  + \left(r_x-1\right)^2 \right] (x,t) dx\notag\\
&+\int_0^t \int \lt[(1+s)^{\frac{2\ga-1}{\ga}- {\ta} }+ (1+s)^{\frac{2\ga-1}{2\ga}-\frac{\ta}{2}}  \bar\rho^{-\lt(\ga-1-\frac{1}{2}\ga\ta\rt)} \rt]\left(x^2 v_x^2+ v^2  \right)dxds\notag\\
&
+\int_0^t   \int \lt[(1+s)^{\frac{\ga-1}{\ga}-\ta} + \bar\rho^{-\lt(\ga-1-\frac{1}{2}\ga\ta\rt)} \rt]x^2\bar{\rho}^{\ga} \left[\left(\frac{r}{x}-1\right)^2  + \left(r_x-1\right)^2 \right] dxds\notag \\
&+(1+t)^{\frac{2\ga-1}{\ga}-\ta}\int x^2 \bar\rho v^2 (x,t) dx
\le  C Q(0),
\end{align}
 where and in the following
$$
\label{}Q(0):=\left\|\lt(  v, x  v_x ,  x\bar\rho^{\frac{1}{2}} v_t \rt)(\cdot, 0) \right\|^2 + \left\|  r_{0x}-1  \right\|_{L^\iy}^2.
$$

To prove this claim, it remains to estimate $K_i$ and $L_i$ in \ef{14eye}.
 First,  it follows from \ef{phy} that for any given constants  $\da\in (0,1]$ and $\beta\in (0, \ga-1)$,
\begin{equation}\label{fact}\begin{split}
\int_{1-\da}^1 \lt(\int_0^x \bar\rho^{-\beta} dy \rt) dx
\le & C \int_{1-\da}^1 \lt(1-\frac{\beta}{\ga-1}\rt)^{-1}\lt[1
 -(1-x)^{1-\frac{\beta}{\ga-1}}\rt]dx\\
 \le & C \frac{\ga-1}{(\ga-1)-\beta} \da= C\lt[(\ga-1)-\beta\rt]^{-1}\da .
\end{split}
\end{equation}
Let $\oa\in(0,1/2)$ be a small  constant to  be determined at the end of this step. It follows from the Cauchy inequality, \ef{lem4est}, the H$\ddot{o}$lder inequality and \ef{fact} that
\begin{equation}\label{14k1}\begin{split}
&\int_0^t (1+s)^{\nu}|K_1|ds  \le \oa \int_0^t (1+s)^{\nu}  \int  \lt(\int_0^x \bar\rho^{-\beta}y (| v| + |y v_y|)dy \rt)^2 dx ds \\
 &  + C\oa^{-1} \int_0^t (1+s) \int v_s^2 dxds
\le  C\oa^{-1} Q(0)+ C \oa \int_0^t (1+s)^{\nu} \int_0^1 \bar\rho^{-\beta} (| v|^2 + |y v_y|^2 )dyds,
\end{split}
\end{equation}
since
\begin{equation*}\label{}\begin{split}
  \int  \lt(\int_0^x \bar\rho^{-\beta}y (| v| + |y v_y|)dy \rt)^2 dx  \le & \int_0^1 \bar\rho^{-\beta} (| v|^2 + |y v_y|^2 )dy \int \lt(\int_0^x \bar\rho^{-\beta} y^2 dy \rt) dx  \\
   \le &  C \int_0^1 \bar\rho^{-\beta} (| v|^2 + |y v_y|^2 )dy  .
\end{split}
\end{equation*}
Similarly, one can obtain
\begin{equation}\label{14l1}\begin{split}
&\int_0^t |L_1|ds
\le   C \oa^{-1} \int_0^t  \int v_s^2 dxds + \oa \int_0^t \int \bar\rho^2 \lt(\int_0^x \bar\rho^{-\beta} y^2 \left(\lt|\frac{r}{y}-1\rt|  + \lt|r_y-1\right| \right) dy \rt)^2 dx ds\\
 & \le   C \oa^{-1} Q(0) + C \oa \int_0^t \int x^2\bar{\rho}^{\ga-\beta} \left[\left(\frac{r}{x}-1\right)^2  + \left(r_x-1\right)^2 \right] dxds ,
\end{split}
\end{equation}
since
$$
 \int \bar\rho^2 \lt( \int_0^x \bar\rho^{-\ga-\beta}dy \rt) dx \le C \int \bar\rho^2  \bar\rho^{-\ga-\beta+(\ga-1)} dx = C \int \bar\rho^{2-\ga}  \bar\rho^{-\beta+(\ga-1)} dx\le C.
$$
 $K_2$ can be rewritten as
\begin{equation*}\label{}\begin{split}
 K_2=& \int_0^{1-\oa} \bar{\rho}^{\ga}\lt(\frac{x^4}{r^4}\rt)_x \int_0^x \lt(\bar\rho^{-\beta}\rt)_y r^2 vdydx \\
 &+ \int_{1-\oa}^1 \bar{\rho}^{\ga}\lt(\frac{x^4}{r^4}\rt)_x\lt[\bar{\rho}^{-\beta}{r^2}v
 - \int_0^x \bar\rho^{-\beta}(r^2 v)_ydy\rt]dx=:K_{21}+K_{22}.
\end{split}
\end{equation*}
Note that
\begin{equation*}\label{}\begin{split}
&|K_{21}| = \lt|\frac{\beta}{\ga} \int_0^{1-\oa} \bar{\rho}^{\ga}\lt(\frac{x^4}{r^4}\rt)_x \int_0^x \bar\rho^{-\beta-(\ga-1)}y \phi r^2 vdydx \rt|\\
\le & C\int_0^{1-\oa} \bar{\rho}^{\ga}x\lt(|r_x-1|+\lt|\frac{r}{x}-1\rt|\rt) x^{-2} \lt( \int_0^x y^6 dy \rt)^{1/2}dx \lt(\int_0^1 v^2 dy\rt)^{1/2}\\
\le  &  C \lt(\int x^2\bar{\rho}^{\ga} \left[\left(\frac{r}{x}-1\right)^2  + \left(r_x-1\right)^2 \right] dx\rt)^{1/2} \lt(\int_0^1 v^2 dy\rt)^{1/2}\\
\le & C \oa^{-1} (1+t)^{-\nu} \int x^2\bar{\rho}^{\ga} \left[\left(\frac{r}{x}-1\right)^2  + \left(r_x-1\right)^2 \right] dx +C \oa (1+t)^\nu \int  v^2 dx
\end{split}
\end{equation*}
and
\begin{equation*}\label{}\begin{split}
&|K_{22}|\le \int_{1-\oa}^1 \bar{\rho}^{\ga}\lt(\frac{x^4}{r^4}\rt)_x\lt[\bar{\rho}^{-\beta}{r^2}v
 - \int_0^x \bar\rho^{-\beta}(r^2 v)_ydy\rt]dx\\
& \le  C \int_{1-\oa}^1 x\bar{\rho}^{\ga-\beta}\lt(|r_x-1|+\lt|\frac{r}{x}-1\rt|\rt)|v|dx \\
& + C \int_{1-\oa}^1 x\bar{\rho}^{\ga}\lt(|r_x-1|+\lt|\frac{r}{x}-1\rt|\rt)\lt(\int_0^x \bar\rho^{-\beta}\lt(|v|+|yv_y|\rt)dy \rt) dx\\
&\le   C \oa^{\frac{\ga-\beta}{2(\ga-1)}} \lt[ (1+t)^{-\nu} \int_{1-\oa}^1 x^2\bar{\rho}^{\ga-\beta}\lt(|r_x-1|^2+\lt|\frac{r}{x}-1\rt|^2\rt)dx
+(1+t)^\nu\int v^2 dx\rt]\\
&+ C \oa^{\frac{\ga}{4(\ga-1)}} \lt[ (1+t)^{-\nu} \int_{1-\oa}^1 x^2\bar{\rho}^{\ga-\beta}\lt(|r_x-1|^2+\lt|\frac{r}{x}-1\rt|^2\rt)dx + (1+t)^\nu \int_0^1  (v^2 + x^2 v_x^2 )dx\rt],
\end{split}
\end{equation*}
due to
\begin{equation*}\label{}\begin{split}
& \int_{1-\oa}^1 \bar\rho^{\frac{\ga}{2}+\beta}\lt(\int_0^x \bar\rho^{-\beta} (| v|^2 + |y v_y|^2 )dy\rt)^2 dx\\
\le & \lt(\int_{1-\oa}^1 \bar\rho^{\frac{\ga}{2}+\beta}\int_0^x \bar\rho^{-2\beta} dydx \rt) \int_0^1  (| v|^2 + |y v_y|^2 )dy \le C \int_0^1  (| v|^2 + |y v_y|^2 )dy.
\end{split}
\end{equation*}
Then,  one  gets, using \ef{lem2est},  that
\begin{equation}\label{14k2}\begin{split}
&\int_0^t (1+s)^\nu|K_2|ds
\le  C\oa^{-1} Q(0) + C \lt(\oa +\oa^{\frac{\ga}{4(\ga-1)}} + \oa^{\frac{\ga-\beta}{2(\ga-1)}} \rt) \\
&\times \lt[  \int_0^t (1+s)^{2\nu} \int_0^1   (v^2 + x^2 v_x^2  )dxds+ \int_0^t   \int x^2\bar{\rho}^{\ga-\beta}\lt(|r_x-1|^2+\lt|\frac{r}{x}-1\rt|^2\rt)dx ds \rt].
\end{split}
\end{equation}
Similarly, $L_2$ can be rewritten as
\begin{equation*}\label{}\begin{split}
L_2 =& \int_0^{1-\oa} \bar{\rho}^{\ga}\lt(\frac{x^4}{r^4}\rt)_x \int_0^x \lt(\bar\rho^{-\beta}\rt)_y(r^3-y^3)dydx \\
 &+ \int_{1-\oa}^1 \bar{\rho}^{\ga}\lt(\frac{x^4}{r^4}\rt)_x\lt[\bar{\rho}^{-\beta}(r^3-x^3)
 - \int_0^x \bar\rho^{-\beta}(r^3-y^3)_ydy\rt]dx=:L_{21}+L_{22}.
\end{split}
\end{equation*}
Note that
\begin{equation*}\label{}\begin{split}
&|L_{21}| = \lt|\frac{\beta}{\ga} \int_0^{1-\oa} \bar{\rho}^{\ga}\lt(\frac{x^4}{r^4}\rt)_x \int_0^x \bar\rho^{-\beta-(\ga-1)}y \phi(r^3-y^3)dydx \rt|\\
\le & C\int_0^{1-\oa} \bar{\rho}^{\ga}x\lt(|r_x-1|+\lt|\frac{r}{x}-1\rt|\rt) x^{-2} \lt( \int_0^x \bar\rho^{-\beta-(\ga-1)}y^2 |r-y|dy\rt)dx\\
\le & C(\omega)   \int x^2\bar{\rho}^{\ga} \left[\left(\frac{r}{x}-1\right)^2  + \left(r_x-1\right)^2 \right] dx
\end{split}
\end{equation*}
and
\begin{equation*}\label{}\begin{split}
|L_{22}| \le & C\int_{1-\oa}^1 x\bar{\rho}^{\ga-\beta}\lt(|r_x-1|+\lt|\frac{r}{x}-1\rt|\rt)|r-x|dx \\
 & +C \int_{1-\oa}^1 x\bar{\rho}^{\ga}\lt(|r_x-1|+\lt|\frac{r}{x}-1\rt|\rt)\lt(\int_0^x \bar\rho^{-\beta}y^2\lt(|r_y-1|+\lt|\frac{r}{y}-1\rt|\rt)dy \rt) dx\\
 \le &  \oa \int x^2\bar{\rho}^{\ga-\beta} \left[\left(\frac{r}{x}-1\right)^2  + \left(r_x-1\right)^2 \right] dx
+  C\oa^{-1}    \int_{1-\oa}^1 \bar\rho^{\ga-\beta} (r-x)^2 dx\\
&+  C   \lt(\int x^2\bar{\rho}^{\ga-\beta} \left[\left(\frac{r}{x}-1\right)^2  + \left(r_x-1\right)^2 \right] dx\rt)^{1/2}\\
 &\times \lt( \int_{1-\oa}^1 \bar{\rho}^{\ga +\beta }\lt(\int_0^x \bar\rho^{-\beta}y^2\lt(|r_y-1|+\lt|\frac{r}{y}-1\rt|\rt)dy   \rt)^2 dx \rt)^{1/2}\\
 \le & C \int \lt(\oa \bar{\rho}^{\ga-\beta} + \oa^{-1}\bar\rho^\ga \rt)x^2 \left[\left(\frac{r}{x}-1\right)^2  + \left(r_x-1\right)^2 \right] dx,
\end{split}
\end{equation*}
where we have used the following simple estimates due to \ef{hardybdry} and \ef{phy}:
\begin{equation*}\label{}\begin{split}
\int_{1-\oa}^1 \bar\rho^{\ga-\beta} (r-x)^2 dx \le & C \int_{1/2}^1 \bar\rho^{\ga-\beta + 2(\ga-1)} \lt[(r-x)^2 +(r_x-1)^2\rt] \\
\le & C \int_{1/2}^1 \bar\rho^{\ga} \lt[(r-x)^2 +x^2 (r_x-1)^2\rt]
\end{split}
\end{equation*}
and
\begin{align*}\label{}
&\int_{1-\oa}^1 \bar{\rho}^{\ga +\beta }\lt(\int_0^x \bar\rho^{-\beta}y^2\lt(|r_y-1|+\lt|\frac{r}{y}-1\rt|\rt)dy   \rt)^2 dx  \\
&\le C\int_{1-\oa}^1\bar{\rho}^{\ga +\beta }\int_0^x \bar\rho^{-2\beta-\ga}(y)dy\int_0^x\bar\rho^{\ga} y^2\lt(|r_y-1|^2+\lt|\frac{r}{y}-1\rt|^2\rt)dydx\\
&\le C\int_{1-\oa}^1\bar{\rho}^{-\beta }dx \int_0^1x^2\bar\rho^{\ga} \lt(|r_x-1|^2+\lt|\frac{r}{x}-1\rt|^2\rt)dx\\
&\le C \int_0^1x^2\bar\rho^{\ga} \lt(|r_x-1|^2+\lt|\frac{r}{x}-1\rt|^2\rt)dx.
\end{align*}
 Thus, it follows from these  and \ef{lem2est} that
\begin{equation}\label{14l2}\begin{split}
\int_0^t  |L_2|ds
\le  C(\omega) Q(0) +C \oa \int_0^t   \int x^2\bar{\rho}^{\ga-\beta}\lt(|r_x-1|^2+\lt|\frac{r}{x}-1\rt|^2\rt)dxds.
\end{split}
\end{equation}
Rewrite  $K_3$ as
\begin{equation*}\label{}\begin{split}
K_3= &- 4 \la_1 \int_0^{1-\oa} \left(\frac{v}{r}\right)_x  \int_0^x \lt(\bar\rho^{-\beta}\rt)_y(r^2 v)dy  dx \\
&- 4 \la_1 \int_{1-\oa}^1 \left(\frac{v}{r}\right)_x \lt[ \bar\rho^{-\beta} r^2 v - \lt(\int_0^x \bar\rho^{-\beta}(r^2 v)_ydy \rt) \rt]dx=:K_{31}+K_{32}.
\end{split}
\end{equation*}
$K_{31}$ can be bounded by
\begin{equation*}\label{}\begin{split}
|K_{31}|= &\lt| 4\frac{\beta}{\ga} \la_1 \int_0^{1-\oa} \left(\frac{v}{r}\right)_x  \int_0^x \bar\rho^{-\beta-(\ga-1)}y \phi r^2 vdy  dx \rt|\\
\le & C   \int_0^{1-\oa} ( |xv_x| +|v| ) x^{-2}  \lt(\int_0^x y^3  |v| dy\rt)  dx\\
\le&  C  \int_0^{1-\oa} ( |xv_x| +|v| )   \lt(\int_0^1 v^2 dy\rt)^{1/2}  dx
\le   C   \int (v^2+x^2 v_x^2) dx,
\end{split}
\end{equation*}
and $K_{32}$ can be bounded by
\begin{equation*}\label{}\begin{split}
|K_{32}|\le &C  \int_{1-\oa}^1 \bar\rho^{-\beta} ( |xv_x| +|v| ) |v| dx + \int_{1-\oa}^1 ( |xv_x| +|v| )  \lt(\int_0^x \bar\rho^{-\beta}\lt(|v|+|yv_y|\rt)dy \rt) dx\\
\le & \oa \int_{1-\oa}^1 \bar\rho^{-\beta} ( |xv_x|^2 +|v|^2 ) dx + \oa^{-1} \int_{1-\oa}^1 \bar\rho^{-\beta} |v|^2 dx \\
&+ \oa^{-1} \int_{1-\oa}^1 ( |xv_x|^2 +|v|^2 ) dx + \oa  \int_0^1 \bar\rho^{-\beta} (| v|^2 + |y v_y|^2 )dy,
\end{split}
\end{equation*}
due to \ef{fact}.
Since $\beta<\ga-1$,  the Hardy inequality \ef{hardybdry} implies that
$$
 \int_{1-\oa}^1 \bar\rho^{-\beta}  |v|^2 dx
\le C \int_{1/2}^1 \bar\rho^{-\beta+2(\ga-1)} (v^2+v_x^2) dx \le C\int (v^2+v_x^2) dx.
$$
These, together with \ef{lem2est'}, yield
\begin{equation}\label{14k3}\begin{split}
&\int_0^t (1+s)^{\nu}|K_3|ds \\
\le & C\oa^{-1}\int_0^t (1+s) \int (x^2v_x^2 + v^2)dxds + C\oa \int_0^t (1+s)^{ \nu}  \int_{1-\oa}^1 \bar\rho^{-\beta} ( |xv_x|^2 +|v|^2 ) dxds\\
\le & C(\omega) Q(0) + C\oa \int_0^t (1+s)^{\nu}  \int_{1-\oa}^1 \bar\rho^{-\beta} ( |xv_x|^2 +|v|^2 ) dxds.
\end{split}
\end{equation}
Similarly, $L_3$ is rewritten as
\begin{equation*}\label{}\begin{split}
L_3= &- 4 \la_1 \int_0^{1-\oa} \left(\frac{v}{r}\right)_x  \int_0^x \lt(\bar\rho^{-\beta}\rt)_y(r^3-y^3)dy  dx \\
&- 4 \la_1 \int_{1-\oa}^1 \left(\frac{v}{r}\right)_x \lt[ \bar\rho^{-\beta}(r^3-x^3) - \lt(\int_0^x \bar\rho^{-\beta}(r^3-y^3)_ydy \rt) \rt]dx=:L_{31}+L_{32}.
\end{split}
\end{equation*}
Clearly, $L_{31}$ and $L_{32}$ can be bounded by
$$
|L_{31}|\le  C   \int (v^2+x^2 v_x^2) dx +  C   \int \bar\rho^{\ga}(r-x)^2 dx
$$
and
$$
|L_{32}| \le   C\int_{1-\oa}^1 \bar{\rho}^{-\beta}\lt(|xv_x|+\lt|v\rt|\rt)|r-x|dx
  +C \int_{1-\oa}^1 \lt(|xv_x|+\lt|v\rt|\rt)\lt(\int_0^x \bar\rho^{-\beta}y|r_y-1|dy \rt) dx.
$$
The second term in $L_{32}$ is bounded by
\begin{equation*}\label{}\begin{split}
 & \int_{1-\oa}^1 \lt(|xv_x|+\lt|v\rt|\rt)\lt(\int_0^x \bar\rho^{-\beta}y|r_y-1|dy \rt) dx\\
\le& C\oa^{\frac{1}{2}}(1+t)^{2\nu}\int_{1-\oa}^1 ( |xv_x|^2 +|v|^2 ) dx\\
& +C\oa^{-\frac{1}{2}}
(1+t)^{-1-\frac{\beta}{\ga}+\iota}\int_{0}^1 \bar\rho^{-\beta} y^2 (r_y-1)^2 dy
\int_{1-\oa}^1 \int_0^x \bar\rho^{-\beta}dydx\\
\le & C\oa^{\frac{1}{2}}\lt[(1+t)^{2\nu}\int_{1-\oa}^1 ( |xv_x|^2 +|v|^2 ) dx
+(1+t)^{-1-\frac{\beta}{\ga}+\iota}\int_{0}^1 \bar\rho^{-\beta} y^2 (r_y-1)^2 dy\rt],
\end{split}
\end{equation*}
due to \ef{fact}. The first term in $L_{32}$ can be bounded as follows. Since $\beta< \ga-1$, it follows from \ef{phy} and the Hardy inequality \ef{hardybdry} that for $\ga>4/3$,
  \begin{equation*}\label{}\begin{split}
&\int_{1-\oa}^1 \bar{\rho}^{-\beta}\lt(|xv_x|+\lt|v\rt|\rt)|r-x|dx \\
\le & C \oa^{\frac{h}{2(\ga-1)}} \lt[ (1+t)^{\nu} \int_{1-\oa}^1 \bar{\rho}^{-\beta}\lt(v_x^2+v^2\rt)dx
+    (1+t)^{-\nu}\int_{1-\oa}^1 \bar\rho^{-\beta-h}|r-x|^2dx\rt]\\
\le & C \oa^{\frac{h}{2(\ga-1)}} \lt[ (1+t)^{\nu} \int_{1-\oa}^1 \bar{\rho}^{-\beta}\lt(v_x^2+v^2\rt)dx
 +\int_{1/2}^1 \bar{\rho}^{\ga-\beta}\lt[(r-x)^2+x^2 (r_x-1)^2 \rt]dx \rt.\\
 &\lt. + (1+t)^{-\frac{\nu \ga}{h+2-\ga}} \int_{1/2}^1 \bar{\rho}^{-\beta}\lt[(r-x)^2+x^2 (r_x-1)^2 \rt]dx \rt],
\end{split}
\end{equation*}
where
$h=\min\lt\{\beta/8, \  \  (\ga-1-\beta)/4\rt\}$ (it should be noted that $\beta+h<\ga-1$), and we have used the estimate
\begin{equation*}\label{}\begin{split}
&\int_{1-\oa}^1 \bar\rho^{-\beta-h}|r-x|^2dx \le
\int_{1/2}^1 \bar\rho^{2(\ga-1)-\beta-h}\lt[|r-x|^2+(r_x-1)^2 \rt]dx\\
\le & \lt(\int_{1/2}^1 \bar{\rho}^{-\beta}\lt[(r-x)^2+x^2 (r_x-1)^2 \rt]dx\rt)^{\frac{h+2-\ga}{\ga}} \lt(\int_{1/2}^1 \bar{\rho}^{\ga-\beta}\lt[(r-x)^2+x^2 (r_x-1)^2 \rt]dx\rt)^{\frac{2(\ga-1)-h}{\ga}}.
\end{split}
\end{equation*}
Consequently, taking into account of  \ef{lem1est}, \ef{lem2est} and \ef{fact}, one gets that
\begin{equation}\label{14l3}\begin{split}
  \int_0^t |L_3| ds \le & C(\oa) Q(0)+  C \lt[\oa^{\frac{1}{2}} + \oa^{\frac{h}{2(\ga-1)}} \rt]  \lt[\int_0^t \int \lt[(1+s)^{2\nu}+(1+s)^\nu \bar\rho^{-\beta}\rt] \rt.\\
 &\lt. \times  ( |xv_x|^2 +|v|^2 ) dxds + \sup_{[0,t]}\int  \bar{\rho}^{-\beta}\lt[(r-x)^2+x^2 (r_x-1)^2 \rt]dx  \rt.\\
 &\lt.+ \int_0^t \int  \bar{\rho}^{\ga-\beta}\lt[(r-x)^2+x^2 (r_x-1)^2 \rt]dx\rt].
\end{split}
\end{equation}
We finally derive from \ef{14eye} and \ef{14k1}-\ef{14l3}, by choosing $\oa\in (0, {1}/{2})$ suitably small, that
\bee\label{}\begin{split}
&\int \lt[(1+t)^{2\nu}\bar\rho^\ga
+(1+t)^{\kappa} +\bar{\rho}^{-\beta}\rt] x^2\left[\left(\frac{r}{x}-1\right)^2  + \left(r_x-1\right)^2 \right] (x,t) dx\\
&+(1+t)^{2\nu}\int x^2 \bar\rho v^2 (x,t) dx+\int_0^t   \int \lt[(1+s)^{\kappa} + \bar\rho^{-\beta} \rt]x^2\bar{\rho}^{\ga} \left[\left(\frac{r}{x}-1\right)^2  + \left(r_x-1\right)^2 \right] dxds\\
&+ \int_0^t \int \lt[(1+s)^{2\nu}+ (1+s)^{\nu}  \bar\rho^{-\beta} \rt]\left(x^2 v_x^2+ v^2  \right)dxds
\le  C Q(0) .
\end{split}\eee
Due to \ef{alphaiota}, this completes the proof of \ef{step4}.

{\em Step 4}.
Multiply equation \ef{later1} by $(1+t)^{\frac{2\ga-1}{\ga}-\ta}$ and integrate the product to deduce that
$$
(1+t)^{\frac{2\ga-1}{\ga}-\ta}\int \lt[\bar{\rho} v_t^2 +  \bar{\rho}^\ga \left(v^2 + x^2 v_x^2 \right) \right] dx   +  \int_0^t (1+s)^{\frac{2\ga-1}{\ga}-\ta} \int \left( x^2 v_{sx}^2+  v_s^2  \right)dxds
\le C Q(0) .
$$
In a similar way to the derivation of \ef{rminusx} and \ef{jump}, one can show
$$
   \int(r-x)^2(x,t)dx
 \le  C(1+t)^{-\frac{3(\ga-1)}{\ga}+\theta} Q(0)
$$
and
$$
   \int \lt(v^2 + x^2 v_x^2\rt)(x,t)dx
 \le C (1+t)^{-\frac{2\ga-1}{\ga}+\ta} Q(0).
  $$
This finishes the proof of \ef{estlem51}.

Moreover, it follows from \ef{estlem51} that for $x\in [0,1]$,
\begin{align}
xv^2(x,t)= & \int_0^x (y v^2(y,t))_y dy \le   2 \lt(\int v^2(y,t)dy\rt)^{1/2}\lt(\int y^2 v_y^2(y,t)dy\rt)^{1/2}  \notag\\
& + \int v^2(y,t)dy \le    C (1+t)^{-\frac{2\ga-1}{\ga}+\ta} Q(0). \label{8/12-2}
\end{align}
Similarly,
$$
x\lt(r(x,t)-x\rt)^2
\le   C (1+t)^{-\frac{2(\ga-1)}{\ga}+\ta} Q(0).
$$
This finish the proof of \ef{8/12-1}.

\hfill $\Box$

\subsection{Higher-order estimates}\label{sec3.4}
In this subsection, we derive the higher-order part of the {\it a priori } estimates for the strong solution $(r,v)$ on the time interval $[0, T]$ defined in Definition \ref{definitionss},
under the assumption \ef{rx} and \ef{vx}.   To obtain the higher-order estimates, we define
\be\label{mathG}
\mathcal{G} : =  \ln r_x + 2 \ln \lt(\frac{r}{x}\rt).
\ee
This transformation between $\mathcal{G}$ and $r$ is one-to-one, and we can solve for $r$ in terms of $\mathcal{G}$ by
\be\label{rg}
r(x, t)=\lt(3\int_0^x y^2 \exp\lt\{\mathcal G(y, t)\rt\}dy\rt)^{1/3} \  \ {\rm for} \ \  x\in \bar I \ \ {\rm and} \ \  \ t\ge 0.
\ee
Indeed, we will show in Section \ref{sec3.4.1} that $\mathcal{G}\sim r_x-1$, $\mathcal{G}_t \sim v_x$, $\mathcal{G}_x\sim r_{xx}$ and $\mathcal{G}_{tx}\sim v_{xx}$.
Then equation \ef{419a}  can be written in the  form of
\begin{equation}\label{7-2}\begin{split}
 \mu \mathcal{G}_{xt} +\ga \left(\frac{x^2 \bar{\rho}}{r^2 r_x } \right)^{\ga } \mathcal{G}_x
    =  \frac{x^2}{r^2} \bar{\rho}  v_t - \left[ \left(\frac{x^2 }{r^2 r_x } \right)^{\ga }  -\left(\frac{x }{r } \right)^{4}  \right]  x \phi \bar\rho.
\end{split}\end{equation}
(This is the same as \ef{viscosityequation}, we recall it here for the convenience of readers).

\subsubsection{Preliminaries for higher-order estimates}\label{sec3.4.1}
 The main goal of this subsubsection is to derive some preliminary estimates for the strong solution $(r,v)$ on the time interval $[0, T]$ defined in Definition \ref{definitionss},
 and to prove  the equivalence of the functionals $\mathfrak{E}(t)$ and $\mathcal{E}(t)$,
under the {\it a priori} assumptions \ef{rx} and \ef{vx}.

We illustrate how to use $\mathcal{G}$ and its derivatives to control
$r$ and $v$ and their derivatives by identifying the principal parts of $\mathcal G$, $\mathcal G_t$, $\mathcal G_x$    and $\mathcal G_{xt}$.
Note that
\begin{equation}\label{tlg1} \mathcal{G}=(r_x-1) + 2 \lt(\frac{r}{x}-1\rt)
+O\lt(|r_x-1|^2 + \lt|\frac{r}{x}-1\rt|^2\rt), \end{equation}
\begin{equation}\label{tlg2} \mathcal{G}_x=\frac{r_{xx}}{r_x}+ 2\frac{x}{r}\lt(\frac{r}{x}\rt)_x=\lt[r_{xx} + 2\lt(\frac{r}{x}\rt)_x\rt]+ \lt(\frac{1}{r_x}-1\rt)r_{xx}+2 \lt(\frac{r}{x}-1\rt)\lt(\frac{r}{x}\rt)_x, \end{equation}
\begin{equation}\label{tlg3} \mathcal{G}_t=\frac{v_x}{r_x}+2\frac{v}{r}=\lt(v_{x} + 2\frac{v}{x}\rt)+ \lt(\frac{1}{r_x}-1\rt)v_{x}+2 \lt(\frac{r}{x}-1\rt)\lt(\frac{v}{x}\rt), \end{equation}
\begin{align}\label{tlg5} \mathcal{G}_{xt}=\lt(\frac{v_x}{r_x}+2\frac{v}{r}\rt)_x =&  \left(v_x +2 \frac{v}{x}\right)_x + \left[ \left(\frac{1 }{r_x}-1\right) v_{xx} + 2\left(\frac{r}{x}-1\right) \left(\frac{v}{x}   \right)_x \right]
  \notag
  \\& -\left[\frac{r_{xx}}{r_x^2}v_x+2\left(\frac{x}{r}\right)^2\left(\frac{r}{x}\right)_x \frac{v}{x} \right]. \end{align}
Thus, it follows from \ef{hhreg1}-\ef{hhreg3} that for $t\in [0, T]$,
\begin{align}
& \lt(r_x, \ r/x, \ v_x, \ v/x, \  \mathcal{G}, \ \mathcal{G}_t    \rt) \in L^\iy(I), \ \    \lt( r/x,  \ v/x , \  v/r \rt) \in H^1(I), \ \ \mathcal{G}_{xt} \in L^2(I),   \label{9-21-1}\\
&  (r, \ v) \in H^2([0 ,a]) \ \ {\rm for} \ \ a \in (0, 1),\label{9-21-2} \\
&  (r, \ v) \in H^2(I), \ \  {\rm if} \ \  \mathcal{G}_x \in L^2.  \label{9-21-3}
\end{align}
In fact, \ef{9-21-3} can be derived easily by noting that
$$r_{xx}=r_x\lt[\mathcal{G}_x- 2 \frac{x}{r}\lt(\frac{r}{x}\rt)_x \rt] \ \ {\rm and} \ \
v_{xx}=r_x\lt[\mathcal{G}_{xt}+ \frac{v_x}{r_x^2}r_{xx} - 2  \lt(\frac{v}{r}\rt)_x \rt]. $$

With the regularity \ef{9-21-1}-\ef{9-21-3}, we  have the following Lemmas.

\begin{lem}\label{lem2.3}
Suppose that \ef{rx} holds for a suitably small $\ea_0$.
Then  for $t\in [0, T]$,
\begin{equation}\label{gjvx}\begin{split}
&\lt\|\lt(v_x, v/x\rt) \rt\|^2  \le 4 \lt\|\mathcal{G}_{t}\rt\|^2  ,
\end{split}\end{equation}
\begin{equation}\label{gjrx}\begin{split}
&\lt\|\lt(r_x-1, {r}/{x}-1\rt) \rt\|^2  \le 4 \lt\|\mathcal{G} \rt\|^2  ,
\end{split}\end{equation}
\begin{equation}\label{gjrxx}\begin{split}
&\lt\|\lt( r_{xx}, \ (r/x)_x \rt)\rt\|^2  \le 4 \lt\|\mathcal{G}_{x}\rt\|^2 ,
\end{split}\end{equation}
\begin{equation}\label{gjvxx}\begin{split}
&\lt\|\lt( v_{xx}, \ (v/x)_x \rt)\rt\|^2 \le  4 \lt\|\mathcal{G}_{tx}\rt\|^2 + c \lt\|(v_x, v/x)\rt\|_{L^\iy}^2   \lt\|\mathcal{G}_{x}\rt\|^2  .
\end{split}\end{equation}
Here the estimates \ef{gjrxx} and \ef{gjvxx} hold if $ \|\mathcal{G}_{x} \|< \iy$.
\end{lem}
{\em Proof}. The proof consists of three steps.

{\em Step 1}. In this step, we prove \ef{gjvx} and \ef{gjrx}.  Let $\varepsilon\in (0, 1/4)$ be an arbitrary constant, and  $\chi_\varepsilon\in[0,1]$ be a cut-off function satisfying
\begin{equation}\label{cutoff1}\begin{split}
&\chi_\varepsilon=0 \ \ {\rm on} \ \ [0,\varepsilon]\cup [1-\varepsilon, 1], \ \     \chi_\varepsilon=1 \  \ {\rm on} \ \  [2\varepsilon,1-2\varepsilon] ,\\
& \chi_\varepsilon=x/\varepsilon-1 \ \ {\rm on} \ \ [\varepsilon,2\varepsilon], \ \    \ \   \ \  \chi_\varepsilon=(1-x)/\varepsilon-1 \ \ {\rm on} \ \ [1-2\varepsilon,1-  \varepsilon] .
\end{split}\end{equation}
Note that $v_x=x(v/x)_x+v/x$ and $ v/x \in H^1$,  we use the integration by parts to get
$$
 \int \chi_\varepsilon  v_x  \frac{v}{x}dx=  \int \chi_\varepsilon \lt(\frac{v}{x}\rt)^2dx
+\frac{1}{2}\int  \chi_\varepsilon x \lt(\lt(\frac{v}{x}\rt)^2\rt)_xdx
 =  \frac{1}{2}\int \chi_\varepsilon \lt(\frac{v}{x}\rt)^2dx-\frac{1}{2}\int x\chi'_\varepsilon \lt(\frac{v}{x}\rt)^2dx.
 $$
Thus, we have
\begin{equation*}\label{}\begin{split}
 & \int \chi_\varepsilon\lt|v_x+2 \frac{v}{x}\rt|^2dx
 =\int \chi_\varepsilon\lt[ v_x^2 +4\lt(\frac{v}{x}\rt)^2 \rt]dx +4\int \chi_\varepsilon  v_x \frac{v}{x} dx \\
= & \int \chi_\varepsilon\lt[ v_x^2 +6\lt(\frac{v}{x}\rt)^2 \rt]dx -2\int \chi_\varepsilon' x \lt(\frac{v}{x}\rt)^2dx \\
=&\int \chi_\varepsilon\lt[ v_x^2 +6\lt(\frac{v}{x}\rt)^2 \rt]dx -2\int_\varepsilon^{2\varepsilon}  \frac{ x}{\varepsilon} \lt(\frac{v}{x}\rt)^2dx +2\int_{1-2\varepsilon}^{1-\varepsilon}  \frac{ x}{\varepsilon} \lt(\frac{v}{x}\rt)^2dx \\
\ge & \int \chi_\varepsilon\lt[ v_x^2 +6\lt(\frac{v}{x}\rt)^2 \rt]dx -4\int_\varepsilon^{2\varepsilon}  \lt(\frac{v}{x}\rt)^2dx .
 \end{split}\end{equation*}
Note that
\begin{equation*}\label{}\begin{split}
v_x+2\frac{v}{x} =\mathcal{G}_{t}-\lt(\frac{1}{r_x}-1\rt)v_x-2\lt(\frac{x    }{r}-1\rt)\frac{v}{x}
\end{split}\end{equation*}
Thus,
\begin{equation*}\label{}\begin{split}
& \int \chi_\varepsilon\lt[ v_x^2 +6\lt(\frac{v}{x}\rt)^2 \rt]dx
\le
 \int \chi_\varepsilon\lt|v_x+2 \frac{v}{x}\rt|^2dx + 4\int_\varepsilon^{2\varepsilon}  \lt(\frac{v}{x}\rt)^2dx \\
 \le & 2 \int \chi_\varepsilon \mathcal{G}_{t}^2 dx+ c\ea_0^2\int \chi_\varepsilon\lt[ v_x^2 +\lt(\frac{v}{x}\rt)^2 \rt]dx
 +4\int  \lt(\frac{v}{x}\rt)^2dx  ;
 \end{split}\end{equation*}
which implies that
\begin{equation}\label{ht23}\begin{split}
\int \chi_\varepsilon\lt[ \frac{1}{2} v_x^2 + \frac{11}{2}\lt(\frac{v}{x}\rt)^2 \rt]dx \le
2 \int  \mathcal{G}_{t}^2 dx
 +4\int  \lt(\frac{v}{x}\rt)^2dx,
\end{split}\end{equation}
provided that \ef{rx} holds for a suitably small number   $\ea_0$. Since $(  v_x, v/x) \in L^2$ due to \ef{9-21-1},    we can obtain by letting $\varepsilon\to 0$ and using the dominated convergence theorem that
\begin{equation*}\label{}\begin{split}
\int  \lt[ \frac{1}{2} v_x^2 + \frac{11}{2}\lt(\frac{v}{x}\rt)^2 \rt]dx \le
2 \int  \mathcal{G}_{t}^2 dx
 +4\int  \lt(\frac{v}{x}\rt)^2dx,
\end{split}\end{equation*}
which means
\begin{equation*}\label{sad}\begin{split}
\int  \lt[  v_x^2 +  \lt(\frac{v}{x}\rt)^2 \rt]dx \le
4 \int  \mathcal{G}_{t}^2 dx.
\end{split}\end{equation*}
This finishes the proof of \ef{gjvx}. Clearly, \ef{gjrx} follows  from similar arguments.

{\em Step 2}. In this step, we prove \ef{gjrxx} under the assumption $\|\mathcal{G}_x\|<\iy$. Due to $\mathcal{G}_x\in L^2(I)$ and \ef{9-21-3}, we have
$r\in H^2(I)$. Then, a similar argument to the proof of \ef{gjvx}  will lead to the proof of \ef{gjrxx},  which goes as follows.  For the cut-off function $\chi_\varepsilon$ defined in \ef{cutoff1}, we use  $r\in H^2(I)$  which ensures $(r/x)_x \in H^1([\varepsilon,1])$, and the integration by parts to get
\bee\label{9152}\begin{split}
\int \chi_\varepsilon r_{xx} \lt(\frac{r}{x}\rt)_xdx = & \int \chi_\varepsilon \lt(x\frac{r}{x}\rt)_{xx}\lt (\frac{r}{x}\rt)_xdx
=\int \chi_\varepsilon \lt\{\frac{1}{2} x\lt[\lt(\lt(\frac{r}{x}\rt)_x\rt)^2\rt]_x+2\lt|\lt(\frac{r}{x}\rt)_x\rt|^2\rt\}dx\\
 =  & \frac{3}{2}\int \chi_\varepsilon \lt|\lt(\frac{r}{x}\rt)_x\rt|^2dx -\frac{1}{2} \int x\chi_\varepsilon'  \lt|\lt(\frac{r}{x}\rt)_x\rt|^2dx,
 \end{split}\eee
which implies
\bee\label{9151}\begin{split}
\int \chi_\varepsilon \lt|r_{xx}+2\lt(\frac{r}{x}\rt)_x\rt|^2dx= & \int \chi_\varepsilon\lt( r_{xx}^2+10\lt|\lt(\frac{r}{x}\rt)_x\rt|^2\rt)dx- 2\int x\chi_\varepsilon'  \lt|\lt(\frac{r}{x}\rt)_x\rt|^2dx\\
\ge & \int \chi_\varepsilon\lt( r_{xx}^2+10\lt|\lt(\frac{r}{x}\rt)_x\rt|^2\rt)dx - 4 \int_\varepsilon^{2\varepsilon} \lt|\lt(\frac{r}{x}\rt)_x\rt|^2dx .
 \end{split}\eee
This, together with \ef{tlg2} and \ef{rx}, gives that for small $\epsilon_0$ in \ef{rx},
\be\label{9151} \int \chi_\varepsilon\lt(\frac{1}{2}r_{xx}^2+\frac{19}{2}\lt|\lt(\frac{r}{x}\rt)_x\rt|^2\rt) dx\le 2 \int   \mathcal{G}_x^2dx+ 4 \int  \lt|\lt(\frac{r}{x}\rt)_x\rt|^2dx.\ee
Therefore, by virtue of \ef{9-21-1} and \ef{9-21-3} (which implies $(r_{xx}, (r/x)_x)\in L^2$), we obtain \ef{gjrxx} with the help of the dominated convergence theorem.

{\em Step 3}. In this step, we prove \ef{gjvxx} under the assumption $\|\mathcal{G}_x\|<\iy$.   In view of \ef{tlg5}, it follows from the Cauchy inequality and \eqref{rx} that
\begin{equation*}\label{nj4}\begin{split}
  \mathcal{G}_{xt}^2 \ge   \frac{1}{2}\left|\left(v_x +2 \frac{v}{x}\right)_x \right|^2 -c \ea_0^2 \left(v_{xx}^2 +  \left|\left(\frac{v}{x}\right)_x\right|^2\right)
   - c\lt\|\lt(v_x,\frac{v}{x}\rt)\rt\|_{L^\iy}^2\left(r_{xx}^2 +  \left|\left(\frac{r}{x}\right)_x\right|^2\right).
    \end{split}\end{equation*}
Due to $\mathcal{G}_x\in L^2(I)$ and \ef{9-21-3}, we have
$v\in H^2(I)$.  So,   we can use the same way as that for the derivation of \ef{9151} to obtain
\begin{equation}\label{9.22.1}\begin{split}
& \int \chi_\varepsilon\lt( \frac{1}{2} v_{xx}^2 + \frac{19}{2}\lt|\lt(\frac{v}{x}\rt)_x\rt|^2 \rt)dx \\
 \le &
2 \int   \mathcal{G}_{tx}^2 dx
 +4\int  \lt|\lt(\frac{v}{x}\rt)_x\rt|^2  dx
 + c\lt\|\lt(v_x,\frac{v}{x}\rt)\rt\|_{L^\iy}^2\int  \left(r_{xx}^2 +  \left|\left(\frac{r}{x}\right)_x\right|^2\right)dx ,
\end{split}\end{equation}
where $\chi_\varepsilon$ is the cut-off function  defined in \ef{cutoff1}.  With the aid of \ef{9-21-1} and \ef{9-21-3}, we see that all the quantities appearing on the right-hand side of \ef{9.22.1} are finite.
Since $(v_{xx}, (v/x)_x) \in L^2$ due to \ef{9-21-1} and \ef{9-21-3}, we can get, by letting $\varepsilon\to 0$ and using the dominated convergence theorem, that
\begin{equation}\label{3.3.12Aug6}\begin{split}
 \int  \lt(  v_{xx}^2 +  \lt|\lt(\frac{v}{x}\rt)_x\rt|^2 \rt)dx
 \le
4 \int   \mathcal{G}_{xt}^2 dx
 + c\lt\|\lt(v_x,\frac{v}{x}\rt)\rt\|_{L^\iy}^2\int  \left(r_{xx}^2 +  \left|\left(\frac{r}{x}\right)_x\right|^2\right)dx .
\end{split}\end{equation}
This, together with \ef{gjrxx}, gives \ef{gjvxx}.

\hfill $\Box$

\begin{rmk}   In view of \ef{9-21-1}, we know that
$\mathcal{G}_{xt}\in L^2(I)$,  $ t\in [0, T] $.
 However, this does not mean that $\mathcal{G}_{x}(t)\in L^2(I)$,  $ t\in [0, T] $, unless we assume that $\mathcal{G}_{x}(0)\in L^2(I)$, because
 $$\mathcal{G}_{x}(t)=\mathcal{G}_{x}(0)+\int_0^t \mathcal{G}_{xs}(x, s)ds ,  \ t\in [0, T]. $$
$\mathcal{G}_{x}(0)\in L^2(I)$ is an additional regularity assumption of the initial data other than $\mathfrak{E}(0)<\infty$. By \ef{gjrxx}, \ef{tlg2} and \ef{9.21.6}, we know that the condition $\mathcal{G}_{x}(0)\in L^2(I)$ is equivalent to $ r_{xx}   (x, 0)\in L^2(I)$ if $\mathfrak{E}(0) $ is small.
\end{rmk}

\begin{rmk} We can see  from the proof of  \ef{gjvxx}, in particular, \ef{3.3.12Aug6},  that
\begin{align}\label{Aug6-1}
\lt\|\lt( v_{xx}, \ (v/x)_x \rt)\rt\|^2 \le & 4 \lt\|\mathcal{G}_{tx}\rt\|^2 + c \lt\|(v_x, v/x)\rt\|_{L^\iy}^2 \lt\|\lt( r_{xx}, \ (r/x)_x \rt)\rt\|^2 ,
\end{align}
if $\|r_{xx}\|<\iy$.  To bound $\|( v_{xx},  (v/x )_x )\|$, we have to control the product of  $ \| ( r_{xx},  (r/x  )_x )\|$ and $ \|(v_x, v/x) \|_{L^\iy}$, in addition to the bound of $ \|\mathcal{G}_{xt} \|$.   Indeed, we prove in Lemma \ref{lem10} that although $ \| ( r_{xx}, \ (r/x  )_x)(\cdot, t) \|$ may grow with respect to time,
$ \|(v_x, v/x)(\cdot, t) \|_{L^\iy}$ decays faster  so that the product is bounded.
\end{rmk}

\begin{rmk} Similar to \ef{gjvxx}, we can obtain that for any $a\in (0,1)$,
\begin{equation*}\label{}\begin{split}
\lt\|\lt( v_{xx}, \ (v/x)_x \rt)\rt\|_{L^2([a,1])}^2 \le  4 \lt\|\mathcal{G}_{tx}\rt\|^2_{L^2([a,1])} + c \lt\|(v_x, v/x)\rt\|_{L^\iy ([a,1])  }^2 \lt\|\mathcal{G}_{x}\rt\|^2_{L^2([a,1])},
\end{split}\end{equation*}
which implies
\begin{equation}\label{gjvxxa}\begin{split}
a^2 \lt\|\lt( v_{xx}, \ (v/x)_x \rt)\rt\|_{L^2([a,1])}^2 \le  4 \lt\|x\mathcal{G}_{tx}\rt\|^2  + c \lt\|(xv_x, v)\rt\|_{L^\iy }^2 \lt\|\mathcal{G}_{x}\rt\|^2 ,
\end{split}\end{equation}
provided that $\lt\|\mathcal{G}_{x}\rt\|<\iy$.
\end{rmk}

\begin{lem}\label{lemhh1} Let $\da>0$ be any fixed constant. Suppose that \ef{rx} holds for a suitable small $\ea_0$, then for $t\in [0 ,T]$,
\begin{align}
&\lt\|\bar\rho^{\da} \lt( r_{xx}, (r/x)_x \rt) \rt\|^2 +\lt\|x \bar\rho^{\da-(\ga-1)/2} \lt( r/x \rt)_x\rt\|^2 \le c \lt\|\bar\rho^{\da} \mathcal{G}_{x}\rt\|^2, \label{weightrxx} \\
 &\lt\|\bar\rho^{\da} \lt( v_{xx}, \lt({v}/{x}\rt)_x \rt)\rt\|^2   \le c \lt\|\bar\rho^{\da} \mathcal{G}_{xt} \rt\|^2 +c\lt\|\lt(v_x, {v}/{x}\rt)\rt\|_{L^\iy}^2  \lt\|\bar\rho^{\da} \mathcal{G}_{x}\rt\|^2, \label{weightvxx} \\
&\lt\|\bar\rho^{\da} \lt(r_{x}-1,  {r}/{x}-1 \rt)\rt\|^2 +\lt\| \bar\rho^{\da-(\ga-1)/2} (r-x) \rt\|^2 \le c \lt\|\bar\rho^{\da} \mathcal{G} \rt\|^2 , \label{weightrx} \\
&\lt\|\bar\rho^{\da} \lt( v_{x}, v/x \rt)\rt\|^2 +\lt\| \bar\rho^{\da-(\ga-1)/2} v \rt\|^2 \le c \lt\|\bar\rho^{\da} \mathcal{G}_{t} \rt\|^2  . \label{weightvx}
\end{align}
Here the estimates \ef{weightrxx}  and \ef{weightvxx}  hold if $\lt\|\bar\rho^{\da} \mathcal{G}_{x}\rt\|<\iy$.
\end{lem}
{\em Proof}. The idea of the proof of this lemma is similar to that of Lemma \ref{lem2.3}. Due to $\bar\rho(1)=0$, we only need to cut off the origin. Let $\varepsilon\in (0, 1/4)$ be an arbitrary constant, and  $\eta_\varepsilon\in[0,1]$ be a cut-off function satisfying
\begin{equation*}\label{}\begin{split}
&\eta_\varepsilon=0 \ \ {\rm on} \ \ [0,\varepsilon],\  \
& \eta_\varepsilon=x/\varepsilon-1 \ \ {\rm on} \ \ [\varepsilon,2\varepsilon]  \ \  {\rm and} \ \    \chi_\varepsilon=1 \  \ {\rm on} \ \  [2\varepsilon,1].
\end{split}\end{equation*}
In a similar way to deriving \ef{ht23}, we can get
\begin{equation*}\label{ht23Aug8}\begin{split}
 \int \eta_\varepsilon \bar\rho^{2\da} \lt[ \frac{1}{2} v_x^2 + \frac{11}{2}\lt(\frac{v}{x}\rt)^2 \rt]dx + 4\frac{\da}{\ga}\int \eta_\varepsilon \phi \bar\rho^{2\da-(\ga-1)} v^2 dx
\le
2 \int  \bar\rho^{2\da}  \mathcal{G}_{t}^2 dx
 +4\int \bar\rho^{2\da}  \lt(\frac{v}{x}\rt)^2dx,
\end{split}\end{equation*}
due to \ef{rhox}. Letting $\varepsilon\to 0$ to give \ef{weightvx}.
Clearly, \ef{weightrxx}-\ef{weightrx} follow from similar arguments.

\hfill$\Box$


\begin{lem}
Suppose that \ef{rx} holds for a suitably small $\ea_0$.
Then  for any $a\in (0,1)$,
\begin{equation}\label{f2}\begin{split}
\lt\|v\rt\|_{L^\iy}^2  \le 2\lt\|v\rt\|\|v_x\|,  \ \   t\in [0, T],
\end{split}\end{equation}
\begin{equation}\label{f1}\begin{split}
\lt\|x v_x\rt\|_{L^\iy}^2 \le c \lt\|xv_x\rt\|\lt( \lt\|x\mathcal{G}_{tx} \rt\|+\lt\|v_x\rt\|+ \lt\|v/x\rt\|\rt),  \ \   t\in [0, T],
\end{split}\end{equation}
\begin{equation}\label{nj2}\begin{split}
 &\lt\|v_x \rt\|_{L^\iy([0, a ])}^2 \le (1/a)\lt\| v_x  \rt\|^2_{L^2\lt([0, a]\rt)} + 2\lt\|v_x\rt\|_{L^2\lt([0, a]\rt)}
\lt\| v_{xx}\rt\|_{L^2\lt([0, a]\rt)},  \ \   t\in [0, T],\\
\end{split}\end{equation}
\begin{equation}\label{nj2.new}\begin{split}
 &\lt\|v/x \rt\|_{L^\iy([0, a ])}^2 \le (1/a)\lt\| v/x  \rt\|^2_{L^2\lt([0, a]\rt)} +
 2\lt\|v/x\rt\|_{L^2\lt([0, a]\rt)}\lt\| (v/x)_x\rt\|_{L^2\lt([0, a]\rt)},  \ \   t\in [0, T].
\end{split}\end{equation}
\end{lem}
{\em Proof}.   Clearly, \ef{f2} follows from the  boundary condition $v(0,t)=0$ and the H$\ddot{\rm o}$lder inequality.
For $xv_x$, notice that
$$
\lt(x v_x\rt)^2 =r_x^2 \lt(x \frac{v_x}{r_x}\rt)^2
= 2 r_x^2 \int_0^x \lt(y \frac{v_y}{r_y}\rt) \lt(y \frac{v_y}{r_y}\rt)_y dy
 \le c \lt\|x \frac{v_x}{r_x}\rt\|  \lt(\lt\|x\lt( \frac{v_x}{r_x}\rt)_x\rt\| +\lt\|\frac{v_{x}}{r_x}\rt\| \rt)
$$
and
$$
\lt\|x\lt( \frac{v_x}{r_x}\rt)_x\rt\|=\lt\|x\lt( \mathcal{G}_t - 2\frac{v}{r} \rt)_x\rt\| \le \lt\|x\mathcal{G}_{tx} \rt\| + 2\lt\|x \lt(\frac{v}{r} \rt)_x\rt\|.
$$
Thus,
\begin{equation*}\label{}\begin{split}
\lt\|x v_x\rt\|_{L^\iy}^2 \le c \lt\|xv_x\rt\|\lt( \lt\|x\mathcal{G}_{tx} \rt\|+\lt\|v_x\rt\|+ \lt\|v/x\rt\|\rt),
\end{split}\end{equation*}
which verifies \ef{f1}.   \ef{nj2} and \ef{nj2.new} follow from simple calculations.
\hfill$\Box$

\begin{lem} Let $\da$ be a fixed positive constant. Then for $t\in [0, T]$,
\begin{equation}\label{girl}\begin{split}
&\lt\|r-x\rt\|_{L^\iy}^2  \le 2\lt\|r-x\rt\|\|r_x-1\|,
\end{split}\end{equation}
\begin{equation}\label{woman}\begin{split}
&\lt\|x^{{3}/{2}}\bar\rho^\da(r_x-1)\rt\|_{L^\iy}^2   \le
3\lt\|x\bar\rho^{\da}(r_x-1)\rt\|^2 + 2\lt\|x^3 \bar\rho^{2\da}(r_x-1)r_{xx}\rt\|_{L^1},
\end{split}\end{equation}
\begin{equation}\label{woman.new}\begin{split}
&\lt\|x^{{3}/{2}}\bar\rho^\da(r/x-1)\rt\|_{L^\iy}^2   \le
3\lt\|x\bar\rho^{\da}(r/x-1)\rt\|^2 + 2\lt\|x^3 \bar\rho^{2\da}(r/x-1)(r/x)_x\rt\|_{L^1}.
\end{split}\end{equation}
Here \ef{woman} holds if the quantities appearing on its right-hand side are finite.
\end{lem}
{\em Proof}. Clearly, \ef{girl} follows from \ef{Aug9-2} and the H$\ddot{\rm o}$lder inequality.
For any $x\in[0,1]$,
\begin{equation*}\label{}\begin{split}
&x^3 \bar\rho^{2\da}(r_x-1)^2= 2\int_0^x  \lt(y^3 \bar\rho^{2\da}(r_y-1)^2\rt)_ydy \\
\le &  3 \int_0^x y^2 \bar\rho^{2\da}(r_y-1)^2dy + 2\int_0^x y^3 \bar\rho^{2\da}(r_y-1)r_{yy}dy,
\end{split}\end{equation*}
due to \ef{rhox}. This gives \ef{woman}. Similarly, we can obtain \ef{woman.new}.

\hfill$\Box$

The following lemma is on the equivalence of the functionals $\mathcal{E}(t)$ and $\mathfrak{E}(t)$ and is the key to the verification of the {\it a priori} assumptions \ef{rx} and \ef{vx}.
\begin{lem}\label{boundsforrv} Suppose that \ef{rx} and \ef{vx} hold for  suitably small numbers  $\ea_0$ and $\ea_1$. Then,
\begin{align}\label{bfrv}
\lt\|\lt(r_x-1, {r}/{x}-1, v_x, {v}/{x}\rt)(\cdot,t)\rt\|_{L^\iy}^2 \le C\mathcal{E}(t), \ \  t\in [0, T],
\end{align}
\be\label{equivalence}
c \mathcal{E}(t)\le \mathfrak{E}(t)\le C \mathcal{E}(t), \ \  t\in [0, T].
\ee
\end{lem}
{\em Proof}.  The proof of \ef{bfrv} consists of two steps, in which the $L^\iy$-bounds on the intervals
$I_1=[0,1/2]$ and $I_2=[1/2,1]$ will be shown, respectively. Once  \ef{bfrv} is proved, \ef{equivalence} follows then from the definitions of $\mathcal{E}(t)$ and $\mathfrak{E}(t)$
by noticing that $v(0, t)=r(0,t)=0$.

{\em Step 1} (away from the boundary). It follows from \ef{7-2}, \ef{tlg2} and \ef{mathcalE} that
\be\label{Aug11.1}
\| \mathcal{G}_{xt}\|^2 \le C \lt(\|\bar\rho^\ga\mathcal{G}_x\|^2 + \|\bar\rho  v_t\|^2 + \|\bar\rho(r-x, xr_x-x)\|^2   \rt) \le C \mathcal{E},
\ee
which, together with \ef{weightvxx} and \ef{tlg2}, implies
\bee
\|\bar\rho^{\ga-1/2}(v_{xx}, (v/x)_x)(\cdot,t)\|^2 \le C \mathcal{E}(t).
\eee
Thus,
\bee
\|(r_{xx}, (r/x)_x,  v_{xx}, (v/x)_x)(\cdot,t)\|^2_{L^2(I_1)}\le C \mathcal{E}(t).
\eee
In view of \ef{hardyorigin}, we then have
\begin{equation}\label{Aug10-1}\begin{split}
\|(r_{x}-1, r/x-1,  v_{x}, v/x)(\cdot,t)\|^2_{L^2(I_1)}\le C \|(xr_{x}-x, r-x,  xv_{x}, v)(\cdot,t)\|^2_{L^2(I_1)} \\
+ C \|(xr_{xx}, x(r/x)_x,  xv_{xx}, x(v/x)_x)(\cdot,t)\|^2_{L^2(I_1)}
 \le C \mathcal{E}(t).
\end{split}\end{equation}
Hence,
\begin{equation}\label{Aug9-4}\begin{split}
&\|(r_{x}-1, r/x-1,  v_{x}, v/x)(\cdot,t)\|^2_{L^\iy(I_1)}
\\
  \le & C\|(r_{x}-1, r/x-1,  v_{x}, v/x)(\cdot,t)\|^2_{H^1(I_1)}
 \le C \mathcal{E}(t).
\end{split}\end{equation}

{\em Step 2} (away from the origin). It follows from   \ef{mathcalE}  and \ef{Aug10-1} that
\begin{equation}\label{Aug10-2}\begin{split}
\|(r_{x}-1, r/x-1,  v_{x}, v/x)(\cdot,t)\|^2
 \le C \mathcal{E}(t),
\end{split}\end{equation}
which implies, with the aid of \ef{mathcalE}, that
\begin{align}
& \|( {r}/{x}-1)(\cdot,t)\|_{L^\iy(I_2)}^2\le 4  \|( {r}-x)(\cdot,t)\|_{L^\iy(I_2)}^2 \le C \lt\|\lt(r-x \rt)(\cdot,t)\rt\|_{H^1(I_2)}^2 \le C \mathcal{E}(t), \label{Aug10-3}\\
& \|(v/x)(\cdot,t)\|_{L^\iy(I_2)}^2\le 4 \| v(\cdot,t)\|_{L^\iy(I_2)}^2 \le C \lt\|v(\cdot,t)\rt\|_{H^1(I_2)}^2 \le C \mathcal{E}(t). \label{Aug10-4}
 \end{align}
Clearly,
\be\label{Aug10-5}
\|(r_x-1)(\cdot,t)\|_{L^\iy(I_2)}^2\le \mathcal{E}(t).
\ee
In view of  \ef{f1}, \ef{Aug10-2} and \ef{Aug11.1}, we see that
\begin{align}\label{Aug10-6}
\|v_x(\cdot, t)\|_{L^\iy(I_2)}^2 \le 4 \lt\|(x v_x)(\cdot,t)\rt\|_{L^\iy(I_2)}^2 \le 4 \lt\|(x v_x)(\cdot,t)\rt\|_{L^\iy}^2\le C   \mathcal{E}(t).
\end{align}
So, \ef{bfrv} is a consequence of \ef{Aug9-4}, \ef{Aug10-3}-\ef{Aug10-6}.

\hfill$\Box$

\subsubsection{Part I:   global existence and decay of strong solutions}\label{sec3.4.2}
In this subsubsection, we prove the global existence and large time decay of the strong solution for suitably small $\mathfrak{E}(0)$.

\begin{lem}\label{lem5} Suppose that \ef{rx} and \ef{vx} hold. Then there exist positive constants  $C$, $C(\theta)$ and $C(\ta,a)$ independent of $t$ such that
  for  any $\theta\in (0,  \ {2(\ga-1)}/({3\ga}))$ and $a\in (0,1)$,
 \begin{align}
&\left\|\bar\rho^{\ga-{1}/{2}}\lt(r_{xx}, (r/x)_x\rt)(\cdot, t) \rt\|^2
 \le C \mathcal{E}(0),   \label{lem5est}\\
& (1+t)^{(\ga-1)/\ga-\ta}  \left\| \lt(r_{xx}, (r/x)_x\rt)(\cdot, t) \rt\|^2_{L^2([0,a])} \le  C(\ta, a) \lt(\mathcal{E}(0)+\|r_{0x}-1\|^2_{L^\iy}\rt), \label{lem5est8/9} \\
&(1+t)^{({\ga-1})/{\ga}-\ta}   \left\| \bar\rho^{ {1}/{2}} v_t(\cdot, t) \rt\|^2  +\int_0^t (1+s)^{({\ga-1})/{\ga}-\ta} \lt\|\lt( v_{x},  {v}/{x} , v_{sx},  {v_s}/{x}  \rt)(\cdot,s)\rt\|^2  ds\notag\\
&\quad   \le C(\theta)\lt(\mathcal{E}(0)+\|r_{0x}-1\|_{L^{\infty}}^2\rt),   \label{lem5est'}\\
& (1+t)^{(\ga-1)/\ga-\ta}  \left\| \lt(v_{xx}, (v/x)_x\rt)(\cdot, t) \rt\|^2_{L^2([0,a])} \le  C(\ta, a) \lt[\mathcal{E}(0)+\|r_{0x}-1\|^2_{L^\iy}\rt], \label{lem5est'8/9}
\end{align}
for all $t\in [0, T]$.
\end{lem}
{\em Proof}.  The proof consists of four steps.

{\em Step 1}. In this step, we prove \ef{lem5est} and \ef{lem5est8/9}. Multiplying equation \eqref{7-2}   by
$ \bar\rho^{2\ga-1}\mathcal{G}_x $ and integrating the product with respect to the spatial variable, one  gets,  using the Cauchy inequality,  that
\begin{equation}\label{catch}\begin{split}
 & \frac{\mu}{2}\frac{d}{dt} \int  \bar\rho^{2\ga-1}   \mathcal{G}_{x}^2dx +\frac{\ga}{2} \int  \left(\frac{x^2}{r^2 r_x } \right)^{\ga } \bar\rho ^{3\ga-1} \mathcal{G}_x^2 dx \\
    \le  & C\int   v_t^2dx  +C  \int x^2\bar{\rho}^\ga\left[ \left(\frac{r}{x}-1\right)^2  + \left(r_x-1\right)^2 \right] dx .
\end{split}\end{equation}
It follows from \ef{catch},  \eqref{lem2est} and \eqref{lem4est}   that
\begin{equation}\label{bt1}\begin{split}
 &  \int  \bar\rho^{2\ga-1}   \mathcal{G}_{x}^2 (x,t)dx + \int_0^t  \int  \bar\rho ^{3\ga-1} \mathcal{G}_x^2 dx ds \le C\mathcal{E}(0).
\end{split}\end{equation}
This, together with \ef{weightrxx}, gives \ef{lem5est}.

In a similar way to deriving \ef{catch}, we can get
\begin{equation*}\label{catch8/9}\begin{split}
 & \frac{\mu}{2}\frac{d}{dt} \int  \bar\rho^{3\ga-2}   \mathcal{G}_{x}^2dx +\frac{\ga}{2} \int  \left(\frac{x^2}{r^2 r_x } \right)^{\ga } \bar\rho ^{4\ga-2} \mathcal{G}_x^2 dx \\
    \le  & C\int   v_t^2dx  +C  \int x^2\bar{\rho}^\ga\left[ \left(\frac{r}{x}-1\right)^2  + \left(r_x-1\right)^2 \right] dx .
\end{split}\end{equation*}
Multiply the inequality above by $(1+t)^{({\ga-1})/{\ga}-\ta}$ and integrate the product to give
\begin{align}
 & (1+t)^{\frac{\ga-1}{\ga}-\ta} \int \bar\rho^{3\ga-2}    \mathcal{G}_{x}^2(x,t)dx +\int_0^t (1+s)^{\frac{\ga-1}{\ga}-\ta} \int \bar\rho^{4\ga-2} \mathcal{G}_x^2 dx ds\notag\\
   \le & C\int_0^t  (1+s)^{\frac{\ga-1}{\ga}-\ta}  \int   v_s^2dx ds  +C \int_0^t  (1+s)^{\frac{\ga-1}{\ga}-\ta} \int x^2\bar{\rho}^\ga\left[ \left(\frac{r}{x}-1\right)^2  + \left(r_x-1\right)^2 \right] dx ds \notag\\
   & +   C\int_0^t  (1+s)^{-\frac{1}{\ga}-\ta} \int \bar\rho^{3\ga-2} \mathcal{G}_x^2 dx ds  \le
    C(\ta) \lt(\mathcal{E}(0) + \|  r_{0x}-1\|_{L^\iy}^2 \rt), \label{Aug9-3}
\end{align}
due to \ef{estlem51}, \ef{bt1} and the following estimate:
\begin{align*}
&\int_0^t  (1+s)^{-\frac{1}{\ga}-\ta} \int \bar\rho^{3\ga-2} \mathcal{G}_x^2 dx ds\\
 = &\int_0^t  (1+s)^{-\frac{1}{\ga}-\ta} \lt(\int \bar\rho^{2\ga-1} \mathcal{G}_x^2 dx \rt)^{\frac{1}{\ga}}\lt(\int \bar\rho^{3\ga-1} \mathcal{G}_x^2 dx \rt)^{\frac{\ga-1}{\ga}} ds \\
\le& C  \int_0^t  (1+s)^{-1-\ga\ta} ds \sup_{s\in [0,t]} \int \bar\rho^{2\ga-1} \mathcal{G}_x^2 (x,s) dx
+ \int_0^t   \int \bar\rho^{3\ga-1} \mathcal{G}_x^2 dx ds.
\end{align*}
It gives from \ef{weightrxx} and \ef{Aug9-3} that
\begin{equation}\label{didadida}\begin{split}
 (1+t)^{\frac{\ga-1}{\ga}-\ta} \int \bar\rho^{3\ga-2}    \lt(r_{xx}^2 + \lt|\lt({r}/{x}\rt)_x\rt|^2 \rt)(x,t)dx  \le  C(\ta) \lt(\mathcal{E}(0)+\|r_{0x}-1\|^2_{L^\iy}\rt).
\end{split}\end{equation}
So, \ef{lem5est8/9} follows directly from \ef{didadida}.

{\em Step 2}. In this step, we prove that
\begin{align}\label{didadida'}
\int_0^t (1+s)^{\frac{\ga-1}{\ga}-\ta} \int \lt(v_{x}^2 + \lt|{v}/{x}\rt|^2 \rt)dx ds \le C(\ta) \lt(\mathcal{E}(0)+\|r_{0x}-1\|^2_{L^\iy}\rt).
\end{align}
 Multiplying the square of   \eqref{7-2} by
$ \bar\rho^{2\ga-2} (1+t)^{(\ga-1)/\ga-\ta}$ and integrating the product with respect to the spatial and temporal variables, we have, using \ef{Aug9-3} and \ef{estlem51}, that
\begin{equation*}\label{}\begin{split}
 & \int_0^t (1+s)^{\frac{\ga-1}{\ga}-\ta} \int  \bar\rho^{2\ga-2}   \mathcal{G}_{xs}^2dx ds
 \le  C \int_0^t (1+s)^{\frac{\ga-1}{\ga}-\ta} \int \lt(\bar\rho^{4\ga-2} \mathcal{G}_x^2 + v_s^2 \rt) dx ds
 \\
   &  +C \int_0^t  (1+s)^{\frac{\ga-1}{\ga}-\ta} \int x^2\bar{\rho}^\ga\left[ \left(\frac{r}{x}-1\right)^2  + \left(r_x-1\right)^2 \right] dx ds  \le  C(\ta) \lt(\mathcal{E}(0)+\|r_{0x}-1\|^2_{L^\iy}\rt).
\end{split}\end{equation*}
This, together with \ef{weightvxx} and \ef{Aug9-3}, implies
\begin{equation*}\label{huala}\begin{split}
\int_0^t (1+s)^{\frac{\ga-1}{\ga}-\ta} \int \bar\rho^{4\ga-2} \lt(v_{xx}^2 + \lt|\lt(\frac{v}{x}\rt)_x\rt|^2 \rt)dx ds \le  C(\ta) \lt(\mathcal{E}(0)+\|r_{0x}-1\|^2_{L^\iy}\rt).
\end{split}\end{equation*}
Clearly,
\begin{equation}\label{8.21.4}\begin{split}
\int_0^t (1+s)^{\frac{\ga-1}{\ga}-\ta} \int_0^{1/2} \lt(v_{xx}^2 + \lt|\lt(\frac{v}{x}\rt)_x\rt|^2 \rt)dx ds \le  C(\ta) \lt(\mathcal{E}(0)+\|r_{0x}-1\|^2_{L^\iy}\rt).
\end{split}\end{equation}
Then, it follows from   \ef{hardyorigin} and \ef{estlem51} that
\begin{equation*}\label{}\begin{split}
\int_0^t (1+s)^{\frac{\ga-1}{\ga}-\ta} \int_0^{1/2} \lt(v_{x}^2 + \lt|{v}/{x}\rt|^2 \rt)dx ds \le  C(\ta) \lt(\mathcal{E}(0)+\|r_{0x}-1\|^2_{L^\iy}\rt).
\end{split}\end{equation*}
We use \ef{estlem51} again to obtain \ef{didadida'}.

{\em Step 3}. In this step, we show that
\begin{equation}\label{8/10-1}\begin{split}
& (1+t)^{({\ga-1})/{\ga}-\ta} \int \bar{\rho}   v_t^2 (x,t) dx
+  \int_0^t (1+s)^{({\ga-1})/{\ga}-\ta}  \int  \left( v_{xs}^2
  + (v_s/x)^2 \right)dxds \\
   \le & C(\theta) \lt(\mathcal{E}(0)+\|r_{0x}-1\|^2_{L^\iy}\rt).
 \end{split}
\end{equation}
Differentiating \eqref{nsp1} with respect to $t$ yields
\begin{equation}\label{nsptime}\begin{split}
 &\bar\rho\left( \frac{x}{r}\right)^2  v_{tt}  -2\bar\rho\left( \frac{x}{r}\right)^3 \frac{v}{x}  v_{t}  -\ga \left[  \left(\frac{x^2}{r^2}\frac{\bar\rho}{ r_x}\right)^\ga \left(2\frac{v}{r}+\frac{v_x}{r_x}\right)   \right]_{x} +4 \left( \frac{x}{r}\right)^5 \frac{v}{x}\left(\bar{\rho}^\ga\right)_x    \\
  = & \mu \left(\frac{v_{xt}}{r_x}+2\frac{v_t}{r}\right)_x -  \mu \left( \frac{v_x^2 }{r_x^2} +2 \frac{v^2}{r^2}  \right)_x.
\end{split}
\end{equation}
Let $\psi$ be a non-increasing function defined on $[0,1]$ satisfying
\begin{equation*}\label{}\begin{split}
 \psi=1 \ \ {\rm on } \ \ [0,1/4], \ \ \psi=0 \ \ {\rm on } \ \  [1/2, 1] \ \ {\rm and} \  \ |\psi'|\le 32.
\end{split}
\end{equation*}
Multiplying   equation \eqref{nsptime} by $ \psi v_t$ and integrating the product with respect to  the  spatial variable, one has, using the integration by parts and the boundary condition $v(0,t)=0$ (so $v_t(0,t)=0$), that
 \begin{equation}\label{levt'}\begin{split}
& \frac{d}{dt}\int \frac{1}{2}\bar{\rho} \psi \left(\frac{x}{r}\right)^2  v_t^2 dx
+ \mu  \int  \left(\frac{v_{xt}}{r_x}+2\frac{v_t}{r}\right)   \left(\psi v_t\right)_x dx
=J_1+J_2+J_3 ,
\end{split}
\end{equation}
where
\begin{equation*}\label{}\begin{split}
J_1:=&\int \frac{v}{r} \bar{\rho}\psi\left(\frac{x}{r}\right)^2   v_{t}^2dx + 4 \int \left(\frac{x}{r}\right)^5v  \phi \bar{\rho}\psi v_t dx \le C
 \int_0^{1/2}  \left( v^2 + v_t^2\right)  dx,\\
J_2:=&-\ga \int   \left(\frac{x^2}{r^2}\frac{\bar\rho}{ r_x}\right)^\ga \left(2\frac{v}{r}+\frac{v_x}{r_x}\right)   \left(\psi v_t\right)_x dx,   \\
J_3:= & \mu   \int  \left[  \frac{v_x^2 }{r_x^2}  +2 \frac{v^2}{r^2}    \right] \left(\psi v_t\right)_x dx .
\end{split}
\end{equation*}
The second term on the left-hand side of \eqref{levt'} can be estimated as follows:
 \begin{equation*}\label{}\begin{split}
 &\int  \left(\frac{v_{xt}}{r_x}+2\frac{v_t}{r}\right)  \left(\psi v_t\right)_x dx
 =  \int \psi\frac{v_{xt}^2}{r_x} dx
 +2 \int \frac{\psi}{r} v_t v_{tx}dx + \int \psi' \left(\frac{v_{xt}}{r_x}+2\frac{v_t}{r}\right) v_t  dx \\
 \ge  & \int \psi\frac{v_{xt}^2}{r_x} dx
-\int \lt(\frac{\psi}{r}\rt)' v_t^2dx
-C \int_0^{1/2}\left( {x^2v_{xt}^2}  +   v_t^2\right)  dx
\\
 \ge  & \int \psi\lt[\frac{v_{xt}^2}{r_x} + \frac{r_x v_t^2}{r^2}\rt] dx
-C \int_0^{1/2}\left( {x^2v_{xt}^2}  +   v_t^2\right)  dx.
\end{split}
\end{equation*}
Then,
 \begin{equation*}\label{}\begin{split}
 \frac{d}{dt}\int \frac{1}{2}\bar{\rho}\psi \left(\frac{x}{r}\right)^2  v_t^2 dx
+ \mu \int \psi \left[  \frac{v_{xt}^2}{r_x}
  + \frac{r_x v_t^2}{r^2} \right]dx
\le   C
 \int_0^{1/2}\left( {x^2v_{xt}^2}  +   v_t^2+v^2\right)  dx +  J_2+ J_3.
\end{split}
\end{equation*}
It therefore follows from the Cauchy  inequality that
 \begin{equation}\label{wawa2}\begin{split}
 \frac{d}{dt}\int \frac{1}{2}\bar{\rho}\psi \left(\frac{x}{r}\right)^2  v_t^2 dx
+ \frac{\mu}{2} \int \psi \left[  \frac{v_{xt}^2}{r_x}
  + \frac{r_x v_t^2}{r^2} \right]dx
\le    C
 \int_0^{1/2}\left( {x^2v_{xt}^2}  +   v_t^2+( v/x)^2+ v_x^2 \right)  dx .
\end{split}
\end{equation}
This, together with \ef{lem4est} and \ef{didadida'}, implies that
\begin{equation*}\label{}\begin{split}
 (1+t)^{\frac{\ga-1}{\ga}-\ta} \int_0^{1/4} \bar{\rho}   v_t^2 dx
+ \int_0^t  (1+s)^{\frac{\ga-1}{\ga}-\ta} \int_0^{1/4}   \left( v_{xs}^2
  + \frac{ v_s^2}{x^2} \right)dxds
  \le   C(\ta) \lt(\mathcal{E}(0)+\|r_{0x}-1\|^2_{L^\iy}\rt).
 \end{split}
\end{equation*}
Using \ef{lem4est} again, we obtain \ef{8/10-1}. So, \ef{lem5est'} follows from  \ef{8/10-1} and \ef{didadida'}.

{\em Step 4}. In this step, we prove \ef{lem5est'8/9}. It follows from \ef{7-2},  \ef{8/10-1}, \ef{estlem51} and \ef{Aug9-3} that
\bee\label{}\begin{split}
& \int \bar\rho^{\ga-2} \mathcal{G}_{xt}^2 (x,t)dx \le    C\int \bar{\rho}  v_t^2 (x,t) dx
+C\int \bar\rho^{3\ga-2} \mathcal{G}_{x}^2(x,t)dx  \\
&   + C \int x^2 \bar{\rho}^\ga \left[ \left(\frac{r}{x}-1\right)^2  + \left(r_x-1\right)^2 \right](x,t) dx  \le   C(\ta) (1+t)^{-\frac{\ga-1}{\ga}+\ta} \lt(\mathcal{E}(0)+\|r_{0x}-1\|^2_{L^\iy}\rt).
\end{split}\eee
This, together with \ef{weightvxx} and \ef{Aug9-3}, implies
\begin{equation}\label{hualaAug10}\begin{split}
 (1+t)^{\frac{\ga-1}{\ga}-\ta} \int \bar\rho^{3\ga-2} \lt(v_{xx}^2 + \lt| ( {v}/{x} )_x\rt|^2 \rt)(x,t)dx  \le  C(\ta) \lt(\mathcal{E}(0)+\|r_{0x}-1\|^2_{L^\iy}\rt).
\end{split}\end{equation}
So, \ef{lem5est'8/9} follows from \ef{hualaAug10}.

\hfill $\Box$

\noindent{\bf Proof of Theorem \ref{mainthm1}}. The proof of the first part of Theorem \ref{mainthm1} follows from the local existence and uniqueness results of strong solutions, of which a sketched proof is given in  Appendix, Part I.  So, the global well-posedness of strong solutions with the estimate \ef{keyconclusion}  can be shown by Lemma \ref{lem5}, together with the equivalence of  $\mathcal{E}(t)$ and $\mathfrak{E}(t)$  shown in \ef{equivalence} and the lower-order estimates obtained in Subsection \ref{sec3.3}, through the standard continuation argument .

Moreover, \ef{Aug23-2} and \ef{Aug23-1} follow from \ef{8/23/1} and \ef{growof2nd8.22} in Lemma \ref{lem10}, which will be proved later.

\hfill$\Box$

To complete  the proof of part $i)$ of Theorem \ref{mainthm}, it suffices to show the following Lemma.

\begin{lem}\label{lem312} For the global strong solution obtained in Theorem \ref{mainthm1}, there exist positive constants  $C(\theta)$ and $C(\ta,a)$ independent of $t$ such that for  any $\theta\in (0,  \ {2(\ga-1)}/({3\ga}))$ and $a\in (0,1)$,
\begin{align}
&(1+t)^{({\ga-1})/{\ga}-\ta} \lt\|\lt(r_x-1, {r}/{x}-1,  v_x, v/{x}\rt)(\cdot,t)\rt\|^2 \le C(\theta)  \mathfrak{E}(0), \label{8/10-2}\\
&(1+t)^{({\ga-1})/{\ga}-\ta} \lt\|\lt(r_x-1, {r}/{x}-1,  v_x, v/{x}\rt)(\cdot,t)\rt\|^2_{H^1([0,a])} \le C(\theta,a)   \mathfrak{E}(0), \label{8/10-3}\\
&(1+t)^{2({\ga-1})/{ \ga}- \ta } \lt\|(r-x)(\cdot,t)\rt\|_{L^\iy}^2
  \le C(\ta)\mathfrak{E}(0),\label{8/11-1}\\
&(1+t)^{({3\ga-2})/({2\ga})- {\ta}} \lt\|(v, x v_x)(\cdot,t)\rt\|_{L^\iy}^2 \le C(\ta) \mathfrak{E}(0),\label{8/11-2}\\
&  (1+t)^{({ \ga-1})/\ga- {\ta}} \lt\| \lt(v_x, v/x\rt)(\cdot,t)\rt\|_{L^\iy}^2 \le C(\ta) \mathfrak{E}(0),\label{8/11-3}\\
&(1+t)^{({\ga-1})/{ \ga} - \ta} \lt\|\bar\rho^{ ({3\ga-2})/4}  (r_x-1, r/x-1)(\cdot,t)\rt\|_{L^\iy}^2 \le C(\ta) \mathfrak{E}(0).\label{8/11-4}
\end{align}
\end{lem}
{\em Proof}. It follows from  \ef{hardyorigin},     \ef{estlem51}, \ef{lem5est8/9} and \ef{lem5est'8/9}  that
\begin{align}
\lt\|\lt(r_x-1, {r}/{x}-1,  v_x, v/{x}\rt)(\cdot,t)\rt\|_{L^2  \lt(\lt[0, {1}/{2}\rt]\rt)}^2 \le C \|(xr_{x}-x, r-x,  xv_{x}, v)(\cdot,t)\|^2_{L^2([0,1/2])} \notag\\
+ C \|(xr_{xx}, x(r/x)_x,  xv_{xx}, x(v/x)_x)(\cdot,t)\|^2_{L^2([0,1/2])}
  \le   C(\theta) (1+t)^{-({\ga-1})/{\ga}+\ta}  \mathfrak{E}(0),\notag
\end{align}
which, together with \ef{estlem51}, gives \ef{8/10-2}. Clearly, \ef{8/10-3} follows from \ef{8/10-2}, \ef{lem5est8/9} and \ef{lem5est'8/9};  \ef{8/11-1} from  \ef{girl},  \ef{estlem51} and \ef{8/10-2}; the estimate for $v$ in \ef{8/11-2} from  \ef{f2},  \ef{estlem51} and \ef{8/10-2}.  Due to  \ef{7-2},  \ef{estlem51}  and \ef{Aug9-3}, we have
\bee\label{}\begin{split}
\| x \mathcal{G}_{xt}(\cdot,t)\|^2 \le & C \lt(\|x\bar\rho^\ga\mathcal{G}_x(\cdot,t)\|^2 + \|x\bar\rho  v_t(\cdot,t)\|^2 + \|\bar\rho(r-x, xr_x-x)(\cdot,t)\|^2   \rt) \\ \le &  C(\theta) (1+t)^{-({\ga-1})/{\ga}+\ta}    \mathfrak{E}(0),
\end{split}\eee
which implies, using \ef{f1}, \ef{estlem51} and \ef{8/10-2}, that
$$(1+t)^{({3\ga-2})/({2\ga})- {\ta}} \lt\| xv_x(\cdot,t)\rt\|_{L^\iy}^2   \le C(\ta) \mathfrak{E}(0).$$
This verifies \ef{8/11-2}. \ef{8/11-3} follows from \ef{8/11-2}, \ef{8/10-3} and the fact $\|\cdot\|_{L^\iy([0,1/2])}\le C \|\cdot\|_{H^1([0,1/2])}$.

It follows from \ef{estlem51} that
\begin{align}
   \lt\|x\bar\rho^{(\ga-1)/2}(r_x-1)(\cdot,t)\rt\| \le  &
 \lt\|x(r_x-1)(\cdot,t)\rt\|^{1/\ga}  \lt\|x\bar\rho^{\ga/2}(r_x-1)(\cdot,t)\rt\|^{(\ga-1)/\ga} \notag\\
  \le  & C(\ta)  (1+t)^{-(\ga-1)/\ga +\ta/2}  \mathfrak{E}(0) , \notag
\end{align}
which implies,  using \ef{woman}, the H$\ddot{o}$lder inequality and \ef{lem5est}, that
\begin{align}
 \lt\|x^{{3}/{2}}\bar\rho^{(3\ga-2)/4}(r_x-1)(\cdot, t)\rt\|_{L^\iy}^2 \le C\lt\|x\bar\rho^{(\ga-1)/2}(r_x-1)(\cdot, t)\rt\|  \lt\|\bar\rho^{(2\ga-1)/2}r_{xx}(\cdot, t)\rt\|  \notag\\
 +  C  \lt\|x(r_x-1)(\cdot, t)\rt\|^2
  \le  C(\ta)  (1+t)^{-(\ga-1)/\ga +\ta}  \mathfrak{E}(0) . \notag
\end{align}
This, together with \ef{8/10-3}, gives
$$
 \lt\|\bar\rho^{(3\ga-2)/4}(r_x-1)(\cdot, t)\rt\|_{L^\iy}^2
  \le    C(\ta)  (1+t)^{-(\ga-1)/\ga +\ta}  \mathfrak{E}(0) . $$
Similarly, we can use \ef{woman.new} to get
$$
 \lt\|\bar\rho^{(3\ga-2)/4}(r/x-1)(\cdot, t)\rt\|_{L^\iy}^2
  \le    C(\ta)  (1+t)^{-(\ga-1)/\ga +\ta}  \mathfrak{E}(0) . $$
This finishes the proof of \ef{8/11-4}.

\hfill $\Box$

\subsubsection{Part II: faster decay }\label{sec3.4.3}
In this subsection, we prove part $ii)$ of Theorem \ref{mainthm}  under the assumption
\be\label{finitenessofF} \mathfrak{F}_\alpha(0)<\infty, \ \  \ \  \alpha \in (0, \ga). \ee
The estimates in this subsection   are for the global  strong solution of \ef{419} as stated in Theorem \ref{mainthm1}.

To obtain the faster time decay estimates of the higher-order norms, we rewrite  equation \eqref{7-2} in the form of
 \begin{equation}\label{x1}\begin{split}
 \mathfrak{P}(x,t) + \mu \mathcal{G}_{xt} = \frac{x^2}{r^2} \bar{\rho}  v_t , \ \ {\rm where} \ \  \mathfrak{P}(x,t):=\ga \left(\frac{x^2 \bar{\rho}}{r^2 r_x } \right)^{\ga } \mathcal{G}_x + \left[ \left(\frac{x^2 }{r^2 r_x } \right)^{\ga }  -\left(\frac{x }{r } \right)^{4}  \right]  x \phi \bar\rho.
\end{split}\end{equation}
It should be noted that
\be\label{8.19.1}
 \mathfrak{P}_t =\ga \bar{\rho}^\ga \left(\frac{x^2 }{r^2 r_x } \right)^{\ga } \mathcal{G}_{xt} + \mathfrak{P}_1, \ \ {\rm where} \ \  \mathfrak{P}_1 : = \ga  \bar{\rho}^\ga \lt[\left(\frac{x^2 }{r^2 r_x } \right)^{\ga } \rt]_t\mathcal{G}_{x} - \left[ \left(\frac{x^2 }{r^2 r_x } \right)^{\ga } -\left(\frac{x }{r } \right)^{4}  \right]_t  x \phi \bar\rho.
\ee
This equation is convenient for us to derive  the time decay estimates for $r_{xx}$ and  $v_{xx}$ with weights.

\begin{lem}\label{lem8.19} Let $\alpha\in (0, \ga)$ and $\mathfrak{F}_\alpha(0)<\infty$. For the global strong solution obtained in Theorem \ref{mainthm1}, there exist positive constants  $C(\alpha)$ and $C(\ta,\alpha)$ independent of $t$ such that  for any $0<\theta< \min\{2(\ga-1)/(3\ga), \ 2(\ga-\alpha)/\ga\}$,
\begin{align}
& \left\|\bar\rho^{(2\ga-1-\alpha)/2}\lt(r_{xx}, (r/x)_x\rt)(\cdot, t) \rt\|^2
 \le C(\alpha) \mathfrak{F}_\alpha (0),   \label{8.20-1} \\
&(1+t)^{({\ga-1})/{ \ga} - \ta} \lt\|\bar\rho^{ ({3\ga-2-\alpha})/4}  (r_x-1, r/x-1)(\cdot,t)\rt\|_{L^\iy}^2 \le C(\ta, \alpha)\mathfrak{F}_\alpha (0).\label{8/11-4.Aug21}
\end{align}
\end{lem}
{\em Proof}.  In a similar way to deriving \ef{bt1}, we have
\begin{equation}\label{bt1.8.19}\begin{split}
&  \int  \bar\rho^{2\ga-1-\alpha}   \mathcal{G}_{x}^2 (x,t)dx + \int_0^t  \int  \bar\rho ^{3\ga-1-\alpha} \mathcal{G}_x^2 dx ds
\le    C\int_0^t   \int   v_s^2dxds  \\
& \qquad  +C  \int_0^t  \int x^2\bar\rho^{\ga+1-\alpha} \left[ \left({r}/{x}-1\right)^2  + \left(r_x-1\right)^2 \right] dx ds
   \le  C(\alpha) \mathfrak{F}_\alpha(0),
\end{split}\end{equation}
which, together with \ef{weightrxx}, gives \ef{8.20-1}.  The proof for \ef{8/11-4.Aug21} is the same as that for
 \ef{8/11-4}. We omit the detail here.

 \hfill $\Box$

\begin{lem}\label{lem41} Let $\alpha\in [\ga-1, \ga)$  and $\mathfrak{F}_\alpha(0)<\infty$. For the global strong solution obtained in Theorem \ref{mainthm1}, there exist positive constants  $C(\alpha,\ta)$ and $C(\alpha,\ta, a)$ independent of $t$ such that  for any
$0<\theta< \min\{2(\ga-1)/(3\ga), \ 2(\ga-\alpha)/\ga\}$ and $a\in (0,1)$,
 \begin{align}
  &(1+t)^{\frac{\kappa }{2}+\frac{4\ga-3}{2\ga}-\frac{3}{2}\ta}   \left\| \bar\rho^{ {1}/{2}} v_t(\cdot, t) \rt\|^2  +\int_0^t (1+s)^{\frac{\kappa }{2}+\frac{4\ga-3}{2\ga}-\frac{3}{2}\ta}  \lt\|\lt( v_{x},  {v}/{x} , v_{sx},  {v_s}/{x}  \rt)(\cdot,s)\rt\|^2  ds\notag\\
& \quad  \le C(\alpha, \ta) \lt(
  {\mathfrak{E}}(0) +1 \rt) {\mathfrak{F}}_\alpha(0),   \label{8-21-2}\\
& (1+t)^{\min\lt\{\frac{\kappa }{2}+\frac{4\ga-3}{2\ga}-\frac{3}{2}\ta, \ \frac{\kappa }{4}+\frac{10\ga-9}{4\ga}-\frac{9}{4}\ta\rt\}} \left\| \lt(r_{xx}, (r/x)_x, v_{xx}, (v/x)_x \rt)(\cdot, t) \rt\|^2_{L^2([0,a])}
\notag\\
& \quad \le C(\alpha, \ta,a ) \lt(
  {\mathfrak{E}}(0) +1 \rt) {\mathfrak{F}}_\alpha(0). \label{8-21-1}
\end{align}
Here $\kappa=0$ when $\alpha = \ga-1$, and  $\kappa = (1/\ga) \min\{  {\alpha-(\ga-1)} , \   {\ga-1} \} -\theta$ when $\alpha\in (\ga-1, \ga)$.
\end{lem}
{\em Proof}. The proof consists of four steps.

{\em Step 1}. In this step, we prove
\begin{align}
&(1+t)^{\ka} \int \bar\rho^{\ga}    \mathcal{G}_{x}^2(x,t)dx +\int_0^t (1+s)^{\ka} \int \bar\rho^{2\ga} \mathcal{G}_x^2 dx ds\le
    C(\ta,\alpha)  \mathfrak{F}_\alpha(0), \label{8.20.1}\\
& \int_0^t (1+s)^{\ka} \int  \lt(\mathcal{G}_{sx}^2 + v_x^2 +|v/x|^2  \rt)dx ds\le
    C(\ta,\alpha)  \mathfrak{F}_\alpha(0).  \label{8.20.2}
\end{align}

When $\alpha=\ga-1$, \ef{bt1.8.19} implies \ef{8.20.1}. When $\alpha\in (\ga-1, \ga)$, we use the same way as that for the derivation of \ef{Aug9-3} to obtain
\begin{align}
 & (1+t)^{\ka} \int \bar\rho^{\ga}    \mathcal{G}_{x}^2(x,t)dx +\int_0^t (1+s)^{\ka} \int \bar\rho^{2\ga} \mathcal{G}_x^2 dx ds\notag\\
   \le & C\int_0^t  (1+s)^{\kappa}  \int   v_s^2dx ds  +C \int_0^t  (1+s)^{\ka} \int x^2\bar{\rho}^\ga\left[ \left(\frac{r}{x}-1\right)^2  + \left(r_x-1\right)^2 \right] dx ds \notag\\
   & +   C\int_0^t  (1+s)^{\ka-1} \int \bar\rho^{\ga} \mathcal{G}_x^2 dx ds \le
    C(\ta,\alpha)  \mathfrak{F}_\alpha(0),  \notag
\end{align}
due to
$$
   \int \bar\rho^{\ga} \mathcal{G}_x^2 dx \le   \lt( \int \bar\rho^{2\ga-1-\alpha} \mathcal{G}_x^2 dx \rt)^{(2\ga-1-\alpha)/\ga}\lt( \int \bar\rho^{3\ga-1-\alpha} \mathcal{G}_x^2 dx \rt)^{( \alpha + 1 -\ga)/\ga}
$$
and
\begin{align*}
  \int_0^t  (1+s)^{\ka-1} & \int \bar\rho^{\ga} \mathcal{G}_x^2 dx ds \le     C \int_0^t  \int  \bar\rho ^{3\ga-1-\alpha} \mathcal{G}_x^2 dx ds \\
 &   +C \int_0^t  (1+s)^{-1-\ga\ta/(2\ga-1-\alpha)} ds \sup_{s\in [0,t]} \int \bar\rho^{2\ga-1-\alpha} \mathcal{G}_x^2(x,s) dx.
\end{align*}
This verifies \ef{8.20.1}.  It follows from \ef{7-2},  \ef{8.20.1} and  \ef{estlem51} that
$$ \int_0^t (1+s)^{\ka} \int  \mathcal{G}_{sx}^2  dx ds\le
    C(\ta,\alpha)  \mathfrak{F}_\alpha(0) ,$$
which, together with \ef{didadida'}, gives    \ef{8.20.2}.

{\em Step 2}. In this step, we prove that
\begin{equation}\label{f3}\begin{split}
&(1+t)^{\frac{\kappa }{2}+\frac{4\ga-3}{2\ga}-\frac{3}{2}\ta}\int  x^2\bar\rho^{4\ga-2}\mathcal{ G}_{x}^2(x,t) dx
+\int_0^t (1+s)^{\frac{\kappa }{2}+\frac{4\ga-3}{2\ga}-\frac{3}{2}\ta}\int x^2\bar\rho^{3\ga-2} \mathcal{G}_{xs}^2 dx ds\\
& \le C(\alpha, \ta) \widetilde{\mathfrak{F}}_\alpha(0), \ \  \ \  {\rm  where } \ \ \widetilde{\mathfrak{F}}_\alpha(0):= {\mathfrak{F}}_\alpha(0) +
  {\mathfrak{E}}(0) {\mathfrak{F}}_\alpha(0) .
\end{split}\end{equation}

Multiplying  \ef{x1} by $x^2\bar\rho^{2\ga-2} \mathfrak{P}_t$ and integrating the resulting equation give that
 $$
\frac{1}{2}\frac{d}{dt}\int x^2\bar\rho^{2\ga-2}\mathfrak{P}^2 dx + \mu\int x^2 \bar\rho^{2\ga-2}\mathfrak{P}_t \mathcal{G}_{xt} dx
    = \int x^2 \bar\rho^{2\ga-1}\mathfrak{P}_t \frac{x^2}{r^2} v_t dx,
$$
which implies, using \ef{8.19.1}, that
\begin{equation*}\label{}\begin{split}
&\frac{1}{2}\frac{d}{dt}\int x^2\bar\rho^{2\ga-2}\mathfrak{P}^2 dx
+\mu \ga\int x^2\bar\rho^{3\ga-2} \left(\frac{x^2 }{r^2 r_x } \right)^{\ga } \mathcal{G}_{xt}^2 dx \\
    = &\ga \int x^2 \bar\rho^{3\ga-1} \left(\frac{x^2 }{r^2 r_x } \right)^{\ga } \mathcal{G}_{xt}  \frac{x^2}{r^2}   v_t dx + \int x^2 \bar\rho^{2\ga-1}\mathfrak{P}_{1} \frac{x^2}{r^2}    v_t dx
    -\mu\int x^2 \bar\rho^{2\ga-2}\mathfrak{P}_{1} \mathcal{G}_{xt} dx .
\end{split}\end{equation*}
Thus, one has
\begin{equation*}\label{}\begin{split}
\frac{1}{2}\frac{d}{dt}\int x^2\bar\rho^{2\ga-2}\mathfrak{P}^2 dx
+\frac{\mu \ga}{2}\int x^2\bar\rho^{3\ga-2} \left(\frac{x^2 }{r^2 r_x } \right)^{\ga } \mathcal{G}_{xt}^2 dx\le C \int \lt( x^2 \bar\rho^{\ga-2}\mathfrak{P}_{1}^2 +v_t^2 \rt) dx \\
 \le  C \int \lt(x^2 v_x^2 +v^2\rt)\bar{\rho}^{3\ga-2}\mathcal{G}_x^2 dx + C \int \lt(x^2 v_x^2 +v^2+v_t^2\rt) dx.
\end{split}\end{equation*}
Combining this with \ef{Aug9-3} shows that
\begin{equation}\label{ppqq}\begin{split}
&\frac{1}{2}\frac{d}{dt}\int x^2\bar\rho^{2\ga-2}\mathfrak{P}^2 dx
+\frac{\mu \ga}{2}\int x^2\bar\rho^{3\ga-2} \left(\frac{x^2 }{r^2 r_x } \right)^{\ga } \mathcal{G}_{xt}^2 dx \\
\le & C\lt( \lt\|x v_x\rt\|_{L^\iy}^2 +C\lt\|v\rt\|_{L^\iy}^2\rt) \int \bar{\rho}^{3\ga-2}\mathcal{G}_x^2 dx + C \int \lt(x^2 v_x^2 +v^2 + v_t^2\rt) dx \\
\le & C(\ta) \mathfrak{E}(0) (1+t)^{-(\ga-1)/\ga+\ta} \lt( \lt\|x v_x\rt\|_{L^\iy}^2 +\lt\|v\rt\|_{L^\iy}^2\rt)  + C \int \lt(x^2 v_x^2 +v^2 + v_t^2\rt) dx.
\end{split}\end{equation}
It then follows from  \ef{f1} and \ef{f2}  that
\begin{equation*}\label{}\begin{split}
&\frac{1}{2}\frac{d}{dt}\int x^2\bar\rho^{2\ga-2}\mathfrak{P}^2 dx
+\frac{\mu \ga}{2}\int x^2\bar\rho^{3\ga-2} \left(\frac{x^2 }{r^2 r_x } \right)^{\ga } \mathcal{G}_{xt}^2 dx\\
\le & C(\ta) \mathfrak{E}(0) (1+t)^{\frac{\kappa }{2}-\frac{4\ga-3}{2\ga}+\frac{3}{2}\ta}\int \lt( \mathcal{G}_{tx}^2 +v_x^2 + \lt({v}/{x}\rt)^2 \rt)dx  \\
&+  C(\ta) \mathfrak{E}(0) (1+t)^{\frac{1}{2\ga}-\frac{\kappa }{2}+\frac{\ta}{2}} \int \lt(x^2 v_x^2 +v^2\rt)  dx + C \int \lt(x^2 v_x^2 +v^2 + v_t^2\rt) dx.
\end{split}\end{equation*}
Multiplying the inequality above by $(1+t)^{\frac{\kappa }{2}+\frac{4\ga-3}{2\ga}-\frac{3}{2}\ta}$ and integrating the product give that
\begin{equation}\label{8.21.1}\begin{split}
&(1+t)^{\frac{\kappa }{2}+\frac{4\ga-3}{2\ga}-\frac{3}{2}\ta}\int x^2\bar\rho^{2\ga-2}\mathfrak{P}^2(x,t) dx
+\int_0^t (1+s)^{\frac{\kappa }{2}+\frac{4\ga-3}{2\ga}-\frac{3}{2}\ta}\int x^2\bar\rho^{3\ga-2} \mathcal{G}_{xs}^2 dx ds \\
\le & C \int x^2\bar\rho^{2\ga-2}\mathfrak{P}^2(x,0) dx+  C(\ta) \mathfrak{E}(0)  \int_0^t (1+s)^{\kappa}\int \lt( \mathcal{G}_{sx}^2 +v_x^2 + (v/x)^2 \rt)dx ds \\
& +  \int_0^t \lt[C(\ta) \mathfrak{E}(0)(1+s)^{\frac{2\ga-1}{\ga}-{\ta}} + C (1+s)^{\frac{\kappa }{2}+\frac{4\ga-3}{2\ga}-\frac{3}{2}\ta}\rt]\int \lt(x^2 v_x^2 +v^2 + v_s^2 \rt) dxds    \\
& + C \int_0^t (1+s)^{(1+s)^{\frac{\kappa }{2}+\frac{2\ga-3}{2\ga}-\frac{3}{2}\ta}}\int x^2\bar\rho^{2\ga-2}\mathfrak{P}^2 dx ds.
\end{split}\end{equation}
Note that
\begin{align}
&\int x^2\bar\rho^{2\ga-2}\mathfrak{P}^2 dx \le C \int x^2\bar\rho^{4\ga-2}\mathcal{ G}_{x}^2 dx + C \int x^2 \bar\rho^\ga \left[ \left(\frac{r}{x}-1\right)^2  + \left(r_x-1\right)^2 \right]dx , \notag\\
& \frac{\kappa }{2}+\frac{2\ga-3}{2\ga}-\frac{3}{2}\ta  \le \frac{\ga-1}{\ga}-{\ta}, \ \ \ \   \frac{\kappa }{2}+\frac{4\ga-3}{2\ga}-\frac{3}{2}\ta \le \frac{2\ga-1}{\ga}-{\ta}    . \label{8.21.3}
 \end{align}
Then,  it follows from \ef{8.21.1}, \ef{8.20.2}, \ef{estlem51}  and  \ef{Aug9-3} that
$$
(1+t)^{\frac{\kappa }{2}+\frac{4\ga-3}{2\ga}-\frac{3}{2}\ta}\int x^2\bar\rho^{2\ga-2}\mathfrak{P}^2 dx
+\int_0^t (1+s)^{\frac{\kappa }{2}+\frac{4\ga-3}{2\ga}-\frac{3}{2}\ta}\int x^2\bar\rho^{3\ga-2} \mathcal{G}_{xs}^2 dx ds \le C(\alpha, \ta) \widetilde{\mathfrak{F}}_\alpha(0).
$$
 This, together with \ef{estlem51}, implies \ef{f3}.

 {\em Step 3}. In this step, we prove that
\begin{align}
& \int_0^t (1+s)^{\frac{\kappa }{2}+\frac{4\ga-3}{2\ga}-\frac{3}{2}\ta}\int \lt( v_x^2 +(v/x)^2\rt) dx ds  \le C(\alpha, \ta) \widetilde{\mathfrak{F}}_\alpha(0)  , \label{8.21.2}\\
& (1+t)^{\frac{\kappa }{2}+\frac{4\ga-3}{2\ga}-\frac{3}{2}\ta} \int \bar{\rho}   v_t^2 (x,t) dx
+  \int_0^t (1+s)^{\frac{\kappa }{2}+\frac{4\ga-3}{2\ga}-\frac{3}{2}\ta}  \int  \left( v_{xs}^2
  + (v_s/x)^2 \right)dxds \notag\\
&   \le   C(\alpha, \ta) \widetilde{\mathfrak{F}}_\alpha(0)  . \label{8/10-1Aug21}
 \end{align}
As  a consequence of \ef{8.21.2} and \ef{8/10-1Aug21}, we get \ef{8-21-2}.

It follows from \ef{gjvx}, \ef{hardyorigin} and \ef{tlg3} that
\begin{align}
& \int \lt( v_x^2 +(v/x)^2\rt) dx \le 4 \int_0^{1/2} \mathcal{G}_{t}^2 dx     + 4 \int_{1/2}^1 \mathcal{G}_{t}^2 dx
\le  C \int_0^{1/2} x^2 \lt( \mathcal{G}_{t}^2 + \mathcal{G}_{tx}^2 \rt)dx    + 16 \int_{1/2}^1 x^2 \mathcal{G}_{t}^2 dx \notag\\
& \le C  \int x^2 \mathcal{G}_{t}^2 dx +  C \int_0^{1/2} x^2  \bar\rho^{3\ga-2}  \mathcal{G}_{tx}^2  dx
 \le   C  \int \lt(x^2 v_x^2 + v^2\rt)  + C \int x^2  \bar\rho^{3\ga-2}  \mathcal{G}_{tx}^2  dx. \label{hh8.21}
\end{align}
This, together with \ef{f3}, \ef{8.21.3} and \ef{estlem51}, gives \ef{8.21.2}.  With the aid of \ef{8.21.2}, we can use the same way as to the derivation of \ef{8/10-1} to obtain \ef{8/10-1Aug21}.

{\em Step 4}. In this step, we prove
\begin{align}
&(1+t)^{\min\lt\{\frac{\kappa }{2}+\frac{4\ga-3}{2\ga}-\frac{3}{2}\ta, \ \frac{\kappa }{4}+\frac{10\ga-9}{4\ga}-\frac{9}{4}\ta\rt\}}\int  \bar\rho^{4\ga-2}\mathcal{ G}_{x}^2(x,t) dx
 \le C(\alpha, \ta) \widetilde{\mathfrak{F}}_\alpha(0) ,\label{f3.8.21} \\
&(1+t)^{\min\lt\{\frac{\kappa }{2}+\frac{4\ga-3}{2\ga}-\frac{3}{2}\ta, \ \frac{\kappa }{4}+\frac{10\ga-9}{4\ga}-\frac{9}{4}\ta\rt\}}\int  \bar\rho^{2\ga-2}\mathcal{ G}_{x t}^2(x,t) dx
 \le C(\alpha, \ta) \widetilde{\mathfrak{F}}_\alpha(0) .\label{8.21.5}
\end{align}
With \ef{f3.8.21} and \ef{8.21.5}, we can obtain \ef{8-21-1}  by use of  \ef{weightrxx} and \ef{weightvxx}.

Let $\bar{\psi}$ be a non-increasing function defined on $[0,1]$ satisfying
$$ \bar\psi=1 \ \ {\rm on } \ \ [0,1/8], \ \ \bar\psi=0 \ \ {\rm on } \ \  [1/4,1] \ \ {\rm and} \  \ |\bar\psi'|\le 32.$$
  Following  the derivation of \ef{ppqq}, one can obtain
\begin{equation*}\label{}\begin{split}
&\frac{1}{2}\frac{d}{dt}\int \bar\psi\bar\rho^{2\ga-2}\mathfrak{P}^2 dx
+\frac{\mu \ga}{2}\int \bar\psi \bar\rho^{3\ga-2} \left(\frac{x^2 }{r^2 r_x } \right)^{\ga } \mathcal{G}_{xt}^2 dx\\
\le & C(\ta) \mathfrak{E}(0)(1+t)^{-(\ga-1)/\ga+\ta}  \lt\|(v_x, v/x)\rt\|_{L^\iy\lt(\lt[0,1/4\rt]\rt)}^2 + C \int \lt( v_x^2 + (v/x)^2+ v_t^2\rt) dx.
\end{split}\end{equation*}
In view of \ef{nj2} and \ef{nj2.new}, we see
\begin{equation*}\label{}\begin{split}
  &\lt\|v_x\rt\|_{L^\iy\lt(\lt[0,1/4\rt]\rt)}^2 +\lt\|v/x\rt\|_{L^\iy\lt(\lt[0,1/4\rt]\rt)}^2 \\
 \le & C \int \left( v_{x}^2 + \left| {v}/{x}\right|^2 \right) dx
 + C \lt[\int \left( v_{x}^2 + \left| {v}/{x}\right|^2 \right) dx\rt]^{1/2}\lt[\int_0^{1/2}  \left(  v_{xx}^2 +  \left|\left( {v}/{x}\right)_x\right|^2\rt)dx \rt]^{1/2} ,
\end{split}\end{equation*}
which implies,
\begin{equation*}\label{}\begin{split}
&\frac{1}{2}\frac{d}{dt}\int \bar\psi\bar\rho^{2\ga-2}\mathfrak{P}^2 dx
+\frac{\mu \ga}{2}\int \bar\psi \bar\rho^{3\ga-2} \left(\frac{x^2 }{r^2 r_x } \right)^{\ga } \mathcal{G}_{xt}^2 dx\\
\le & C(\ta) \mathfrak{E}(0)(1+t)^{\frac{5-6\ga}{4\ga}-\frac{\ka}{4}+\frac{5}{4}\ta}    \int_0^{1/2} \lt( v_{xx}^2 +|( v/x)_x|^2 \rt) dx \\
& + C(\ta)\mathfrak{E}(0) (1+t)^{\frac{3-2\ga}{4\ga}+\frac{\ka}{4}+\frac{3}{4}\ta}  \int \lt( v_x^2 + (v/x)^2 \rt) dx
+ C \int \lt( v_x^2 + (v/x)^2+ v_t^2\rt) dx.
\end{split}\end{equation*}
Similar to \ef{f3}, one can use \ef{8.21.2} and \ef{8.21.4} to obtain
\begin{equation*}\label{}\begin{split}
&(1+t)^{\frac{\kappa }{4}+\frac{10\ga-9}{4\ga}-\frac{9}{4}\ta}\int  \bar\psi \bar\rho^{4\ga-2}\mathcal{ G}_{x}^2(x,t) dx
+\int_0^t (1+s)^{\frac{\kappa }{4}+\frac{10\ga-9}{4\ga}-\frac{9}{4}\ta}\int \bar\psi\bar\rho^{3\ga-2} \mathcal{G}_{xs}^2 dx ds \\
& \le C(\alpha, \ta) \widetilde{\mathfrak{F}}_\alpha(0),
\end{split}\end{equation*}
due to
$$\frac{\kappa }{4}+\frac{10\ga-9}{4\ga}-\frac{9}{4}\ta  \le \frac{2\ga-1}{\ga}-{\ta} .$$
This, together with \ef{f3}, implies \ef{f3.8.21}.  Finally, we can use \ef{7-2},  \ef{8/10-1Aug21}, \ef{estlem51} and \ef{f3.8.21} to show \ef{8.21.5}.

\hfill$\Box$

\begin{lem}\label{lem44}  Let $\alpha\in [\ga-1, \ga)$  and $\mathfrak{F}_\alpha(0)<\infty$. For the global strong solution obtained in Theorem \ref{mainthm1}, there exist positive constants  $C(\alpha,\ta)$ and $C(\alpha,\ta, a)$ independent of $t$ such that  for any
$0<\theta< \min\{2(\ga-1)/(3\ga), \ 2(\ga-\alpha)/\ga\}$ and $a\in (0,1)$,
\begin{align}
&(1+t)^{ \frac{\kappa }{2}+\frac{4\ga-3}{2\ga}-\frac{3}{2}\ta }   \lt\|\lt(  r_x-1, r/x-1 \rt)(\cdot,t)\rt\|_{L^2([0,a])}^2   \le C(  \alpha, \ta, a)\widetilde{\mathfrak{F}}_\alpha(0), \label{Aug21.3}\\
&(1+t)^{ \frac{\kappa }{2}+\frac{4\ga-3}{2\ga}-\frac{3}{2}\ta }   \lt\|\lt(   v_x, v/{x}\rt)(\cdot,t)\rt\|^2   \le C(  \alpha, \ta)\widetilde{\mathfrak{F}}_\alpha(0), \label{Aug21.1}\\
&(1+t)^{\min\lt\{\frac{\kappa }{2}+\frac{4\ga-3}{2\ga}-\frac{3}{2}\ta, \ \frac{\kappa }{4}+\frac{10\ga-9}{4\ga}-\frac{9}{4}\ta\rt\}} \lt\|\lt(r_x-1, {r}/{x}-1,  v_x, v/{x}\rt)(\cdot,t)\rt\|^2_{H^1([0,a])} \notag\\
& \quad \le C(  \alpha, \ta, a)\widetilde{\mathfrak{F}}_\alpha(0), \label{Aug21.2}\\
&(1+t)^{\frac{8\ga-5}{4\ga}+\frac{\ka}{4}-\frac{5}{4}\ta } \lt\|v(\cdot,t)\rt\|_{L^\iy}^2 \le C(  \alpha, \ta)\widetilde{\mathfrak{F}}_\alpha(0),\label{Aug21.4}\\
& (1+t)^{\frac{1}{2} b_1 }
\|xv_x(\cdot,t)\|^2_{L^\iy} +  (1+t)^{\frac{1}{2}\min\{b_1,  b_2\}  } \lt\| \lt( v_x, v/x\rt) (\cdot,t)\rt\|_{L^\iy}^2 \le C(  \alpha, \ta)\widetilde{\mathfrak{F}}_\alpha(0),\label{Aug21.5}\\
&(1+t)^{\frac{\ka}{2}+ \frac{2\ga-1}{2\ga}-\frac{\ta}{2}} \lt\|\bar\rho^{ \ga/2 }  (r_x-1, r/x-1)(\cdot,t)\rt\|_{L^\iy}^2 \le C(  \alpha, \ta)\widetilde{\mathfrak{F}}_\alpha(0).\label{Aug21.6}
\end{align}
Here $\widetilde{\mathfrak{F}}_\alpha(0)= {\mathfrak{F}}_\alpha(0) +
  {\mathfrak{E}}(0) {\mathfrak{F}}_\alpha(0)$,  $\kappa=0$ when $\alpha = \ga-1$,   $\kappa = (1/\ga) \min\{  {\alpha-(\ga-1)} , \   {\ga-1} \} -\theta$ when $\alpha\in (\ga-1, \ga)$,
\begin{align}
  b_1=&\min\lt\{ \max\lt\{  \lt(\frac{\kappa }{2}+\frac{4\ga-3}{2\ga}-\frac{3}{2}\ta\rt)\frac{\aa+1}{2\ga-1+\aa}, \  \frac{3 }{2}\kappa +\frac{2\ga-1}{2\ga}-\frac{ \ta}{2}\rt\}, \rt.\notag\\
& \lt. \qquad \ \ \frac{\kappa }{2}+\frac{4\ga-3}{2\ga}-\frac{3}{2}\ta  \rt\} + \frac{2\ga-1}{\ga} -\ta ,\label{8.23.b1}\\
  b_2 =  &   \min\lt\{\frac{\kappa }{2}+\frac{4\ga-3}{2\ga}-\frac{3}{2}\ta, \ \frac{\kappa }{4}+\frac{10\ga-9}{4\ga}-\frac{9}{4}\ta\rt\} + \frac{\kappa }{2}+\frac{4\ga-3}{2\ga}-\frac{3}{2}\ta. \label{8.21.b2}
\end{align}
\end{lem}
{\em Proof}.  It follows from  \ef{hardyorigin} and \ef{tlg1} that
\begin{equation*}\label{}\begin{split}
   \int \bar\rho^{\ga} \mathcal{G}^2 dx    \le  &    C  \int_0^{1/2} \mathcal{G}^2 dx    +  4 \int_{1/2}^1 x^2 \bar\rho^{\ga} \mathcal{G}^2 dx    \le C \int_0^{1/2}  x^2 \lt( \mathcal{G}^2 + \mathcal{G}_x^2\rt) dx
  +  4 \int_{1/2}^1 x^2 \bar\rho^{\ga} \mathcal{G}^2 dx \\
\le & C \int    x^2  \bar\rho^\ga \mathcal{G}^2 dx  + C \int_0^{1/2} x^2 \bar\rho^{4\ga-2} \mathcal{G}_x^2  dx \\
\le  & C \int \bar\rho^\ga \lt(|xr_x-x|^2 + |r-x|^2\rt) dx  +C \int x^2 \bar\rho^{4\ga-2} \mathcal{G}_x^2  dx  ,
\end{split}\end{equation*}
which implies, using \ef{weightrx}, \ef{estlem51} and \ef{f3}, that
$$
(1+t)^{\frac{\kappa }{2}+\frac{4\ga-3}{2\ga}-\frac{3}{2}\ta}  \int \bar\rho^{\ga} \lt(|r_x-1|^2 + |r/x-1|^2\rt) (x,t)dx  \le C(\ta, \alpha)\widetilde{\mathfrak{F}}_\alpha(0).
$$
This gives \ef{Aug21.3}. It follows from \ef{7-2}, \ef{estlem51} and \ef{f3} that
$$
(1+t)^{\frac{\kappa }{2}+\frac{4\ga-3}{2\ga}-\frac{3}{2}\ta}\int  x^2\bar\rho^{2\ga-2}\mathcal{ G}_{tx}^2(x,t) dx
 \le C(\alpha, \ta) \widetilde{\mathfrak{F}}_\alpha(0),
$$
which, together with \ef{hh8.21} and \ef{estlem51}, gives \ef{Aug21.1}. \ef{Aug21.2} is   a conclusion of \ef{Aug21.3}, \ef{Aug21.1} and \ef{8-21-1}. \ef{Aug21.4} follows from \ef{f2}, \ef{Aug21.1} and \ef{estlem51}.

Next, we prove \ef{Aug21.5}. It follows from  \ef{f3} and \ef{bt1.8.19} that
\begin{equation}\label{8.23.1}\begin{split}
&(1+t)^{  \lt(\frac{\kappa }{2}+\frac{4\ga-3}{2\ga}-\frac{3}{2}\ta\rt)\frac{\aa+1}{2\ga-1+\aa}}\int  x^2 \bar\rho^{2\ga} \mathcal{G}_x^2 (x,t) dx
 \le C(\alpha, \ta) \widetilde{\mathfrak{F}}_\alpha(0) ,
\end{split}\end{equation}
due to
$$
\int x^2 \bar\rho^{2\ga} \mathcal{G}_x^2  dx\le
\lt(\int x^2 \bar\rho^{4\ga-2} \mathcal{G}_x^2  dx\rt)^{\frac{\aa+1}{2\ga-1+\aa}}\lt(\int x^2 \bar\rho^{2\ga-1-\alpha} \mathcal{G}_x^2  dx\rt)^{\frac{2\ga-2}{2\ga-1+\aa}}.
$$
In a similar way to deriving \ef{f3}, we have
$$
 (1+t)^{\frac{3 }{2}\kappa +\frac{2\ga-1}{2\ga}-\frac{ \ta}{2}}\int  x^2\bar\rho^{2\ga}\mathcal{ G}_{x}^2(x,t) dx
+\int_0^t (1+s)^{ \frac{3 }{2}\kappa +\frac{2\ga-1}{2\ga}-\frac{ \ta}{2}}\int x^2\bar\rho^{ \ga} \mathcal{G}_{xs}^2 dx ds
 \le C(\alpha, \ta) \widetilde{\mathfrak{F}}_\alpha(0),
$$
since \ef{ppqq} can be replaced by
\begin{equation*}\label{}\begin{split}
&\frac{1}{2}\frac{d}{dt}\int x^2 \mathfrak{P}^2 dx
+\frac{\mu \ga}{2}\int x^2\bar\rho^{\ga} \left(\frac{x^2 }{r^2 r_x } \right)^{\ga } \mathcal{G}_{xt}^2 dx \\
\le & C\lt( \lt\|x v_x\rt\|_{L^\iy}^2 +C\lt\|v\rt\|_{L^\iy}^2\rt) \int \bar{\rho}^{\ga}\mathcal{G}_x^2 dx + C \int \lt(x^2 v_x^2 +v^2 + v_t^2\rt) dx \\
\le & C(\ta) \mathfrak{E}(0) (1+t)^{-\kappa} \lt( \lt\|x v_x\rt\|_{L^\iy}^2 +\lt\|v\rt\|_{L^\iy}^2\rt)  + C \int \lt(x^2 v_x^2 +v^2 + v_t^2\rt) dx.
\end{split}\end{equation*}
This, together with \ef{8.23.1}, gives
\begin{equation}\label{8.23.2}\begin{split}
&(1+t)^{\max\lt\{  \lt(\frac{\kappa }{2}+\frac{4\ga-3}{2\ga}-\frac{3}{2}\ta\rt)\frac{\aa+1}{2\ga-1+\aa}, \  \frac{3 }{2}\kappa +\frac{2\ga-1}{2\ga}-\frac{ \ta}{2}\rt\}}\int  x^2 \bar\rho^{2\ga} \mathcal{G}_x^2 (x,t) dx
 \le C(\alpha, \ta) \widetilde{\mathfrak{F}}_\alpha(0) .
\end{split}\end{equation}
So, it yields from \ef{7-2}, \ef{8-21-2}, \ef{8.23.2} and \ef{estlem51}   that
\begin{equation}\label{8.22.3}\begin{split}
&(1+t)^{b_1-\frac{2\ga-1}{\ga} + \ta}\int  x^2 \mathcal{ G}_{tx}^2(x,t) dx
 \le C(\alpha, \ta) \widetilde{\mathfrak{F}}_\alpha(0) ,
\end{split}\end{equation}
where $b_1$ is defined by \ef{8.23.b1}. This, together with \ef{f1}, \ef{estlem51} and \ef{Aug21.1}, gives
$$
 (1+t)^{ b_1/2}
\|xv_x(\cdot,t)\|^2_{L^\iy}
 \le C(\alpha, \ta) \widetilde{\mathfrak{F}}_\alpha(0).
$$
Moreover, it follows from \ef{nj2}, \ef{nj2.new}, \ef{Aug21.1} and \ef{Aug21.2} that
\begin{equation}\label{8/21/2}\begin{split}
&(1+t)^{ b_2/ 2}
\|(v_x, v/x)(\cdot,t)\|^2_{L^\iy([0,1/4])}
 \le C(\alpha, \ta) \widetilde{\mathfrak{F}}_\alpha(0) ,
\end{split}\end{equation}
where $b_2$ is defined by \ef{8.21.b2}. Therefore, we have \ef{Aug21.5}.

In a similar way to deriving \ef{8/21/2}, one can get
\begin{equation}\label{8/21/3}\begin{split}
(1+t)^{b_2 /2  }
\|(r_x-1,r/x-1)(\cdot,t)\|^2_{L^\iy([0,1/4])}
 \le C(\alpha, \ta) \widetilde{\mathfrak{F}}_\alpha(0),
\end{split}\end{equation}
where $b_2$ is defined in \ef{8.21.b2}. It follows from   \ef{woman}, \ef{estlem51}, \ef{weightrxx} and \ef{8.20.1} that
\begin{align}
\lt\|x^{{3}/{2}}\bar\rho^{\ga/2}(r_x-1)(\cdot, t)\rt\|_{L^\iy}^2 \le C\lt\|x\bar\rho^{\ga /2}(r_x-1)(\cdot, t)\rt\|  \lt\|  \bar\rho^{\ga /2}r_{xx}(\cdot, t)\rt\|  \notag\\
 + C  \lt\|x\bar\rho^{\ga /2}(r_x-1)(\cdot, t)\rt\|^2
\le C(\alpha, \ta) \widetilde{\mathfrak{F}}_\alpha(0) (1+t)^{-\frac{\ka}{2} - \frac{2\ga-1}{2\ga} + \frac{\ta}{2}}  . \notag
\end{align}
Clearly, we can get the same estimate for $r/x-1$. Due to $ {\ka} + ({2\ga-1})/\ga - {\ta}  \le b_2 $, we then obtain \ef{Aug21.6}.

 \hfill $\Box$

\begin{rmk}\label{rmk8.22} As a consequence of Lemma \ref{lem44}, we have for any
$0<\theta<  2(\ga-1)/(3\ga)$,
\begin{align}
&(1+t)^{ \frac{11\ga-10}{4\ga} -\frac{5}{2}\ta } \lt\| \lt( v_{xx}, (v/x)_x \rt) (\cdot,t)\rt\|_{L^2([0,1/2])}^2 \le C(  \ta)\widetilde{\mathfrak{F}}_{\ga-\ga\ta}(0) , \label{8-23-1}\\
&(1+t)^{ \frac{9\ga-6}{4\ga} - \frac{3}{2}\ta} \lt\| \lt( xv_x, v \rt) (\cdot,t)\rt\|_{L^\iy}^2 \le C(  \ta)\widetilde{\mathfrak{F}}_{\ga-\ga\ta}(0),\label{8-22-1}
\\
& (1+t)^{  \frac{5\ga-4}{2\ga} - 2\ta   } \lt\| x \mathcal{G}_{tx} (\cdot,t)\rt\|^2 \le C(  \ta)\widetilde{\mathfrak{F}}_{\ga-\ga\ta}(0).\label{8-22-2}
\end{align}
Indeed, \ef{8-23-1} follows from \ef{Aug21.2},  \ef{8-22-1} from \ef{Aug21.4} and \ef{Aug21.5}, and  \ef{8-22-2} from   \ef{8.22.3}.
\end{rmk}

\subsubsection{Part III: further regularity}\label{sec3.4.4}
In this subsection, we further study the higher regularity of the strong solution obtained in Theorem \ref{mainthm1} and prove part $iii)$ and part $iv)$ of Theorem \ref{mainthm}.

\begin{lem}\label{lem10}  Let $\alpha\in [\ga, \ 2\ga-1]$  and $\mathfrak{F}_\alpha(0)<\infty$. For the global strong solution obtained in Theorem \ref{mainthm1}, there exist positive constants  $C$ and $C(\ta)$ independent of $t$ such that for any $0<\theta<  2(\ga-1)/(3\ga)$,
\begin{align}
& \lt\| \bar\rho^{(2\ga-1-\alpha)/2} \lt( r_{xx},  (r/x)_x \rt)(\cdot, t)\rt\|^2
\le   C \mathfrak{F}_\alpha (0)  +  C(\ta) (1+t)^{(\alpha-\ga + \ta\ga)/(\alpha-1)}\mathfrak{E}(0), \label{8.22.1} \\
& \lt\|\lt( r_{xx}, \ (r/x)_x \rt)(\cdot, t)\rt\|^2 \le C \mathfrak{F}_{2\ga-1} (0)  +  C (\ta) (1+t)^{\frac{1}{2}+\frac{\ga}{2\ga-2}\ta} \mathfrak{E}(0), \label{8/23/1}\\
&\lt\|\lt( v_{xx}, \ (v/x)_x \rt)(\cdot, t)\rt\|^2 \le C(\theta) (1+t)^{-\frac{7\ga-6}{4\ga} + 4 \ta }  \mathfrak{F}_{2\ga-1} (0) \lt(1+ \mathfrak{F}_{2\ga-1} (0)   \rt)  \lt(1+ \mathfrak{E} (0)  \rt)  , \label{growof2nd8.22}
\end{align}
provided that $\mathfrak{F}_{2\ga-1} (0)<\iy$ in \ef{8/23/1} and \ef{growof2nd8.22}.
\end{lem}
{\em Proof}. In a similar way to the derivation of \ef{bt1}, we have
\begin{equation*}\label{}\begin{split}
 &  \int \bar\rho^{2\ga-1-\alpha} \mathcal{G}_{x}^2(x,t)dx + \int_0^t  \int  \bar\rho^{3\ga-1-\alpha}\mathcal{G}_x^2 dx ds\\
  \le & C \mathfrak{F}_\alpha (0)   + C\int_0^t \int x^2 \bar\rho ^{\ga+1-\alpha} \left[ \left({r}/{x}-1\right)^2  + \left(r_x-1\right)^2 \right] dxds.
\end{split}\end{equation*}
It follows from the H$\ddot{o}$lder inequality that
\begin{equation*}\label{}\begin{split}
 & \int \bar\rho ^{\ga+1-\alpha} \left[ \left(r-x\right)^2  + \left(xr_x-x\right)^2 \right] dx\\
 & \le \lt( \int   \left[ \left(r-x\right)^2  + \left(xr_x-x\right)^2 \right] dx\rt)^{({\alpha-1})/{\ga}} \lt( \int \bar\rho ^{\ga} \left[ \left(r-x\right)^2  + \left(xr_x-x\right)^2 \right] dx\rt)^{({\ga+1-\alpha})/{\ga}}\\
 &\le (1+t)^{-({\ga-1})/{\ga}+\theta}\lt( (1+t)^{({\ga-1})/{\ga}-\theta} \int   \left[ \left(r-x\right)^2  + \left(xr_x-x\right)^2 \right] dx\rt)^{({\alpha-1})/{\ga}}\\
  &\times  \lt( (1+t)^{({\ga-1})/{\ga}-\theta} \int \bar\rho ^{\ga} \left[ \left(r-x\right)^2  + \left(xr_x-x\right)^2 \right] dx\rt)^{({\ga+1-\alpha})/{\ga}},
\end{split}\end{equation*}
which, together with \ef{estlem51} and  the Young inequality, implies that
\begin{equation*}\label{}\begin{split}
&\int_0^t \int x^2 \bar\rho ^{\ga+1-\alpha} \left[ \left({r}/{x}-1\right)^2  + \left(r_x-1\right)^2 \right] dxds\\
&\le C\int_0^t(1+s)^{({\ga-1})/{\ga}-\theta}\int \bar\rho ^{\ga} \left[ \left(r-x\right)^2  + \left(xr_x-x\right)^2 \right] dxds\\
&+C (1+t)^{1-\lt( {\ga-1} - \ga \theta\rt)/(\alpha-1) }\sup_{s\in [0,t]}\lt\{(1+s)^{({\ga-1})/{\ga}-\theta} \int    \left[ \left(r-x\right)^2  + \left(xr_x-x\right)^2 \right](x,s) dx\rt\} \\
&\le C(\ta) \mathfrak{E}(0)(1+t)^{(\alpha-\ga + \ta\ga)/(\alpha-1)  }.
\end{split}\end{equation*}
Thus,
\begin{equation}\label{8.22.2}\begin{split}
 &  \int \bar\rho^{2\ga-1-\alpha} \mathcal{G}_{x}^2(x,t)dx + \int_0^t  \int  \bar\rho^{3\ga-1-\alpha}\mathcal{G}_x^2 dx ds \le
   C \mathfrak{F}_\alpha (0)  +  C(\ta) \mathfrak{E}(0)(1+t)^{\frac{\alpha-\ga + \ta\ga}{\alpha-1}}.
\end{split}\end{equation}
This, together with   \ef{weightrxx} and \ef{gjrxx}, gives  \ef{8.22.1}.

Choose $\alpha=2\ga-1$ in \ef{8.22.2} to give
 \begin{equation}\label{8/23/2}\begin{split}
   \int  \mathcal{G}_{x}^2(x,t)dx   \le
   C \mathfrak{F}_{2\ga-1} (0)  +  C (\ta) (1+t)^{\frac{1}{2}+\frac{\ga}{2\ga-2}\ta} \mathfrak{E}(0).
\end{split}\end{equation}
This, together with  \ef{gjrxx}, gives \ef{8/23/1}.  It follows from \ef{gjvxxa}, \ef{8/23/2}, \ef{8-22-1} and \ef{8-22-2} that
 \begin{equation*}\label{}\begin{split}
&   \lt\|\lt( v_{xx}, \ (v/x)_x \rt)(\cdot, t)\rt\|_{L^2([1/2,1])}^2 \\
   \le  &   C(  \ta)    \widetilde{\mathfrak{F}}_{\ga-\ga\ta}(0)
   \lt\{(1+t)^{- \frac{ 5\ga-4  }{2\ga } + {2}\ta  }   +    (1+t)^{- \frac{9\ga-6}{4\ga}  + \frac{3}{2} \ta }
   \lt(  \mathfrak{F}_{2\ga-1} (0)   +  (1+t)^{\frac{1}{2}+\frac{\ga\ta}{2\ga-2}} \mathfrak{E}(0)\rt) \rt\} \\
  \le  & C(  \ta)    \widetilde{\mathfrak{F}}_{\ga-\ga\ta}(0) \lt(\mathfrak{F}_{2\ga-1} (0)  +1 \rt)(1+t)^{- \frac{7\ga-6}{4\ga} + 4 \ta }  .
\end{split}\end{equation*}
This, together with \ef{8-23-1}, implies \ef{growof2nd8.22}.

\hfill$\Box$

\begin{lem}\label{thelastlemmaforii}  Suppose that $\|x\bar\rho^{ {1}/{2}} v_{tt}(\cdot, 0)\|^2<\infty  $.  For the global strong solution obtained in Theorem \ref{mainthm1}, there exists a positive constant  $C$  independent of $t$ such that
\be\label{furthregularityofR}
\|x\bar\rho^{ {1}/{2}} v_{tt}(\cdot, t)\|^2+\int_0^{\infty}  \lt\|(v_{ss}, x v_{ssx})(\cdot,s)\rt\|^2 ds\le C \mathfrak{E}(0)+ C\|x\bar\rho^{ {1}/{2}} v_{tt}(\cdot, 0)\|^2  ,  \ t\ge 0.  \ee
\end{lem}
{\em Proof.}  Multiplying $\pl^2_t  \ef{nsp1} $ by $r^2v_{tt}$ and integrating the resulting equation both in $x$ and $t$, one can show  \ef{furthregularityofR}. Indeed, the derivation of \ef{furthregularityofR} is similar to that of \ef{eg2}, so we omit the details here.

\hfill $\Box$

\section{Proof of  Theorem \ref{mainthm2}}\label{sec4}
Due to Theorem \ref{mainthm} and Theorem \ref{mainthm1},  the triple $(\rho, u, R(t))$ ($t\ge 0$) defined by \ef{vacuumboundary} and \ef{solution} gives the unique global strong solution to the free boundary problem \ef{103}.
The decay estimates $i)$ and $ii)$ in Theorem \ref{mainthm2}    follow from the corresponding ones in Theorem \ref{mainthm}, by noting that
$$| \rho(r(x, t) ,t)-\bar\rho(x) | \le
C \bar\rho(x) \lt(|r_x(x,t)-1| + |x^{-1}r(x,t)-1|\rt)$$
and
$$u_r(r,t)=\frac{v_x(x,t)}{r_x(x,t)} \ \  {\rm and} \  \  \frac{u(r,t)}{r}= \frac{x}{r(x,t)} \frac{v(x,t)}{x}.$$

The $W^{2, \infty}$-estimate of $R(t)$ can be proved as follows. First,   it follows from \ef{lem4est} that
$$\int_0^\infty v_t^2(1, t)dt\le C\int_0^\infty\int_{\frac{1}{2}}^{1} (v_t^2+v_{xt}^2)dxdt\le C\mathfrak{E}(0). $$
One the other hand,  \ef{furthregularity} implies that
$$\int_0^\infty v_{tt}^2(1, t)dt\le C\int_0^\infty\int_{\frac{1}{2}}^{1} (v_{tt}^2+v_{xtt}^2)dxdt\le C \mathfrak{E}(0)+ C\|x\bar\rho^{ {1}/{2}} v_{tt}(\cdot, 0)\|^2. $$
Combining these two estimates with the fact that
$$\ddot{R}^2(t)=v_t^2(1, t)\le v_t^2(1, 0)+2\lt(\int_0^\infty v_t^2(1, t)dt\rt)^{1/2}\lt(\int_0^\infty v_{tt}^2(1, t)dt\rt)^{1/2}$$
gives
\ef{accelaration} immediately.
This finishes the proof of Theorem \ref{mainthm2}.

\hfill $\Box$

\begin{rmk} One may prove the boundedness of $r_{tt}(x, t)$ for any fixed $x\in [0, 1]$ and $t\ge 0$ if $|v_t(x, 0)|$ is finite by an argument  similar to the above.  This implies that every particle moving with the fluid has the bounded acceleration for $t\in (0, \infty)$ if it does so initially.
\end{rmk}

\centerline{Acknowledgement}
This research was partially supported by  the Zheng Ge Ru Foundation, and Hong
Kong RGC Earmarked Research Grants,  a Focus Area Grant from
The Chinese University of Hong Kong,  a grant from Croucher Foundation,  a NSF grant and a  NSFC grant. Zeng was also supported by the Center for Mathematical Sciences and Applications at Harvard University.

\renewcommand{\theequation}{A-\arabic{equation}}
\renewcommand{\thethm}{A-\arabic{thm}}
\setcounter{equation}{0}
\section*{Appendix}  

\subsection*{Part I. Local existence of strong solutions in the functional space $\mathfrak{E}(t)$}
In this part, we prove the local existence of strong solutions to  problem \ef{419} on a time interval $[0, T_*]$  for some $T_*>0$ in the function
space $\{(r, v):  \mathfrak{E}\in C([0, T_*])\}  $ by using a finite difference method, as used
in \cite{Okada, LiXY, LXY, Chengq}, where either a one-dimensional model or a three-dimensional model with spherical symmetry in a cut-off domain excluding a neighborhood of the origin is considered. Ideas used to derive the estimates in Section \ref{sec3} will be employed here to deal with the differences between the problem considered in this paper and those considered in \cite{Okada, LiXY, LXY, Chengq}.  The proof works for the case when $\|r_x(x, 0)-1\|_{L^{\infty}}$ is small. It may be possible to obtain a  local existence theory in our function space by only assuming that  $r_x(x, 0)$ has both lower and  upper positive bounds. However, we do not pursue this generality here which may need extra work  because the main purpose of this paper is to establish the global existence of  strong solutions for small data.


Recall that \ef{419a} reads
 \begin{equation}\label{nsp1-9.3}\begin{split}
& \bar\rho\left( \frac{x}{r}\right)^2  v_t - \lt( \mathfrak{B} + 4\la_1  \frac{v}{r}  \rt)_x =
  \lt\{  \bar\rho^{\ga} \lt[1-  \left(\frac{x^2}{r^2}\frac{1}{ r_x}\right)^\ga  \rt] \rt\}_x  + \left(\bar{\rho}^\ga\right)_x \lt( \frac{x^4}{r^4} -1 \rt)   ,
\end{split}
\end{equation}
where $\mathfrak{B}$ is given by \ef{bdry1}, that is,
$$\mathfrak{B}=\mu \frac{v_x}{r_x}+ \lt(2\la_2- \frac{4}{3} \la_1  \rt) \frac{v}{r}.$$

We use $(r^0, v^0)(x)$ to denote the initial data $(r, v)(x, 0)$(This notation will avoid possible confusions when we define the finite difference scheme.)
 The finite difference scheme is defined as follows.
Let $N$ be a positive integer, and $h={1}/{N}$.  For $n=0, 1, \cdots N$,  set
\be\label{a1}
x_n=nh,  \ \  \bar\rho_n=\bar\rho(x_n),  \ \ \bar q_n=(\bar\rho^{\ga})_x(x_n).
\ee
We can approximate \ef{nsp1-9.3} by the following initial value problem  for the system for $v_1,\cdots,v_{N-1}$ of ordinary differential equations: {
\begin{subequations}\label{a2}\begin{align}
& \bar\rho_n\lt(\frac{x_n}{r_n}\rt)^2\frac{dv_n}{dt} + \frac{1}{h} \lt[ \lt( \mathfrak{B}_{n} + 4\la_1 \frac{v_{n-1}}{r_{n-1}}  \rt)-\lt( \mathfrak{B}_{n+1} + 4\la_1 \frac{v_n}{r_n}  \rt) \rt]= \bar q_n \lt( \frac{x_n^4}{r_n^4} -1 \rt)
  \notag \\
& +  \frac{1}{h} \lt\{\bar\rho_{n+1}^\ga \lt[ 1 - \lt(\frac{h}{r_{n+1}-r_n}\rt)^\ga \lt(\frac{x_n}{r_n}\rt)^{2\ga}  \rt]  -\bar\rho_{n}^\ga \lt[ 1 - \lt(\frac{h}{r_{n }-r_{n-1}}\rt)^\ga \lt(\frac{x_{n-1}}{r_{n-1}}\rt)^{2\ga}  \rt] \rt\} , \label{a2-1}\\
& v_n(0)=v^0(x_n),     \label{a2-2}
\end{align}\end{subequations}}
where
\begin{align}
& r_n(t)=r^0(x_n)+\int_0^t v_n(s)ds ,\label{9.3-2} \\
& \mathfrak{B}_n=\mu\frac{v_n-v_{n-1}}{r_n-r_{n-1}}+ \lt(2\lambda_2-\frac{4}{3}\lambda_1\rt)\frac{v_{n-1}}{r_{n-1}} .     \label{9.3-1}
\end{align}
This system is supplemented by the following conditions to match the boundary conditions $v(0, t)=0$ and
$\mathfrak{B}(1, t)=0$:
\be\label{a3}
v_0(t) =0, \ \  \mathfrak{B}_N(t)=\mu\frac{v_N-v_{N-1}}{r_N-r_{N-1}}+ \lt(2\lambda_2-\frac{4}{3}\lambda_1\rt)\frac{v_{N-1}}{r_{N-1}}=0 .\ee
Clearly, if $r(0, 0)=0$, it follows from $v_0(t)=0$ that
\be\label{a3.1}
r_0(t)=0.
\ee

Next, we use the condition $\mathfrak{B}_N(t)=0$ to determine $v_N$ and $r_N$ in terms of $v_{N-1}$ and $r_{N-1}$. It follows from $(d/ dt) r_n = v_n$ that
$$0=\mathfrak{B}_N(t) = \frac{d}{dt}\lt\{\mu \ln \lt(\frac{r_N-r_{N-1}}{h}\rt) + \lt(2\lambda_2-\frac{4}{3}\lambda_1\rt) \ln \lt(\frac{r_{N-1}}{x_{N-1}}\rt) \rt\},$$
which implies
\be\label{a3.2}
r_N(t)=r_{N-1}(t) + \lt[r^0(x_N)- r^0(x_{N-1})\rt] \lt(\frac{r^0(x_{N-1})}{r_{N-1}(t)}\rt)^{\frac{2\lambda_2}{\mu}-\frac{4\lambda_1}{3\mu}}.
\ee
This, together with \ef{a3}, gives
\be\label{a3.3}
v_N(t)= v_{N-1}(t)-   \lt(\frac{2\lambda_2}{\mu}-\frac{4\lambda_1}{3\mu}\rt) \lt[r^0(x_N)- r^0(x_{N-1})\rt] \lt(\frac{r^0(x_{N-1})}{r_{N-1}(t)}\rt)^{\frac{2\lambda_2}{\mu}-\frac{4\lambda_1}{3\mu}} \frac{v_{N-1}(t)}{r_{N-1}(t)} .
\ee
With \ef{a3}-\ef{a3.3}, ODE system \ef{a2} is closed.

The approximation of the functional $\mathfrak{E}$ defined in \ef{mathmarch} is given by
\begin{align}
&\mathfrak{E}_N(t)   = \max_{1\le n\le N}\lt\{\lt|\frac{r_n(t)-r_{n-1}(t)}{h}-1\rt|^2+\lt|\frac{v_n(t)-v_{n-1}(t)}{h}\rt|^2\rt\}+h\sum_{n=1}^{N-1}\bar\rho_{n}|\dot{v}_n(t)|^2
\notag\\
&\quad+h\sum_{n=1}^{N-1}\bar\rho_n^{2\ga-1}\lt\{\lt|\frac{r_{n+1}(t)-2r_n(t)+r_{n-1}(t)}{h^2}\rt|^2
+\lt|\frac{1}{h}\lt(\frac{r_{n}(t)}{x_{n }}-\frac{r_{n-1}(t)}{x_{n-1}}\rt) \rt|^2\rt\}, \label{a6} \end{align}
where and  thereafter, $\dot g=(d/dt)g$ for a function $g=g(t)$.   It should be noted that $h\sum_{n=1}^{N-1}\bar\rho_{n}|\dot{v}_n(0)|^2$ can be given by the initial data via equation \ef{a2-1}. From now on, we choose a positive integer $N_0$ so large  that
\be\label{a4'} \mathfrak{E}_N(0)\le 2\mathfrak{E}(0), \ \  N\ge N_0.\ee

\begin{lem}\label{lem9.15} Let $\gamma\in (4/3, 2)$. Suppose that $\mathfrak{E}(0)< \iy$ and $r^0(0)=0$.
  Then there exist positive constants $N_1\ge N_0$ and  $ \bar \varepsilon>0$ independent of $N$ such that problem \ef{a2} admits a unique solution $(r_n, v_n)(t)$ on $[0, T^*]$ for some positive constant $T^*$  independent of $N$ satisfying
\be\label{a5}
\mathfrak{E}_N(t)\le K\mathfrak{E}_N(0) \le 2 K\mathfrak{E}(0) , \  \  t\in    [0, T^*], \  N\ge N_1 \ee
for some  positive constant  $K$ independent of $N$, provided that
\be\label{qqsn} \lt\|r^0_x-1\rt\|_{L^{\infty}(I)}\le  \bar \varepsilon. \ee
 Moreover, $T^*$ satisfies
 $$ T^*\ge \min\lt\{\frac{\bar c }{ \sqrt{ 2 K \mathfrak{E}(0) }}, \ \ \frac{1}{K(1+2 K \mathfrak{E}(0))} \rt\} $$
 for some positive constant $\bar c$ independent of $N$.
\end{lem}

\begin{rmk}  $\gamma\in (4/3, 2)$ is in general not necessary for the local existence,  $\gamma>1$ should be sufficient.
The reason we put this condition in the lemma is to ensure the existence and uniqueness of the stationary solution, the Lane-Emden solution, to keep the consistence with the global existence theory. \end{rmk}

\noindent{\em Proof of Lemma \ref{lem9.15}}. It follows from the ODE theory that problem \ef{a2} has a solution on a time interval. Let $T_N>0$ be the maximum existence time.  It  follows from $v_0(t)=0$ that
 \be\label{9.7-4}
 \max_{1\le n \le N} \lt|{v_n(t)}/{x_n} \rt|^2 \le   \mathfrak{E}_N(t), \ \ \ \   t\in [0, T_N);
\ee
which, together with $\dot{r}_n =v_n$, implies that for $n=1,\cdots, N$, and  $t\in [0, T_N)$,
$$
\lt|\frac{r_n(t)-r_{n-1}(t)}{h}-1\rt|\le    \lt|\frac{r_0(x_n)-r_{0}(x_{n-1})}{h}-1\rt| + \int_0^t \lt| \frac{v_n(s)-v_{n-1}(s)}{h}  \rt|  ds
\le   \lt\|r^0_x-1\rt\|_{L^{\infty}(I)}+ t \sup_{s\in [0,t]} \sqrt{ \mathfrak{E}_N(s)}  .
$$
This implies, due to $r_0(t)=0$, that for $n=1,\cdots, N$, and  $t\in [0, T_N)$,
$$
\lt|{r_n(t)}/{x_n}-1\rt|\le  \lt\|r^0_x-1\rt\|_{L^{\infty}(I)}\le   \lt\|r^0_x-1\rt\|_{L^{\infty}(I)}+ t \sup_{s\in [0,t]} \sqrt{ \mathfrak{E}_N(s)}   .
$$
Therefore, one can check that, 
for $n=1,\cdots, N$, and $t\in [0, T]$,
\begin{align}\label{9.7-1}
\lt|  \lt(r_n(t)-r_{n-1}(t)\rt)/h-1\rt|\le  2\bar\varepsilon \ \ {\rm and} \ \   \lt|{r_n(t)}/{x_n}-1\rt|\le  2\bar\varepsilon.
\end{align}
 provided 
\be\label{9.29.3}
T \sup_{s\in [0,T]} \sqrt{ \mathfrak{E}_N(s)} \le  \bar\varepsilon,
\ee
and 
\be\label{initialforrox}\lt\|r^0_x-1\rt\|_{L^{\infty}(I)}\le \bar\varepsilon\ee
for some constant $\bar\varepsilon>0$.
In particular, it holds that  for $n=1,\cdots, N$, and $t\in [0, T]$,
\be\label{9.6-1}
\frac{1}{2} \le    \frac{r_n(t)}{x_n}   \le \frac{3}{2} \ \ {\rm and} \ \     \frac{1}{2} \le  \frac{r_n(t)-r_{n-1}(t)}{h}  \le \frac{3}{2} ,
\ee
if \be \bar\varepsilon\le \frac{1}{4}. \ee
With \ef{9.7-4}, \ef{9.7-1} and \ef{9.6-1}, we will prove in {\em Steps 1-3} that for sufficiently large $N$, there exists constant $\bar C>1$ independent of $N$ such that if $T$ satisfies \ef{9.29.3}, then
\begin{align} \label{9-29-2}
  \mathfrak{E}_N(t)
  \le    \bar C  \mathfrak{E}_N(0)  + \bar C \int_0^t \lt[ \mathfrak{E}_N(s) + \mathfrak{E}_N^2 (s) \rt] ds, \ \ \ \  t\in [0, T].
\end{align}
Once this statement is proved, the lemma will follow from an argument which we give in {\em Step 4}.

{\em Step 1}. In this step, we prove that
\begin{align}
&h \sum_{n=1}^{N-1}\bar\rho_n x_{n-1}^2 \dot{v}_n^2(t) + \int_0^t h \sum_{n=1}^{N-1} \lt[ x_{n}^2  \lt(\frac{\dot{v}_{n+1}(s)- \dot{v}_{n}(s)}{h} \rt)^2  + \dot{v}_n^2(s) \rt]ds \notag\\
\le & C  \mathfrak{E}_N(0)
+ C \int_0^t \lt[ \mathfrak{E}_N(s) + \mathfrak{E}_N^2 (s) \rt] ds.
 \label{9-7-1}
\end{align}
 Here and in the rest of this part of Appendix,   $C$  denotes a generic positive constant independent of $t$ and $N$.
It follows from \ef{a2-1}, $\dot{r}_n=v_n$ and \ef{9.6-1}  that
\begin{align}
  \bar\rho_n\lt(\frac{x_n}{r_n}\rt)^2 \ddot{v}_n + \frac{1}{h} \lt[\lt( \dot{\mathfrak{B}}_{n} -\dot{\mathfrak{B}}_{n+1} \rt) + 4\la_1\lt( \frac{\dot{v}_{n-1}}{r_{n-1}}  -   \frac{\dot{v}_n}{r_n}  \rt) \rt]
 = \frac{1}{h} \lt( \mathcal{P}_{n+1} - \mathcal{P}_n \rt) + e_n
  , \label{9.7-5}
\end{align}
where
\begin{align*}
&\mathcal{P}_n = \ga\bar\rho_n  \lt(\frac{h}{r_{n }-r_{n-1}}\rt)^\ga \lt(\frac{x_{n-1}}{r_{n-1}}\rt)^{2\ga}  \lt( \frac{v_{n }-v_{n-1}}{r_{n }-r_{n-1}}  + 2\frac {v_{n-1}}{r_{n-1}}\rt),\\
& |e_n| \le C \lt(  \bar\rho_n \lt|\frac{v_n}{x_n}\rt| |\dot{v}_n| + \frac{1}{h} \lt|\frac{v_n }{r_n }- \frac{v_{n-1} }{r_{n-1} }\rt| \lt|\frac{v_n }{r_n } + \frac{v_{n-1} }{r_{n-1} }\rt|  + \lt| \frac{v_n}{x_n} \rt| \rt).
\end{align*}
\ef{9-7-1} follows from the summation of the product of $r_{n-1}^2 \dot{v}_n$ and \ef{9.7-5}. First, we  analyze the second term on the left-hand side of \ef{9.7-5}.  Notice that
\begin{align*}
& \sum_{n=1}^{N-1} \frac{1}{h}\lt( \dot{\mathfrak{B}}_{n} -\dot{\mathfrak{B}}_{n+1} \rt) r_{n-1}^2 \dot{v}_n =
 \sum_{n=1}^{N-1} \frac{1}{h}\dot{\mathfrak{B}}_{n+1}\lt( r_{n}^2 \dot{v}_{n+1}  - r_{n-1}^2 \dot{v}_{n}   \rt) \\
& =  \sum_{n=1}^{N-1} \dot{\mathfrak{B}}_{n+1}\lt( r_{n}^2 \frac{\dot{v}_{n+1} -\dot{v}_{n}}{h}  + \frac{r_n^2 - r_{n-1}^2 }{h} \dot{v}_{n}   \rt)
\ge  \sum_{n=1}^{N-1} \dot{\mathfrak{B}}_{n+1} \lt( x_n^2 \frac{\dot{v}_{n+1} -\dot{v}_{n}}{h}  +2 x_n \dot{v}_{n}  \rt) \\
&\qquad -  C \sum_{n=1}^{N-1} \lt| \dot{\mathfrak{B}}_{n+1} \rt| \lt(\bar\varepsilon x_n^2 \lt| \frac{\dot{v}_{n+1} -\dot{v}_{n}}{h}\rt|  +  \lt(\bar\varepsilon+  n^{-1} \rt)  x_n  \lt| \dot{v}_{n}\rt|  \rt)
\end{align*}
and
\begin{align*}
&\lt| \dot{\mathfrak{B}}_{n+1} -\lt(\mu \frac{\dot{v}_{n+1}-\dot{v}_n}{h} + \lt(2\la_2-\frac{4}{3}\la_1\rt) \frac{\dot{v}_n}{x_n}  \rt)  \rt|\\
& \le
C \bar\varepsilon \lt( \lt|\frac{\dot{v}_{n+1}-\dot{v}_n}{h}\rt| + \lt|\frac{\dot{v}_n}{x_n}\rt| \rt) + C  \lt( \lt|\frac{ {v}_{n+1}- {v}_n}{h}\rt|^2 + \lt|\frac{ {v}_n}{x_n}\rt|^2 \rt).
\end{align*}
It then yeilds from the Cauchy inequality that,  for any $\da\in (0,1)$,
\begin{align}
& \sum_{n=1}^{N-1} \frac{1}{h}\lt( \dot{\mathfrak{B}}_{n} -\dot{\mathfrak{B}}_{n+1} \rt) r_{n-1}^2 \dot{v}_n  \ge  - C \da^{-1}  \mathfrak{E}_N \sum_{n=1}^{N-1} \lt\{  x_n^2 \lt( \frac{ {v}_{n+1} - {v}_{n}}{h} \rt)^2 + {v}_{n}^2  \rt\}\notag\\
& \quad + \sum_{n=1}^{N-1} \lt\{\mu x_n^2 \lt( \frac{\dot{v}_{n+1} -\dot{v}_{n}}{h} \rt)^2 + \lt(4\la_2- \frac{8}{3}\la_1\rt) \dot{v}_{n}^2 +   \lt(    4\la_2 + \frac{4}{3}\la_1 \rt)x_n \frac{\dot{v}_{n+1} -\dot{v}_{n}}{h}\dot{v}_{n} \rt\} \notag\\
 &\quad - C     \sum_{n=1}^{N-1} \lt(n^{-1} + \da^{-1} n^{-2}\rt) \dot{v}_{n}^2   - C ( \bar\varepsilon +  2\da)\sum_{n=1}^{N-1} \lt\{  x_n^2 \lt( \frac{\dot{v}_{n+1} -\dot{v}_{n}}{h} \rt)^2 + \dot{v}_{n}^2  \rt\}. \label{9.29.1}
\end{align}
Similarly,
\begin{align}
&\sum_{n=1}^{N-1} \frac{1}{h}\lt( \frac{\dot{v}_{n-1}}{r_{n-1}}  -   \frac{\dot{v}_n}{r_n}  \rt) r_{n-1}^2 \dot{v}_n
= \sum_{n=1}^{N-1} r_{n-1} \frac{\dot{v}_{n-1} -\dot{v}_{n}}{h} \dot{v}_n
+ \sum_{n=1}^{N-1}  \frac{r_{n-1}}{r_n} \frac{r_{n} -r_{n-1}}{h} \dot{v}_n^2 \notag\\
\ge & -\sum_{n=0}^{N-2} \lt\{ r_n   \frac{\dot{v}_{n+1} -\dot{v}_{n}}{h}\dot{v}_{n} + h   r_n \lt(\frac{\dot{v}_{n+1} -\dot{v}_{n}}{h}\rt)^2 \rt\} + \sum_{n=1}^{N-1}  \lt\{  \dot{v}_n^2 - C      \lt( n^{-1} +   \bar\varepsilon\rt)\dot{v}_n^2 \rt\}. \label{9.29.2}
\end{align}
So, if $\bar\varepsilon$ is small (the smallness is independent of $N$) and $N$ is large (which implies $h$ is small),  it follows from \ef{9.29.1} and \ef{9.29.2} that
\begin{align}
& \sum_{n=1}^{N-1} \frac{1}{h} \lt[\lt( \dot{\mathfrak{B}}_{n} -\dot{\mathfrak{B}}_{n+1} \rt) + 4\la_1\lt( \frac{\dot{v}_{n-1}}{r_{n-1}}  -   \frac{\dot{v}_n}{r_n}  \rt) \rt]r_{n-1}^2 \dot{v}_n  \ge  - C \sum_{n=1}^{N-1}  \frac{1}{n} \dot{v}_{n}^2 \notag\\
 & \quad + \frac{1}{2} \sigma \sum_{n=1}^{N-1}\lt[ x_{n}^2  \lt(\frac{\dot{v}_{n+1} - \dot{v}_{n} }{h} \rt)^2  + \dot{v}_n^2 \rt]
 -  C   \mathfrak{E}_N \sum_{n=1}^{N-1} \lt\{  x_n^2 \lt( \frac{ {v}_{n+1} - {v}_{n}}{h} \rt)^2 + {v}_{n}^2  \rt\} . \label{9.7-6}
\end{align}
 Here $\sigma=\min\{2\la_1/3, \ \la_2\}$. (Indeed, the first term on the second line of \ef{9.7-6} is obtained by considering the following two cases: $n=1,\cdots, N-2$, and $n=N-1$, separately.)  The first term on the right-hand side of \ef{9.7-6} can be bounded by
\begin{align*}
 \sum_{n=1}^{N-1}  \frac{1}{n} \dot{v}_{n}^2 \le \sum_{n=1}^{[N/2]}    \dot{v}_{n}^2 + \frac{2}{N} \sum_{n=[N/2]+1}^{N-1}  \dot{v}_{n}^2\le C \sum_{n=1}^{[N/2]} \bar\rho_n   \dot{v}_{n}^2 + \frac{2}{N} \sum_{n=[N/2]+1}^{N-1}  \dot{v}_{n}^2 \le \sum_{n=1}^{N-1}  \lt(  C\bar\rho_n +  \frac{2}{N} \rt) \dot{v}_{n}^2 .
\end{align*}
This, together with \ef{9.7-6}, gives that for sufficiently large $N$,
\begin{align}
& \sum_{n=1}^{N-1} \frac{1}{h} \lt[\lt( \dot{\mathfrak{B}}_{n} -\dot{\mathfrak{B}}_{n+1} \rt) + 4\la_1\lt( \frac{\dot{v}_{n-1}}{r_{n-1}}  -   \frac{\dot{v}_n}{r_n}  \rt) \rt]r_{n-1}^2 \dot{v}_n  \ge  - C \sum_{n=1}^{N-1}  \bar\rho_n \dot{v}_{n}^2 \notag\\
 & \quad + \frac{1}{4} \sigma \sum_{n=1}^{N-1}\lt[ x_{n}^2  \lt(\frac{\dot{v}_{n+1} - \dot{v}_{n} }{h} \rt)^2  + \dot{v}_n^2 \rt]
 -  C  \mathfrak{E}_N \sum_{n=1}^{N-1} \lt\{  x_n^2 \lt( \frac{ {v}_{n+1} - {v}_{n}}{h} \rt)^2 + {v}_{n}^2  \rt\} .\label{9.7-7}
\end{align}
In a similar but much easier way, we can deal with the other terms in \ef{9.7-5} and  obtain \ef{9-7-1}.

{\em Step 2}. In this step, we prove that
\begin{align}
& h \sum_{n=1}^{N-1}\bar\rho_n   \dot{v}_n^2(t) + \int_0^t h  \lt[\sum_{n=1}^{N} \lt(\frac{\dot{v}_{n }(s)- \dot{v}_{n-1}(s)}{h} \rt)^2  + \sum_{n=1}^{N} \frac{\dot{v}_n^2 (s)}{ x_{n}^2   }  \rt]ds \notag\\
\le &   C  \mathfrak{E}_N(0)  + C \int_0^t \lt[ \mathfrak{E}_N(s) + \mathfrak{E}_N^2 (s) \rt] ds.  \label{9-7-2}
\end{align}
Rewrite \ef{9.7-5} as
\begin{align}
 \bar\rho_n\lt(\frac{x_n}{r_n}\rt)^2 \ddot{v}_n + \frac{1}{h} \mu \lt(\mathcal{Q}_n - \mathcal{Q}_{n+1}  \rt)
=  \frac{1}{h} \lt( \mathcal{P}_{n+1} - \mathcal{P}_n \rt) + \bar e_n , \label{9-7-3}
\end{align}
where
\begin{align*}
&\mathcal{Q}_n= \frac{\dot{v}_n-\dot{v}_{n-1}}{r_n-r_{n-1}}+ 2 \frac{\dot{v}_{n-1}}{r_{n-1}} -\frac{\lt( {v}_n- {v}_{n-1}\rt)^2}{\lt(r_n-r_{n-1}\rt)^2}- 2 \frac{ {v}_{n-1}^2}{r_{n-1}^2} ,\\
& |\bar e_n| \le  C \lt(  \bar\rho_n |\dot{v}_n| \lt| {v_n}/{x_n} \rt| + \lt| {v_n}/{x_n} \rt|  \rt) .
\end{align*}
Let $\psi$ be a non-increasing function defined on $[0,1]$ satisfying
\begin{equation*}\label{}\begin{split}
 \psi=1 \ \ {\rm on } \ \ [0,1/4], \ \ \psi=0 \ \ {\rm on } \ \  [1/2, 1] \ \ {\rm and} \  \ |\psi'|\le 32.
\end{split}
\end{equation*}
Set $\psi_n=\psi(x_n)$. Similarly, we consider the summation of the product of \ef{9-7-3} and $\psi_n \dot{v}_n$. To deal with the second term on the left-hand side of \ef{9-7-3},  one  notices that
\begin{align*}
\sum_{n=1}^{N-1}\frac{1}{h} \lt(\mathcal{Q}_n - \mathcal{Q}_{n+1}  \rt)\psi_n \dot{v}_n
= &\sum_{n=1}^{N-1}\frac{1}{h} \mathcal{Q}_n \lt(\psi_n \dot{v}_n - \psi_{n-1} \dot{v}_{n-1}  \rt) \\
=&  \sum_{n=1}^{N-1}  \mathcal{Q}_n \lt(\psi_n \frac{\dot{v}_n -   \dot{v}_{n-1} }{h} + \frac{\psi_n-\psi_{n-1}}{h} \dot{v}_{n-1} \rt);
\end{align*}
and
\begin{align*}
&\sum_{n=1}^{N-1}  \frac{\dot{v}_{n-1}}{r_{n-1}}  \psi_n \frac{\dot{v}_n -   \dot{v}_{n-1} }{h}
= \sum_{n=1}^{N-1} \frac{1}{h} \lt(  \psi_n \frac{\dot{v}_n \dot{v}_{n-1} }{r_{n-1}} -\psi_{n+1} \frac{\dot{v}_n ^2}{r_n}\rt)
 =    -\sum_{n=1}^{N-1}  \frac{\dot{v}_{n-1}}{r_{n-1}}  \psi_n \frac{\dot{v}_n -   \dot{v}_{n-1} }{h} \\
 &\quad - h \sum_{n=1}^{N-1} \frac{\psi_n}{r_{n-1}} \lt( \frac{ \dot{v}_{n-1} -\dot{v}_n }{ h} \rt)^2
+ \sum_{n=1}^{N-1}  \frac{\psi_n-\psi_{n+1}}{h r_n} \dot{v}_n^2 + \sum_{n=1}^{N-1} \frac{ \psi_n }{r_{n-1}} \frac{r_n-r_{n-1}}{h r_n  }   \dot{v}_n^2 \\
 & =  -\sum_{n=1}^{N-1}  \frac{\dot{v}_{n-1}}{r_{n-1}}  \psi_n \frac{\dot{v}_n -   \dot{v}_{n-1} }{h}  - h \sum_{n=2}^{N-1} \frac{\psi_n}{r_{n-1}} \lt( \frac{ \dot{v}_{n-1} -\dot{v}_n }{ h} \rt)^2 + \sum_{n=1}^{N-1}  \frac{\psi_n-\psi_{n+1}}{h r_n} \dot{v}_n^2  \\
  &\quad + \sum_{n=2}^{N-1} \frac{ \psi_n }{r_{n-1}} \frac{r_n-r_{n-1}}{h r_n  }   \dot{v}_n^2,
\end{align*}
which implies
\begin{align}
2 \sum_{n=1}^{N-1}  \frac{\dot{v}_{n-1}}{r_{n-1}}  \psi_n \frac{\dot{v}_n -   \dot{v}_{n-1} }{h}
\ge  &    \sum_{n=2}^{N-1}   \psi_n  \frac{r_n-r_{n-1}}{h }  \frac{  \dot{v}_n^2 }{r_n^2 }  -   h \sum_{n=2}^{N-1} \frac{\psi_n}{r_{n-1}} \lt( \frac{ \dot{v}_{n-1} -\dot{v}_n }{ h} \rt)^2 \notag\\
&+   \sum_{n=1}^{N-1}  \frac{\psi_n-\psi_{n+1}}{h r_n} \dot{v}_n^2 . \notag
\end{align}
Then,  \ef{9-7-2} follows from \ef{9-7-1}, \ef{a3}, \ef{a3.2}, \ef{a3.3} and simple calculations.

{\em Step 3}. In this step, we prove
\begin{align} \label{9-29-1}
& \mathfrak{E}_N(t)  + h\sum_{n=1}^{N-1}\bar\rho_n^{2\ga-1}\lt\{\lt|\frac{v_{n+1}(t)-2v_n(t)+v_{n-1}(t)}{h^2}\rt|^2 +\lt|\frac{1}{h}\lt(\frac{v_{n }(t)}{x_{n }}-\frac{v_{n-1}(t)}{x_{n-1}}\rt) \rt|^2\rt\}\notag\\
  \le  &  C  \mathfrak{E}_N(0)  +   C \int_0^t \lt[ \mathfrak{E}_N(s) + \mathfrak{E}_N^2 (s) \rt] ds, \ \ \ \  t\in [0, T].
\end{align}
To this end, we rewrite \ef{a2-1} as
\begin{align}
\mu \frac{\dot{\mathcal{G}}_n-\dot{\mathcal{G}}_{n+1}}{h}
= &
-\bar\rho_n\lt(\frac{x_n}{r_n}\rt)^2\dot{v}_n +  \bar\rho_{n}^\ga \frac{\exp\{-\ga \mathcal{G}_{n+1}\} - \exp\{-\ga \mathcal{G}_{n}\}}{h} \notag\\
& + \frac{\bar\rho_{n+1}^\ga - \bar\rho_{n}^\ga}{h}  \lt[ 1 - \lt(\frac{h}{r_{n+1}-r_n}\rt)^\ga \lt(\frac{x_n}{r_n}\rt)^{2\ga}  \rt]
 + \bar q_n \lt( \frac{x_n^4}{r_n^4} -1 \rt) = :\ell_n.  \label{9.9.1}
\end{align}
Here for $n=1,\cdots,N-1$,
 \begin{align}
& \mathcal{G}_n= \ln\lt(\frac{r_n-r_{n-1}}{h}\rt) + 2\ln\lt(\frac{r_{n-1}}{x_{n-1}}\rt),  \label{9.16-1}\\
&  |\ell_n |
\le   C  \bar\rho_n \lt\{ |\dot{v}_n| +  \bar\rho_{n}^{\ga-1} \lt|\frac{\mathcal{G}_{n}-\mathcal{G}_{n+1}}{h}\rt|  +\lt|\frac{r_{n+1}-r_{n}}{h}-1\rt|+\lt|\frac{r_{n}}{x_n}-1\rt|\rt\} . \notag
\end{align}
It is easy to derive from \ef{9.9.1} and \ef{9-7-2} that
\begin{align}
 h \sum_{n=1}^{N-1} \frac{1}{\bar\rho_n} \lt(\frac{\dot{\mathcal{G}}_n(t)-\dot{\mathcal{G}}_{n+1}(t)}{h}\rt)^2
\le C  \mathfrak{E}_N(0)  + C \int_0^t \lt[ \mathfrak{E}_N(s) + \mathfrak{E}_N^2 (s) \rt] ds  , \label{9.10.1}\end{align}
and then
\begin{align}
 h \sum_{n=1}^{N-1}  \bar\rho_n^{2\ga-1} \lt(\frac{ {\mathcal{G}}_n(t)- {\mathcal{G}}_{n+1}(t)}{h}\rt)^2
\le     C  \mathfrak{E}_N(0)  + C \int_0^t \lt[ \mathfrak{E}_N(s) + \mathfrak{E}_N^2 (s) \rt] ds . \label{9.10.2}
\end{align}
Following the arguments in Section \ref{sec3.4.1} (in particular, Lemmas \ref{lemhh1} and \ref{boundsforrv}), we can use  \ef{9-7-2},  \ef{9.10.1} and \ef{9.10.2} to get \ef{9-29-1}.

{\em Step 4}.
Set $T^*=\sup\{t:  \mathfrak{E}_N(t)\le 2\bar C\mathfrak{E}_N(0)\}$,
where $\bar C>1$ is the constant in \ef{9.15-1}.
Obviously $T^*>0$ since $\bar C>1$.
We may assume that $T^*<  {\bar \varepsilon}/ { \sqrt {2\bar C\mathfrak{E}_N(0)}}$.
(Otherwise,  $T^* \ge {\bar \varepsilon}/ { \sqrt {4\bar C\mathfrak{E} (0)}}$
due to $ \mathfrak{E}_N(0) \le 2  \mathfrak{E}(0)$, and the lower bound is independent of $N$, so that   the lemma is proved.)  By the definition of $T^*$,  
$\mathfrak{E}_N(t)\le 2\bar C\mathfrak{E}_N(0)$ for $t\in [0, T^*]$. Then,
$$ T^*\sup_{s\in [0,T^*]} \sqrt{ \mathfrak{E}_N(s)} \le T^* \sqrt{2\bar C\mathfrak{E}_N(0)} \le  \bar\varepsilon.$$
So, $T^*$ satisfies \ef{9.29.3}, and it follows from estimate \ef{9-29-2} and definition of $T^*$ that
\begin{align*}
2\bar C\mathfrak{E}_N(0)=  \mathfrak{E}_N(T^*)
  \le    \bar C  \mathfrak{E}_N(0)  + \bar C \int_0^{T^*} \lt[ \mathfrak{E}_N(s) + \mathfrak{E}_N^2 (s) \rt] ds  \\
  \le  \bar C\mathfrak{E}_N(0) +  \bar C  {T^*} 2 \bar C\mathfrak{E}_N(0) \lt(1 +2 \bar C\mathfrak{E}_N(0)  \rt).
\end{align*}
This implies
$$T^*\ge \frac{1}{2 \bar C  \lt(1 +2 \bar C\mathfrak{E}_N(0)  \rt) }\ge \frac{1}{2 \bar C  \lt(1 + 4 \bar C\mathfrak{E}(0)  \rt) }. $$
The Lemma is proved by taking  $K=2\bar C$.

\hfill
$\Box$

\begin{rmk} Let $ T_*=\min\{T^*,  1\}$.  For the solution $(r_n, v_n)$ obtained in Lemma \ref{lem9.15}, it holds that
\begin{align} \label{9.15-1}
& \mathfrak{E}_N(t)  + h\sum_{n=1}^{N-1}\bar\rho_n^{2\ga-1}\lt\{\lt|\frac{v_{n+1}(t)-2v_n(t)+v_{n-1}(t)}{h^2}\rt|^2 +\lt|\frac{1}{h}\lt(\frac{v_{n }(t)}{x_{n }}-\frac{v_{n-1}(t)}{x_{n-1}}\rt) \rt|^2\rt\}\notag\\
&  \le     C  \mathfrak{E}(0) +   C \mathfrak{E}^2(0)   , \ \ t\in [0, T_*]  .
\end{align}
  Indeed, \ef{9.15-1} follows from \ef{9-29-1} and \ef{a5}.
\end{rmk}

\begin{lem}\label{lem9.16} Suppose that the  assumptions  in Lemma \ref{lem9.15} are satisfied. Then the following estimates are obtained for any $n=1,\cdots, N$, and $t,s \in [0, T_*]$,
\begin{align}
&\frac{1}{2} \le    \frac{r_n(t)}{x_n} \le \frac{3}{2}  , \   \frac{1}{2} \le  \frac{(r_n -r_{n-1})(t)}{h}  \le \frac{3}{2} ,   \   \lt|\frac{v_n(t)}{x_n} \rt| + \lt|\frac{(v_n-v_{n-1})(t)}{h} \rt| \le \sqrt{ C\mathfrak{E}(0) },    \label{9.15.2} \\
& h\sum_{n=1}^N \lt|\frac{r_n(t)-r_{n-1}(t)}{h}\rt|^2+  h\sum_{n=1}^N\lt|\frac{v_{n}(t)-v_{n-1}(t)}{h}\rt|^2
 \le  C\lt( \mathfrak{E}(0)       +1 \rt)  ,\label{9.16.1}\\
& h\sum_{n=1}^{N }   \lt|\frac{1}{h}\lt(\frac{r_{n }(t)}{x_{n }}-\frac{r_{n-1}(t)}{x_{n-1}}\rt) \rt|^2
+ h\sum_{n=1}^{N }\lt|\frac{1}{h}\lt(\frac{v_{n }(t)}{x_{n }}-\frac{v_{n-1}(t)}{x_{n-1}}\rt) \rt|^2 \le C  \sum_{i=0}^2\mathfrak{E}^i(0)   ,\label{9.16.2}\\
&h\sum_{n=1}^{N-1}   \lt|\frac{1}{h}\lt( \bar\rho_{n+1}^{\ga-\frac{1}{2}}\frac{r_{n+1}(t)-r_n(t)}{h} - \bar\rho_{n}^{\ga-\frac{1}{2}}\frac{r_{n }(t)-r_{n-1}(t)}{h}\rt) \rt|^2  \le C(  \mathfrak{E}(0) +1 ), \label{9.16.3} \\
&h\sum_{n=1}^{N-1}   \lt|\frac{1}{h}\lt( \frac{v_{n+1}(t)-v_n(t)}{r_{n+1}(t)-r_n(t)} - \frac{v_{n }(t)-v_{n-1}(t)}{r_{n}(t)-r_{n-1}(t)}\rt) \rt|^2  \le C  \sum_{i=0}^3\mathfrak{E}^i(0)  , \label{9.16.4} \\
& h \sum_{n=1}^N \lt|v_n(t)-v_n(s)\rt|^2 +   h\sum_{n=1}^{N}\lt| \frac{v_{n }(t)}{x_{n }}-\frac{v_{n }(s)}{x_{n }}  \rt|^2  \le C\mathfrak{E}(0)|t-s|,  \label{a9} \\
& h\sum_{n=1}^{N }   \lt|  \frac{r_{n }(t)-r_{n-1}(t)}{h} -\frac{r_{n }(s)-r_{n-1}(s)}{h}\rt|^2  \le C\mathfrak{E}(0)|t-s|^2  ,\label{9.16.5}\\
& h\sum_{n=1}^{N }   \lt| \frac{v_{n }(t)-v_{n-1}(t)}{r_{n}(t)-r_{n-1}(t)} -\frac{v_{n }(s)-v_{n-1}(s)}{r_{n}(s)-r_{n-1}(s)}  \rt|^2  \le C \sum_{i=1}^2 \mathfrak{E}^i(0) |t-s| + C \mathfrak{E}^2(0)|t-s|^2   .  \label{9.16.6}
\end{align}
Here $C$ is a constant independent of $N$.
\end{lem}
{\em Proof}.  Clearly, \ef{9.15.2} follows from \ef{a5}, \ef{9.7-4} and \ef{9.6-1}; and \ef{9.16.1} follows
  from \ef{a5}. For \ef{9.16.2}, it follows from \ef{9.15.2} and \ef{9.15-1} that
\begin{align*}
&h\sum_{n=1}^{N }\lt|\frac{1}{h}\lt(\frac{v_{n }(t)}{x_{n }}-\frac{v_{n-1}(t)}{x_{n-1}}\rt) \rt|^2
=   h \lt( \sum_{n=1}^{[N/2] } +  \sum_{n=[N/2]+1}^{N }\rt)\lt|\frac{1}{h}\lt(\frac{v_{n }(t)}{x_{n }}-\frac{v_{n-1}(t)}{x_{n-1}}\rt) \rt|^2 \\
\le  &  C h \sum_{n=1}^{[N/2] } \bar\rho_n^{2\ga-1}\lt|\frac{1}{h}\lt(\frac{v_{n }(t)}{x_{n }}-\frac{v_{n-1}(t)}{x_{n-1}}\rt) \rt|^2
  + C h \sum_{n=[N/2]+1}^{N }\lt( \lt|\frac{v_n(t)-v_{n-1}(t) }{h}   \rt|^2 + \lt|\frac{v_n(t)}{x_n}\rt|^2\rt)\\
\le & C\mathfrak{E}(0) +  C \mathfrak{E}^2(0).
\end{align*}
Similarly, the estimates for $r_n$ in \ef{9.16.2} follows from \ef{9.15.2} and \ef{a5}. \ef{9.16.3} follows from simple calculations and \ef{a5}. For \ef{9.16.4}, we note that
\begin{align*}
\frac{1}{h}\lt| \frac{v_{n+1}(t)-v_n(t)}{r_{n+1}(t)-r_n(t)} - \frac{v_{n }(t)-v_{n-1}(t)}{r_{n}(t)-r_{n-1}(t)} \rt| =\frac{1}{h} \lt|\lt(\dot{\mathcal{G}}_{n+1} -  2\frac{v_{n}}{r_{n}}\rt) -\lt(\dot{\mathcal{G}}_{n } -  2\frac{v_{n-1}}{r_{n-1}}\rt)\rt|\\
\le \lt|\frac{\dot{\mathcal{G}}_{n+1} -\dot{\mathcal{G}}_{n} }{h}\rt|+ \frac{2}{h} \lt(\frac{x_n}{r_n}\lt|\frac{v_n}{x_n}-\frac{v_{n-1}}{x_{n-1}}\rt| + \lt|\frac{v_{n-1}}{x_{n-1}} \rt| \frac{x_n}{r_n}\frac{x_{n-1}}{r_{n-1}} \lt|\frac{r_n}{x_n}-\frac{r_{n-1}}{x_{n-1}}\rt| \rt).
\end{align*}
Then \ef{9.16.4} follows from \ef{9.10.1}, \ef{9.15.2} and \ef{9.16.2}. For \ef{a9}, notice that
\begin{align*}
h\sum_{n=1}^{N}\lt| \frac{v_{n }(t)}{x_{n }}-\frac{v_{n }(s)}{x_{n }}  \rt|^2
=h\sum_{n=1}^{N}\lt| \int_s^t \frac{\dot{v}_{n }(\tau)}{x_{n }}d\tau \rt|^2 \le h\sum_{n=1}^{N}  \int_s^t \lt|\frac{\dot{v}_{n }(\tau)}{x_{n }}\rt|^2 d\tau |t-s| .
\end{align*}
Then the estimate  for $v_n/x_n$ in \ef{a9} follows from \ef{9-7-2}. With this, the estimate  for $v_n$ in \ef{a9} holds obviously. Similarly, \ef{9.16.5} follows from \ef{9.15.2}; \ef{9.16.6} from \ef{9-7-2} and \ef{9.15.2}.

\hfill$\Box$

For $h={1}/{N}$, we define the functions $(r^h, v^h)(x, t )$  as follows:
\begin{align*}\label{a9}
r^h(x,t)=r_{n-1}(t)+\frac{r_n(t)-r_{n-1}(t)}{h}(x-x_{n-1}),    \\
 v^h(x,t)=v_{n-1}(t)+\frac{v_n(t)-v_{n-1}(t)}{h}(x-x_{n-1}),
\end{align*}
for $x_{n-1}<x<x_n$, $1\le n\le N$  and $0\le t\le T_*$. Then we have
\bee
r^h_t(x,t)=v^h(x,t), \ \
r^h_x(x,t)=\frac{r_n(t)-r_{n-1}(t)}{h} \ \ {\rm and} \ \  v^h_x(x,t)=\frac{v_n(t)-v_{n-1}(t)}{h}.
\eee

\begin{prop}\label{prop9.16} If the assumptions  in Lemma \ref{lem9.15} are satisfied,   then there exist  subsequences  of $\{r^h\}$,  $\{v^h\}$, $\{r^h/x\}$,  $\{v^h/x\}$, $\{\bar\rho^{\ga-1/2} r^h_x\}$ and $\{v^h_x /r^h_x\}$, still labeled by $\{r^h\}$,  $\{v^h\}$, $\{r^h/x\}$,  $\{v^h/x\}$, $\{\bar\rho^{\ga-1/2} r^h_x\}$ and $\{v^h_x /r^h_x\}$ for convenience, such that $\{r^h\}$,  $\{v^h\}$, $\{r^h/x\}$,  $\{v^h/x\}$, $\{\bar\rho^{\ga-1/2} r^h_x\}$ and $\{v^h_x /r^h_x\}$ converge boundedly and almost everywhere in $[0,1]\times[0,T_*]$ for $h\to 0$.
\end{prop}
{\em Proof}. We use similar arguments as in \cite{Okada, LiXY}.  First, we consider $\{r^h\}$ and $\{v^h\}$. It follows from  \ef{9.15.2}  and \ef{9.16.1}  that the functions of the families  $\{r^h\}$ and  $\{v^h\}$ as  functions of $x$ have uniformly bounded total variations with respect to $h$ for each fixed time $t\in [0, T_*]$. Let $t=t_k$ ($k=1, 2, \cdots$) be a countable set which is everywhere
dense in $[0, T_*]$. By Helly's theorem and a diagonal process, from the family of functions $\{v^h\}$, we can select a subsequence, still labeled as $\{v^h\}$ for convenience, converging boundedly and almost everywhere in $x\in I$ on the dense set $\{t_k; k=1, 2, \cdots\}$ in $[0, T^*]$ as $h\to 0$.
Consequently, by the Lebesgue's theorem, the subsequence $\{v^h\}$ converges in $L^2$-norm on $\{t_k; k=1, 2, \cdots\}$. Next, by \ef{a9},  the  continuity in time for the $L^2$-norm of $v^h$, it is standard to show that the
subsequence $\{v^h\}$ converges in $L^2(I)$ uniformly in $t\in [0, T_*]$, as $h\to 0$. So, we can select further a subsequence,  again still labeled as $\{v^h\}$ for convenience, converges almost everywhere in $(x, t)\in I\times [0, T_*]$. Denoting the limiting function of $\{v^h\}$ by $v$.  Since $r^h_t(x, t)=v^h(x, t)$, we know that the corresponding subsequence of $\{r^h\}$ (still labeled as $\{r^h\}$ for convenience) converges almost
everywhere to $r(x, t):=r_0(x)+\int_0^t v(x, s)ds$ in $I\times [0, T_*]$.

Similarly,  one can derive from estimates \ef{9.15.2}, \ef{9.16.2} and \ef{a9} that there exist certain subsequences of $\{r^h/x\}$ and $\{v^h/x\}$,  still labeled as $\{r^h/x\}$ and $\{v^h/x\}$ for convenience, converge  boundedly and almost everywhere  in $(x, t)\in I\times [0, T^*]$ to the functions $r(x,t)/x$ and $v(x,t)/x$, respectively. The estimates \ef{9.15.2}, \ef{9.16.3}, \ef{9.16.4}, \ef{9.16.5} and \ef{9.16.6}  guarantee that one can select certain subsequences of  $\{\bar\rho^{\ga-1/2} r^h_x\}$ and $\{v^h_x /r^h_x\}$ (still labeled as $\{\bar\rho^{\ga-1/2} r^h_x\}$ and $\{v^h_x /r^h_x\}$ for convenience) such that they converge  boundedly and almost everywhere  in $(x, t)\in I\times [0, T_*]$ to the functions $\bar\rho^{\ga-1/2}r_x(x,t)$ and $v_x(x,t)/r_x(x,t)$, respectively.

\hfill $\Box$

Due to the above convergence and uniform estimates on the approximating sequence, one may first verify that the  limiting function $(r, v)$  is a weak solution to  problem \ef{419} in the sense of
\begin{equation*}\begin{split}
&\int_0^{T_*}\int_{I} \lt\{ \bar\rho\left( \frac{x}{r}\right)^2  v\psi_t + \left(\frac{x^2}{r^2}\frac{\bar\rho}{ r_x}\right)^\ga   \psi_x  -  \frac{x^2}{r^4}    \bar\rho \int_0^x 4\pi y^2\bar\rho(y)   dy \rt\} (x, t) dxdt     \\
& =\lt.\int_{I} \lt(\bar\rho\left( \frac{x}{r}\right)^2  v\psi \rt) (x, t)dx\rt|_{t=0}^{T^*}+\int_0^{T^*}\int_{I} \lt\{ \mu\lt(\frac{v_x}{r_x}+2\frac{v}{r}\rt) \psi_x  + 2 \bar\rho\left( \frac{x}{r}\right)^3 \frac{v}{x}  v\psi \rt\}(x, t)dxdt   \end{split}\end{equation*}
for any test function $\psi\in C^1([0, T_*]\times I)$ satisfying $\psi(\cdot, t)\in C_c(I)$.
 Then one can use the standard regularity argument (cf. \cite{LXY}) which is ensured by the above uniform estimates to show it is   also a strong solution
as in Definition \ref{definitionss} satisfying
 $$\mathfrak{E}(t)\in C([0, T_*]) \  \  {\rm and} \ \
 \int_{I} \frac{1}{\bar\rho(x)} \mathfrak{B}_x^2(x, t) dx\le C \mathfrak{E}(0),    \ \ t\in [0, T_*].$$
This in particular implies that $\mathfrak{B}(\cdot, t)\in H^1(I)$ for  $t\in [0, T_*]$, so that one can define the trace
of $\mathfrak{B}$ at $x=1$.

The uniqueness of the strong solution follows from the argument in \cite{LXZ} (Section 11)  of the weighted energy estimates.

\subsection*{Part II. Linearized analysis}
We use the simple case that $\lambda_2=(2/3)\lambda_1=1/3$ to illustrate the main ideas of the linear analysis.
As in \cite{17'}, we set  $w=r/x-1$.  Then
\be\label{9.19.3}
r=x+x w, \ \ r_x=1+ w +x w_x, \ \  v= x w_t, \ \  v_x= w_t+ xw_{tx}, \ \  v_t=x w_{tt}.
\ee
 The linearized problem for \ef{419}  around  $w=0$ (the equilibrium of \ef{419}) reads
\begin{subequations}\label{linearization}\begin{align}
&x\bar\rho w_{tt} - (3\ga-4)(\bar\rho^{\ga})_x w -\ga x^{-3} \lt(\bar\rho^{\ga} x^4 w_x\rt)_x = (xw_{xt}+3w_t)_x,\label{7-2-3}\\
&     \mathfrak{B}_L:=xw_{xt}+w_t=0  \  \ {\rm at} \ \ x=1.
\end{align}\end{subequations}
The condition $v(0, t)=0$ has be incorporated in the transformation from $r$ to $w$ since $v=xw_t$.
(Indeed, \ef{linearization} follows from the facts  that the principal parts of $(r/x)^{-2\ga}$, $r_x^{-\ga}$, $(r/x)^{-4}$, $v_x/r_x$ and $v/r$ are $1 - 2\ga w$, $1 -\ga (w+xw_x)$, $1-4 w$, $xw_{xt}+w_t$ and $w_t$, respectively. The derivation of the left-hand side of \ef{7-2-3} can be also found in \cite{17'}.)
Naturally, we may rewrite \ef{7-2-3} as
\be\label{9.19.1}
x^4\bar\rho w_{tt} + (3\ga-4) \phi x^4 \bar\rho  w -\ga  \lt(  x^4  \bar\rho^{\ga}w_x\rt)_x = x^3 \lt(\mathfrak{B}_L\rt)_x  +2 x^3  w_{tx},
\ee
where $\phi$ is defined in \ef{rhox}, which is bounded from above and below by positive constants,  so that it is easy to see that $\ga>4/3$ is crucial to the linearized stability.

{\em Lower-order estimates}. The basic multipliers for \ef{9.19.1} are $  w_t$ and $  w$ which yield   the boundedness of
\be\label{9.19.4}
 \left\|(r-x, xr_x-x)(\cdot,t)  \right\|^2   +  (1+t )\left\| \left(v, xv_x \right)(\cdot, t) \right\|^2
  + (1+t) \lt\|  x \bar\rho^{1/2} v_t(\cdot, t) \rt\|^2.
\ee
Indeed, we have the following basic estimates for the solution of \ef{linearization} if $\ga>4/3$.
\bee\begin{split}
&\frac{1}{2}\frac{d}{dt} \int x^4 \lt( \bar\rho w_t^2+ (3\ga-4)\phi   \bar\rho w^2 + \ga  \bar\rho^\ga w_x^2  \rt)dx + \int x^2\lt[ (w_t+xw_{tx})^2+ 2   w_t^2\rt]dx =0,  \\
& \frac{1}{2}\frac{d}{dt} \int \lt[x^2(w+xw_{x})^2+ 2 x^2 w^2 + x^4 \bar\rho w w_t \rt]dx
+ \int \lt[ ( 3\ga-4)\phi x^4 \bar\rho w^2 + \ga x^4 \bar\rho^\ga w_x^2  \rt]dx \notag \\
&\qquad =\int x^4 \bar\rho w_t^2 dx,
\end{split}\eee
which implies
\be\label{basicestimate}\begin{split}
& \int \lt[x^2(w+xw_{x})^2+   x^2 w^2  \rt](x,t)dx + (1+t) \int \lt(x^4 \bar\rho w_t^2+    x^4 \bar\rho w^2 +  x^4 \bar\rho^\ga w_x^2  \rt)(x,t)dx\\
& +\int_0^t \int  \lt(  x^4 \bar\rho w^2 +   x^4 \bar\rho^\ga w_x^2  \rt)dx ds
+ \int_0^t  (1+s) \int \lt[x^2 (w_s+xw_{sx})^2+   x^2 w_s^2\rt]  dx ds \\
 \le &  C \int \lt[x^4 w_{x} ^2+   x^2 w^2 + x^4 \bar\rho w_t^2 \rt](x,0)dx .
 \end{split}\ee
Since the coefficients for the equation and boundary conditions in \ef{linearization} are independent of $t$, differentiating \ef{linearization} with respect to $t$, one obtains the same estimates for the corresponding $t$-derivatives of each quantity appearing
in \ef{basicestimate}, in particular,
\begin{align}
&  (1+t) \int x^4 \lt( \bar\rho w_{tt}^2+   \bar\rho w_t^2 +   \bar\rho^\ga w_{tx}^2  \rt)(x,t)dx
+ \int_0^t  (1+s) \int x^2\lt[ (w_{ss}+xw_{ssx})^2+    w_{ss}^2\rt]  dx ds \notag\\
& \le   C \int \lt[x^4 w_{x}^2+   x^2 w^2 + x^4 \bar\rho w_t^2  +x^4 \bar\rho w_{tt}^2  +  x^4 \bar\rho^\ga w_{tx}^2 \rt](x,0)dx . \label{9.19.2}
 \end{align}
Integrating the identity
\begin{align*}\frac{d}{dt}\lt[ (1+t) \int  x^2 \lt[(w_t+xw_{tx})^2+   w_t^2\rt] (x, t)dx \rt] =\int x^2 \lt[(w_t+xw_{tx})^2+    w_t^2\rt] (x, t)dx \\
+ 2 (1+t) \int  x^2 \lt[(w_t+xw_{tx})(w_{tt}+xw_{ttx}) +    w_t w_{tt}\rt] (x, t)dx
\end{align*}
with respect to $t$,  and use \ef{basicestimate} and \ef{9.19.2} to get
\be\label{basicestimate1'}\begin{split}
& (1+t) \int \lt[x^2 (w_t+xw_{tx})^2+   x^2 w_t^2\rt] (x, t)dx \\
\le & C \int \lt[x^4 w_{x}^2+   x^2 w^2 + x^2   w_t^2  +x^4 \bar\rho w_{tt}^2  +  x^4   w_{tx}^2 \rt](x,0)dx.\end{split}\ee
Therefore, the boundedness of  \ef{9.19.4}  is a consequence of \ef{basicestimate}-\ef{basicestimate1'} and \ef{9.19.3}.

{\em Higher-order estimates}. We may rewrite \ef{7-2-3} as
\begin{equation}\label{linearization2}\begin{split}
& G_{xt}+\ga\bar\rho^{\ga}G_x
 =x\bar\rho w_{tt} + \ga \phi x^2 \bar\rho  w_x+(3\ga-4)\phi x \bar\rho w,  \ \  {\rm where} \ \  G=xw_x+3w.
\end{split}\end{equation}
Based on \ef{basicestimate}, \ef{9.19.2} and \ef{hardyorigin}, we may apply a  multiplier $\bar\rho^{a}G_x$ with $a\ge 2\ga-2$ to \ef{linearization2}  to get
\be\label{basicestimate3}\begin{split}
&\int  \bar\rho^{a}G_x^2(x, t)dx+\int_0^t \int \bar\rho^{\ga+a}G_{x}^2 dxds \\
\le &  C \int \lt[ \bar\rho^{a}G_x^2 + x^4 w_{x}^2+   x^2 w^2 + x^4 \bar\rho w_t^2  +x^4 \bar\rho w_{tt}^2  +  x^4 \bar\rho^\ga w_{tx}^2 \rt](x,0)dx.\end{split}\ee
Due to equation \ef{linearization2}, we have
\bee\label{9.19}
\int x^4 \bar\rho w_{tt}^2  (x,0)dx \le  C \int  x^2 \bar\rho^{-1} \lt( |G_{xt}|^2 + \lt|\bar\rho^{\ga}G_x\rt|^2
+| x^2 \bar\rho  w_x |^2+ | x \bar\rho w |^2 \rt) (x, 0) dx,
\eee
which implies that the largest $a$ in \ef{basicestimate3} could be $2\ga-1$. So, we choose $a=2\ga-1$.  Moreover, it follows from \ef{basicestimate3} and  equation \ef{linearization2} that
$$
  \int_0^t \int \bar\rho^{ 3\ga -1 }G_{xt}^2 dxds \le C \int \lt[ \bar\rho^{2\ga-1}G_x^2 + x^4 w_{x}^2+   x^2 w^2 + x^4 \bar\rho w_t^2  +x^4 \bar\rho w_{tt}^2  +  x^4 \bar\rho^\ga w_{tx}^2 \rt](x,0)dx  .$$
This may suggest
\be\label{9.19-1}\begin{split}
&\int  \bar\rho^{2\ga-1}(x^2w^2_{xx}+w_x^2)(x, t)dx+\int_0^{t} \int \bar\rho^{3 \ga -1 }(x^2w_{sxx}+w_{sx}^2) dxds \\
\le &  C \int \lt[ \bar\rho^{2\ga-1}G_x^2 + x^4 w_{x}^2+   x^2 w^2 + x^4 \bar\rho w_t^2  +x^4 \bar\rho w_{tt}^2  +  x^4 \bar\rho^\ga w_{tx}^2 \rt](x,0)dx  .\end{split}\ee
(Indeed, it is one of main ideas to justify this for the corresponding nonlinear equation). The higher-order estimate \ef{9.19-1} can improve the regularity near the origin as follows. It follows from \ef{9.19-1},  \ef{basicestimate} and \eqref{hardyorigin} that
\be\label{9.19-2}
\int_0^{t} \int (w_s^2 + w_{xs}^2) dx ds \le  C \int \lt[ \bar\rho^{2\ga-1}G_x^2 +   x^2 w^2 + x^4 \lt( w_{x}^2 +   \bar\rho w_t^2  +  \bar\rho w_{tt}^2  +   \bar\rho^\ga w_{tx}^2 \rt)\rt](x,0)dx.
\ee
Let $\psi$ be a  non-increasing cut-off function defined on $[0, 1]$ satisfying
$\psi=1$  on $ [0,1/4]$, $ \psi=0$  on $ [{1}/{2}, 1] $, and $|\psi'|\le 32$ .
With \ef{9.19-2} and \ef{9.19.2}, we can integrate the product of $ \ef{7-2-3}_t$ and  $ \psi x w_{tt}$ with respect to $x$ and $t$ to get
$$\int_0^{1/4} x^2  w_{tt}^2 (x, t)dx \le  C \int \lt[ \bar\rho^{2\ga-1}G_x^2 + x^4 w_{x}^2+   x^2 w^2 + x^4 \bar\rho w_t^2  +x^2 \bar\rho w_{tt}^2  +  x^4 \bar\rho^\ga w_{tx}^2 \rt](x,0)dx,$$
which, together with \ef{9.19.2}, gives
\be\label{9.19}
\int  x^2 \bar\rho w_{tt}^2 (x, t)dx \le  C \int \lt[ \bar\rho^{2\ga-1}G_x^2 + x^4 w_{x}^2+   x^2 w^2 + x^4 \bar\rho w_t^2  +x^2 \bar\rho w_{tt}^2  +  x^4 \bar\rho^\ga w_{tx}^2 \rt](x,0)dx.
\ee
As a consequence of \ef{9.19.3}, \ef{9.19-1} and \ef{9.19},    we can obtain  the boundedness of
\be\label{9.19.4}
 \lt\|  \bar\rho^{1/2} v_t(\cdot, t) \rt\|^2 +  \lt\|  \bar\rho^{\ga-1/2}\lt(r_{xx}, \  (r/x)_x\rt)(\cdot, t) \rt\|^2 .
\ee

{\em Conclusion}:  A natural  higher-order   functional for the study of the linear problem \ef{linearization} is:
$$E(t)=\left\|(r-x, xr_x-x)(\cdot,t)  \right\|^2   +  \left\| \left(v, xv_x \right)(\cdot, t) \right\|^2
  +   \lt\|   \bar\rho^{1/2} v_t(\cdot, t) \rt\|^2 +  \lt\|  \bar\rho^{\ga-1/2}\lt(r_{xx}, \  (r/x)_x\rt)(\cdot, t) \rt\|^2.$$

\begin{small}
\noindent  
T. Luo \\
{\it Department of Mathematics and Statistics,
Georgetown University,
Washington, DC, 20057, USA} \\
{E-mail: tl48@georgetown.edu} \\
Z. Xin\\
{\it Institute of Mathematical Sciences,
The Chinese University of Hong Kong, Hong Kong}\\
{E-mail: zpxin@ims.cuhk.edu.hk} \\
H. Zeng
{\it Yau Mathematical Sciences Center,
Tsinghua University,
Beijing, 100084, China$^1$}\\
{E-mail: hhzeng@mail.tsinghua.edu.cn} \\
 {\it Center of Mathematical Sciences and Applications,  Harvard University,   Cambridge, MA 02318, USA} \\
\end{small}


\begin{thebibliography}{100}
\bibitem{ch} S. Chandrasekhar,  An Introduction to the Study of Stellar Structures. University of
Chicago Press, Chicago, 1938.

 \bibitem{Chengq} Chen, Gui-Qiang; Kratka, Milan Global solutions to the Navier-Stokes equations for compressible heat-conducting flow with symmetry and free boundary. Comm. Partial Differential Equations 27 (2002), no. 5-6, 907--943.

\bibitem{10} Coutand, D., Shkoller, S.: Well-posedness in smooth function spaces for the moving-
boundary 1-D compressible Euler equations in physical vacuum. Commun. Pure
Appl. Math. 64, 328-366 (2011)

\bibitem{10'}
Coutand, Daniel; Shkoller, Steve; Well-Posedness in Smooth Function Spaces for the Moving-Boundary Three-Dimensional Compressible Euler Equations in Physical Vacuum. Arch. Ration. Mech. Anal. 206 (2012), no. 2, 515-616.

\bibitem{DLYY}
Y. Deng, T.P. Liu, T. Yang,  Z. Yao,  Solutions of Euler-Poisson equations for gaseous
stars. Arch. Ration. Mech. Anal. 164 (2002), no. 3, 261-285.

\bibitem{DZ}
B. Ducomet, A. Zlotnik: Stabilization and stability for the spherically symmetric Navier-Stokes-Poisson system, Appl.Math.Lett. 18 (2005), 1190-1198

\bibitem{duan}
Q. Duan, On the dynamics of Navier-Stokes equations for a shallow water model, J. Differential Equations, 250 (2011), 2687-2714.


\bibitem{KM}
A. Kufner, L. Maligranda, L.-E. Persson: The Hardy inequality. Vydavatelsksy
Servis, PlzeTn, 2007. About its history and some related results.
\bibitem{fangzhang}D.-Y. Fang and T. Zhang, Global behavior of compressible Navier-Stokes equations with a degenerate viscosity coefficient, Arch. Rational Mech. Anal., 182 (2006), 223-253.


\bibitem{fangzhang1}
D.-Y. Fang and T. Zhang, Global behavior of spherically symmetric Navier-Stokes-Poisson system with degenerate viscosity coefficients, Arch. Rational Mech. Anal., 191 (2009), 195-243.

\bibitem{GLX}
Guo, Zhenhua; Li, Hai-Liang; Xin, Zhouping, Lagrange structure and dynamics for solutions to the spherically symmetric compressible Navier-Stokes equations. Comm. Math. Phys. 309 (2012), no. 2, 371-412



\bibitem{jangnsp}
Jang, J., Local well-posedness of dynamics of viscous gaseous stars. Arch. Rational
Mech. Anal. 195, 797-863 (2010).

\bibitem{jang65}
J. Jang. Nonlinear instability in gravitational Euler-Poisson system for $\gamma=\frac{6}{5}$. Arch.
Ration. Mech. Anal. 188 (2008), no. 2, 265--307.

\bibitem{jangmas}
Jang, J., Masmoudi, N.:Well-posedness for compressible Euler with physical vacuum singularity. Commun. Pure Appl. Math. 62, 1327-1385 (2009)

\bibitem{jm}
Jang, J., Masmoudi, N.: Well-posedness of compressible Euler equations in a physical vacuum, Comm. Pure Appl. Math. 68 (2015),  61--111.


\bibitem{jangtice}
J. Jang, I. Tice, Instability theory of the Navier-Stokes-Poisson equations, Anal. PDE 6 (2013), no. 5, 1121-1181.

\bibitem{17'}
Jang, J. Nonlinear Instability Theory of Lane-Emden stars, Comm. Pure Appl. Math. 67 (2014),   1418--1465.

\bibitem{JXZ}
S. Jiang, Z.  Xin, P. Zhang,
Global weak solutions to 1D compressible isentropic Navier-Stokes equations with density-dependent viscosity.
Methods Appl. Anal.  12  (2005),  no. 3, 239-251.

\bibitem{lebovitz1}
Lebovitz, N.R., Lifschitz, A.: Short-wavelength instabilities of Riemann ellipsoids, Philos. Trans. Roy. Soc. London Ser. A 354(1709), 927-950 (1996).

\bibitem{lebovitz2}
Lebovitz, N.R.: The virial tensor and its application to self-gravitating fluids. Astrophys. J. 134, (1961) 500-536.

\bibitem{liebyau}
Lieb, E.H., Yau, H.T.: The Chandrasekhar theory of stellar collapse as the limit of quantum mechanics. Commun. Math. Phys. 112(1), 147-174 (1987)

\bibitem{linss} S.-S. Lin. Stability of gaseous stars in spherically symmetric motions. SIAM J. Math. Anal. 28 (1997), no. 3, 539-569.

\bibitem{tpliudamping}
Liu, T.-P. Compressible flow with damping and vacuum. Japan J. Appl. Math. 13 (1996), 25-32.


\bibitem{13} Liu, T.-P.; Yang, T.,  Compressible flow with vacuum and physical singularity. Methods Appl. Anal. 7 (2000), 495-509.

\bibitem{LiXY}
T.-P. Liu, Z. Xin and T. Yang, Vacuum states of compressible flow, Discrete and Continuous Dynamical Systems, 4(1998), 1-32.

\bibitem{LXY} T. Luo, Z. Xin and T. Yang,  Interface behavior of compressible Navier-Stokes equations with vacuum, SIAM J. Math. Anal. 31 (6) (2000) 1175-1191.

\bibitem{luosmoller1} Luo, T.; Smoller, J.: Nonlinear dynamical stability of Newtonian rotating and non-rotating white dwarfs and rotating supermassive stars. Comm. Math. Phys. 284 (2008), no. 2, 425-457.

\bibitem{luosmoller2} Luo, T.; Smoller, J.: Existence and non-linear stability of rotating star solutions of the compressible Euler-Poisson equations. Arch. Ration. Mech. Anal. 191 (2009), no. 3, 447-496.

\bibitem{LXZ} T. Luo, Z. Xin and H. Zeng, Well-Posedness  for the Motion  of  Physical Vacuum of  the Three-dimensional Compressible Euler Equations with or without Self-Gravitation,  Arch. Ration. Mech. Anal. 213 (2014), 763-831.

\bibitem{LXZ2} T. Luo, Z. Xin and H. Zeng,  	Nonlinear asymptotic stability of the Lane-Emden solutions for the viscous gaseous star problem with degenerate density dependent viscosities, arXiv:1507.01069.

\bibitem{makino} Makino, T.: On a local existence theorem for the evolution equation of gaseous stars.
Patterns and Waves. Stud. Math. Appl., Vol. 18. North-Holland, Amsterdam, 459-479,
1986



\bibitem{Okada} M. Okada, Free boundary value problems for the equation of one-dimensional motion of viscous gas, Japan J. Appl. Math., 6(1989), 161-177.




\bibitem{Okada1} M. Okada and T. Makino, Free boundary problem for the equations of spherically symmetrical motion of viscous gas, Japan J. Indust. Appl. Math., 10(1993), 219-235.

\bibitem{OSM} S. Matusu-Necasova, M. Okada, T. Makino: Free boundary problem for the equation of
spherically symmetric motion of viscous gas III, Japan J.Indust.Appl.Math. 14 (1997), 199-213

\bibitem{Okada3} Mari Okada, Free boundary problem for one-dimensional motions of compressible gas and vacuum, Japan J. Indust. Appl. Math. 21 (2) (2004) 109-128.


 \bibitem{rein}G. Rein, Non-linear stability of gaseous stars. Arch. Ration. Mech. Anal. 168 (2003), no. 2, 115-130.

\bibitem{94}  P. Secchi, On the evolution equations of viscous gaseous stars. Ann. Scuola Norm. Sup. Pisa Cl. Sci. (4) 18 (1991), no. 2, 295--318.

\bibitem{95} P. Secchi, On the uniqueness of motion of viscous gaseous stars. Math. Methods Appl. Sci. 13 (1990), no. 5, 391--404.

 \bibitem{tokusky} S. H. Shapiro \& S. A. Teukolsky,  Black Holes, White Dwarfs, and Neutron Stars, {\it WILEY-VCH}, (2004)

\bibitem{96}  Strohmer, G.,  Asymptotic estimates for a perturbation of the linearization of an equation for compressible viscous fluid flow. Studia Math. 185 (2008), no. 2, 99--125.
 \bibitem{weinberg} S. Weinberg , Gravitation and Cosmology {\it John Wiley
and Sons}, New York , 1972.


    \bibitem{ya} T. Yang, Singular behavior of vacuum states for compressible fluids. J. Comput. Appl. Math.  190  (2006),  no. 1-2, 211-231.

\bibitem{YYZ}Yang, Tong; Yao, Zheng-an; Zhu, Changjiang,
Compressible Navier-Stokes equations with density-dependent viscosity and vacuum. (English summary)
Comm. Partial Differential Equations 26 (2001), no. 5-6, 965-981.

\bibitem{yangzhu}Yang, Tong; Zhu, Changjiang,
Compressible Navier-Stokes equations with degenerate viscosity coefficient and vacuum. Comm. Math. Phys. 230 (2002), no. 2, 329-363.


\bibitem{zhu}Zhu, Changjiang; Zi, Ruizhao, Asymptotic behavior of solutions to 1D compressible Navier-Stokes equations with gravity and vacuum. Discrete Contin. Dyn. Syst. 30 (2011), no. 4, 1263-1283



\end{thebibliography}
\end{document}